\documentclass[nonblindrev]{informs3}

\OneAndAHalfSpacedXI 

\MANUSCRIPTNO{}

\usepackage{natbib}
\bibpunct[, ]{(}{)}{,}{a}{}{,}%
\def\bibfont{\small}%
\def\BIBand{and}%
 
\TheoremsNumberedThrough     
\ECRepeatTheorems

\EquationsNumberedThrough    

\usepackage[caption=false]{subfig}
\usepackage{amsmath, amsfonts, graphicx, epsfig, epstopdf, xcolor, mathrsfs, theoremref, enumerate, setspace, algorithm, bbm, amssymb, mathrsfs, makecell}
\usepackage[shortlabels]{enumitem}
\usepackage[margin=1in]{geometry}
\usepackage[utf8]{inputenc}
\usepackage[noend]{algpseudocode}
\usepackage{mathtools}
\usepackage{tikz}
\usetikzlibrary{shapes}
\usepackage{booktabs}
\usepackage{numprint} 
\usepackage{diagbox}
\usepackage{hyperref}
\usepackage{chngcntr}

\newcommand{\cD}{\mathcal{D}}

\newcommand{\cS}{\mathcal{S}}

\newcommand{\cK}{\mathcal{K}}

\newcommand{\forAll}{\, \forall \,}
\newcommand{\setbar}{\text{ }|\text{ }}
\newcommand{\myProof}{\noindent\textbf{Proof: }}

\newcommand{\drug}{U}
\newcommand{\conc}{C}

\newcommand{\econ}{E}
\newcommand{\pop}{P}
\newcommand{\epop}{N}
\newcommand{\bin}{Z}
\newcommand{\mC}{B}
\newcommand{\mv}{V}

\newcommand{\tym}{s}
\newcommand{\lastTym}{S}

\newcommand{\tymEnd}{T}
\newcommand{\timeSet}{\cS}
\newcommand{\dSet}{\cD}
\newcommand{\operSet}{\mathscr{O}}
\newcommand{\doperSet}{\widehat{\operSet}}
\newcommand{\typeSet}{\mathcal{Q}}

\newcommand{\DAYS}{\mathrm{DAYS}}

\newcommand{\MEALS}{\mathrm{MEALS}}

\newcommand{\rhs}{\beta}
\newcommand{\para}{\alpha}
\newcommand{\bpara}{\boldsymbol{\para}}
\newcommand{\concPar}{\xi}
\newcommand{\gompOne}{\Lambda}

\newcommand{\dCons}{\eta}
\newcommand{\dExp}{\rho}

\newcommand{\finishTym}{T}

\newcommand{\ipop}{\pi}
\newcommand{\bipop}{\boldsymbol{\ipop}}

\newcommand{\vol}{\mathscr{V}}
\newcommand{\Prob}{\mathbf{Pr}}

\newcommand*{\myQED}{\hfill\ensuremath{\square}}

\makeatletter
\newcommand*\bigcdot{\mathpalette\bigcdot@{.5}}
\newcommand*\bigcdot@[2]{\mathbin{\vcenter{\hbox{\scalebox{#2}{$\m@th#1\bullet$}}}}}
\makeatother
\begin{document}

\RUNAUTHOR{Ajayi, Hosseinian, Schaefer, Fuller}

\RUNTITLE{Combination Chemotherapy Optimization with Discrete Dosing}

\TITLE{Combination Chemotherapy Optimization with Discrete Dosing}

\ARTICLEAUTHORS{%
\AUTHOR{Temitayo Ajayi}
\AFF{Nature Source Improved Plants, Ithaca, NY, 14850, \EMAIL{tayo.ajayi25@gmail.com}} 
\AUTHOR{Seyedmohammadhossein Hosseinian}
\AFF{Department of Computational and Applied Mathematics, Rice University, Houston, TX, 77005, \EMAIL{hosseinian@rice.edu}} 
\AUTHOR{Andrew J. Schaefer}
\AFF{Department of Computational and Applied Mathematics, Rice University, Houston, TX, 77005, \EMAIL{andrew.schaefer@rice.edu}} 
\AUTHOR{Clifton D. Fuller}
\AFF{Department of Radiation Oncology, The University of Texas MD Anderson Cancer Center, Houston, TX, 77030, \EMAIL{cdfuller@mdanderson.org}
}
} 
\ABSTRACT{%
Chemotherapy is one of the primary modalities of cancer treatment. Chemotherapy drug administration is a complex problem that often requires expensive clinical trials to evaluate potential regimens. One way to alleviate this burden and better inform future trials is to build reliable models for drug administration. 
Previous chemotherapy optimization models have mainly relied on optimal control, which does not lend itself to capturing complex and vital operational constraints in chemotherapy planning involving discrete decisions, such as doses via pills and rest periods. In addition, most of the existing models for chemotherapy optimization lack an explicit toxicity measure and impose toxicity constraints primarily through (fixed) limits on drug concentration. The existing stochastic optimization models also focus on maximizing the probability of cure when tumor heterogeneity is uncertain. In this paper, we develop a mixed-integer program for combination chemotherapy (utilization of multiple drugs) optimization that incorporates various important operational constraints and, besides dose and concentration limits, controls treatment toxicity based on its effect on the count of white blood cells. To address the uncertainty of tumor heterogeneity, we propose chance constraints that guarantee reaching an operable tumor size with a high probability in a neoadjuvant setting. We present analytical results pertinent to the accuracy of the model in representing biological processes of chemotherapy and establish its merit for clinical applications through a numerical study of breast cancer. 
}



\maketitle
%

\section{Introduction}
\label{sec:intro}
Chemotherapy, i.e., the administration of cytotoxic drugs, is a prominent cancer treatment modality. In contrast to local treatment methods, such as surgery and radiation therapy, chemotherapy is a systemic treatment that targets cancer cells throughout a patient's body. Hence, it is widely used for patients in advanced stages of cancer; more than 60\% of the patients diagnosed with stage III or IV of breast, colon, rectal, lung, testicular, urinary bladder, and uterine corpus cancers in the United States underwent chemotherapy in 2016~\citep{CancerSurvivorship2021}. Cytotoxic drugs kill cancer cells, which results in tumor shrinkage. However, due to their toxic nature and narrow therapeutic margin, these drugs damage healthy cells as well and come with several, possibly life-threatening side effects~\citep{sideEffects}. The main objective of chemotherapy planning is to determine administration dosage and schedule for cytotoxic drugs such that a significant tumor shrinkage is achieved, or ideally, disappear, while the adverse effects on healthy organs are minimized~\citep{howChemo2021}. Chemotherapy treatment plans are typically evaluated by randomized clinical trials; see e.g.,~\citep{ebata2018randomized, zhao2020efficacy, mariotti2021effect}. 
Such trials are limited in variability due to constraints of permissible treatments and clinical and ethical considerations. Mathematical models of chemotherapy decision-making can alleviate some of these burdens while aiding treatment improvement and evaluation.

Mathematical models for chemotherapy planning must account for the dynamics of tumor evolution as well as the pharmacokinetics (distribution within the body) and pharmacodynamics (effect on tumor and healthy cells) of cytotoxic drugs. These processes take place in continuous time and are naturally described by ordinary differential equations (ODEs). In this regard, chemotherapy planning has been mainly approached as an optimal control (OC) problem in the literature.~\cite{Swan1977} first studied chemotherapy planning as an OC problem, followed by the seminal models proposed by~\cite{martin1990} and~\cite{martin1992}; the objective of these models is to minimize cancer cell population at the end of a treatment period subject to drug concentration limits---as a measure of toxicity---and intermediate tumor size or shrinkage rate. The early OC chemotherapy literature contained minimal details, which allowed many of them to be solved analytically; these include~\citep{Swan1977, zietz1979, Murray1990, Murray1994, Panetta1995, murray1997}. Extending~\cite{martin1990} and~\cite{martin1992}, more complex and realistic OC models for chemotherapy optimization used approximation techniques~\citep{Martin1992Low, Martin1992Resistant, pereira1995, costa1997, nanda2007, dePillis2007, dOnofrio2009OnTumors, itik2009, Harrold2009} and heuristic algorithms~\citep{Iliadis2000, Tan2002, Floares2003, villasana2004, Liang2006, tse2007, alam2013}. 
We refer to~\cite{Shi2014} and~\cite{saville2019} for detailed surveys.

Continuous OC models capture the ``biological'' dynamics of chemotherapy processes well; however, cancer treatment involves important discrete components and operational constraints. For example, some cytotoxic drugs are available in the form of pills and are taken orally. For these drugs, an administration dose must be a multiple of the pill size, e.g., one pill per day, and any deviation from this regimen can easily lead to an underdose or overdose. Oral drugs are often prescribed to be taken with food to help with their digestion and to alleviate their side effects. Metabolic processes can lead to mandated rest periods for certain drugs. These give a discrete nature to drug administration scheduling, which is not captured by continuous OC models. Modeling such operational constraints for chemotherapy planning requires integer control variables; introducing integer variables to OC problems makes them extremely hard to solve~\citep{sager2005numerical}. In addition, the existing chemotherapy optimization models mainly impose treatment toxicity constraints implicitly, through (fixed) limits on drug concentration. This presents another challenge to the applicability of these models to combination chemotherapy, i.e., utilization of multiple cytotoxic drugs, which is the common practice in the presence of drug resistance~\citep{Luqmani2005,hu2016recent}. In fact, in the absence of an explicit toxicity measure, these models do not clarify how the adverse effects of chemotherapy could vary under different combinations of administration regimens for multiple drugs.

Tumor heterogeneity is another important consideration in cancer treatment planning~\citep{Polyak2011,hu2017}. 
Tumors are composed of different cell types with distinct characteristics; tumor heterogeneity is considered one of the main factors of therapeutic resistance~\citep{y2020}. Recent advances in next-generation sequencing (NGS) technologies have made the characterization of the cell composition of a tumor possible; this requires multiple, spatially separated samples from the tumor~\citep{gerlinger2012,PIRAINO2019}. However, sampling from a tumor can lead to needle tract seeding, i.e., implantation of cancerous cells in healthy regions, which may lead to cancer metastasis; a higher risk of seeding is incurred as the number of sampling passes increases~\citep{tyagi2014}. In the absence of multiple biopsies, tumor heterogeneity remains uncertain for treatment planning~\citep{abecassis2019}. The only existing chemotherapy optimization models that consider this uncertainty include~\citep{Coldman1983, Day1986TreatmentChemotherapy, coldman2000}. These models suffer from the same shortcomings as the existing deterministic ones, and they seek to maximize the ``probability of cure'' over the course of a treatment, which challenges their clinical relevance even more than the deterministic models.

To fill the aforementioned gaps between the existing models of chemotherapy optimization and medical practice, we present a mixed-integer linear programming (MILP) model for combination chemotherapy planning, which seeks to find optimal administration dose and schedule for cytotoxic drugs by minimizing cancer cell population at the end of a treatment. We use discretization and linearization techniques to recast ODEs representing the biological and pharmacological processes of chemotherapy into a MILP framework; the flexibility of this framework allows for modeling complex operational constraints of chemotherapy. In particular, we incorporate discrete administration dose and schedule as well as clinically mandated rest periods in our model. We use the white blood cell count as a measure of treatment toxicity. More specifically, we consider the effect of cytotoxic drugs on the count of two major white blood cell types, namely neutrophil and lymphocyte, which account for more than 80\% of the total white blood cells in human body. Neutropenia (low count of neutrophils) is often a toxicity of concern in chemotherapy~\citep{Pizzo1993,Kosaka2015,kasi2018}, which can leave cancer patients susceptible to other, possibly fatal diseases during a chemotherapy treatment. To address the uncertainty of tumor heterogeneity, we propose chance constraints and present a neoadjuvant (prior to a primary surgery) chemotherapy optimization model for treatment planning under this uncertainty. We provide analytical results concerning the accuracy of the model in representing biological processes of chemotherapy. We use the clinical literature and published data for patients with breast cancer to calibrate our model parameters and perform sensitivity analysis to identify the most influential factors in a treatment outcome.

The structure of this paper is as follows: Section~\ref{sec:prelim} presents chemotherapy modeling preliminaries, including tumor and white blood cell population dynamics, pharmacokinetics and pharmacodynamics of cytotoxic drugs, and operational constraints. Section~\ref{sec:models} presents our deterministic and stochastic MILP models for combination chemotherapy planning along with our analytical results. Section~\ref{sec:calib} includes the model calibration details, and Section~\ref{sec:results} contains the results of our computational experiments. Section~\ref{sec:conc} concludes the paper. The proofs and some technical details are provided in the e-companion of the paper. 
\section{Modeling Preliminaries}
\label{sec:prelim}
Throughout this paper, we consider a treatment period $[0, T] \subset \mathbb{R}$ and a set of available cytotoxic drugs $\dSet$. For each drug $d \in \dSet$, the (continuous) functions $\drug_{d}(t)$ and $\conc_{d}(t)$ represent the administration dose and drug concentration, respectively, at time $t \in [0, T]$. We denote the set of cancer cell types by $\typeSet$, and for each cell type $q \in \typeSet$, we use $\epop_{q}(t)$ to represent the corresponding cell count as a function of time. We also introduce the variable functions $\pop_{q}(t) = \ln (\epop_{q}(t)),\forAll q \in \typeSet$. Each cancer cell type is resistant to a (possibly empty) subset of drugs. Note that drug resistance can be present even before chemotherapy starts~\citep{Swierniak2009}. Finally, the white blood cell count at time $t \in [0, T]$ is denoted by $N_w(t)$; we will distinguish between neutrophils and lymphocytes when we present the operational constraints of our models.
\subsection{Cell Population Dynamics}
\label{sec:cellDynamics}
Tumors proliferate by cell division and exhibit exponential growth in early stages, but the growth rate gradually decreases as malignant cells compete for limited nutritional resources. This resembles an S-curve growth, which is most commonly modeled by a Gompertzian function in cancer research; e.g.,~\citep{Laird1964, laird1965dynamics, Norton1988, Harrold2009, Frances2011, Tjorve2017}. The ODE representation of this function is as follows: 
\begin{align} \label{GompertzLaird}
\dot{\epop_{q}}(t) = \gompOne \, \epop_{q}(t) \, \ln \left ( \dfrac{\epop_{q,\infty}}{\epop_{q}(t)} \right ), \ \epop_{q}(0) = \epop_{q,0}, 
\end{align}
where $\epop_{q,0}$ and $\epop_{q,\infty}$ denote the initial population of cancer cell type $q \in \typeSet$ and its steady-state asymptotic limit, respectively, and $\gompOne$ is a shape parameter that dictates the rate at which the population transitions from the initial state to the steady-state limit.

The population dynamics of white blood cells are different. White blood cells are perpetually produced (mainly in bone marrow and the thymus gland) and circulate in the blood; they normally have a lifespan of a few days.~\cite{Iliadis2000} model the white blood cell dynamics as follows:
\begin{align} \label{wbcModel}
\dot{\epop_{w}}(t) = \upsilon_{w} - \nu_{w} \, \epop_{w}(t), \
\epop_{w}(0) = \epop_{w,0},
\end{align}
where $\upsilon_{w}$ and $\nu_{w}$ are the white blood cells' production and turnover rates, respectively, and $\epop_{w,0}$ denotes their (constant) level in the body under normal conditions. It is easy to verify that $\epop_{w}(t) = \epop_{w,0}$ is a solution to Eq.~\eqref{wbcModel} given $\upsilon_{w} = \nu_{w} \epop_{w,0}$.

Eqs.~\eqref{GompertzLaird} and~\eqref{wbcModel} provide the basis for our pharmacodynamics models. 
\subsection{Pharmacokinetics}
\label{sec:pharmacokinetics}
A drug's distribution within the body (pharmacokinetics) is known to be a complex, multi-compartmental, and multi-phase process. In cancer research, however, the dose profile of a cytotoxic drug is often represented by a single-compartmental model, where the drug concentration decays exponentially over time~\citep{martin1992, Martin1994, jacqmin2007, Harrold2009, Frances2011}. The process is described by the following ODE:
\begin{align} \label{doseHistoryODE}
\dot{\conc_{d}}(t) = - \concPar_{d} \, \conc_{d}(t) + \frac{\drug_{d}(t)}{\vol}, \ \conc_{d}(0) = 0,
\end{align}
where $\vol > 0$ represents the volume of the ``effect compartment'' that is used to convert an administered dose to drug concentration, and $\concPar_{d}$ is a constant characterizing the elimination rate of a drug $d \in \dSet$ in the body. In our models, the boundary condition $\conc_{d}(0) = 0$ indicates there is no drug in a patient's body before the start of treatment.

The underlying assumption of Eq.~\eqref{doseHistoryODE} is that the contribution of a newly administered dose of a drug to its concentration profile starts from the peak it generates on the concentration curve. After a single administration, the drug concentration-time curve is highly right-skewed; it reaches its peak in a relatively short time while it takes much longer for the drug to vanish from a patient's body. For the sake of simplicity, Eq.~\eqref{doseHistoryODE} ignores the short time it takes for a drug to reach its maximum concentration after an administration. 
\subsection{Pharmacodynamics}
\label{sec:Pharmacodynamics}
The main paradigm of pharmacodynamics (drug effect) modeling in chemotherapy optimization is based on the seminal works of~\cite{Skipper1964XIII, skipper1967XXI}, which indicate that, given a dose of a cytotoxic drug, it kills a constant fraction of cancer cells. The fractional kill effect of a cytotoxic drug on cancer cells is modeled by adding a bilinear term, composed of the product of drug concentration and cancer cell count with a constant factor, to the cancer evolution model, i.e., Eq.~\eqref{GompertzLaird}. In combination chemotherapy, the effect of multiple drugs is commonly modeled following the additivity principal; see e.g.,~\citep{Martin1992Low, Petrovski2004, tse2007, Frances2011}. In an additive model, the drugs perform as if each acts in isolation, and the effects of all drugs are summed. Our pharmacodynamics model follows~\cite{Frances2011}, who also assume an exponential decay on drug effectiveness over time, due to the resistance developed in cancer cells when exposed to a drug, as follows: 
\begin{align} \label{GompertzLaird+ph}
\dot{\epop_{q}}(t) = \gompOne \, \epop_{q}(t) \, \ln \left ( \dfrac{\epop_{q,\infty}}{\epop_{q}(t)} \right ) - \sum\limits_{d \in \dSet} \dCons_{d,q} \, \exp(-\dExp_{d,q} \, t) \, \epop_{q}(t) \, \econ_{d}(t), \ \epop_{q}(0) = \epop_{q,0},
\end{align}
where $\dCons_{d,q}$ is the fractional kill effect parameter of a drug $d \in \dSet$ on a cancer cell type $q \in \typeSet$, the parameter $\rho_{d,q}$ determines how drug effectiveness decays over time, and $\econ_{d}(t)$ denotes the {\em effective} concentration of a drug $d \in \dSet$ as a function of time. The effective concentration $\econ_{d}(t)$ indicates the amount that the drug concentration exceeds some threshold $\rhs_{d,\mathrm{eff}}$, below which the drug is ineffective therapeutically~\citep{Iliadis2000,Tan2002,Harrold2009}. By this definition,
\begin{align*} \label{eq:effConc}
\econ_{d}(t) = \max\{0,~\conc_{d}(t) - \rhs_{d,\mathrm{eff}}\}.
\end{align*}
Observe that, with the logarithmic transformation $\pop_{q}(t) = \ln (\epop_{q}(t)),\forAll q \in \typeSet$, Eq.~\eqref{GompertzLaird+ph} can be written as a linear equation, as follows:
\begin{equation*} \label{drugEffectAdditive}
\dot{\pop_{q}}(t) = \gompOne \, \Big ( \ln (\epop_{q,\infty}) - \pop_{q}(t) \Big ) - \sum\limits_{d \in \dSet} \dCons_{d,q} \, \exp(-\dExp_{d,q} \, t) \, \econ_{d}(t), \ \pop_{q}(0) = \ln (\epop_{q,0}).
\end{equation*}

We model the fractional kill effect of cytotoxic drugs on white blood cells in a similar manner:
\begin{align} \label{wbcModelDrugs}
\dot{\epop_{w}}(t) = \upsilon_{w} - \nu_{w} \, \epop_{w}(t) - \sum\limits_{d \in \dSet} \dCons_{d, w} \, \epop_{w}(t) \, C_{d}(t - t_{w}), \ t \geq t_{w}, 
\end{align}
where $\dCons_{d, w}$ is the fractional kill effect parameter of a drug $d \in \dSet$ on white blood cells, and $t_{w}$ denotes the delay in the response of white blood cells to cytotoxic drugs~\citep{Iliadis2000}. Note that, in Eq.~\eqref{wbcModelDrugs}, we make conservative assumptions that white blood cells do not develop resistance to cytotoxic drugs over time and the toxic effect of a drug exists even if its concentration is below the threshold $\rhs_{d,\mathrm{eff}}$. During the time interval $[0,t_{w})$, the white blood cell population dynamics is governed by Eq.~\eqref{wbcModel}.
\subsection{Operational Constraints}
\label{sec:operational}
Operational constraints enforce clinically permissible treatment plans. We consider a partition of the treatment period $[0, \tymEnd]$ into $M$ days, each denoted by $D_{m},~m \in \DAYS = \{1,\ldots,M\}$. Because certain oral drugs are consumed with meals, e.g., capecitabine \citep{Segal2014}, three time points (periodic with respect to the days) are designated as meal times; this set of time points is denoted by $\MEALS$.

Below, we describe operational constraints captured by our model; some of these constraints are explicitly included in the chemotherapy optimization literature.

\begin{enumerate}
\item \textbf{Maximum concentration} \citep{Martin1992Low,Iliadis2000,Baker2006}: For a drug $d \in \dSet$, let $\rhs_{d,\mathrm{conc}}$ denote the maximum permissible concentration in a patient's body; the corresponding constraint is $\conc_{d}(t) \leq \rhs_{d, \mathrm{conc}}, \forAll t \in [0,\tymEnd].$
\item \textbf{Maximum infusion rate} \citep{Hande1998,Reigner2001,Baker2006,Ershler2006,Palmeri2008}: Let $\rhs_{d,\mathrm{rate}}$ denote the maximum permissible infusion rate for a drug $d \in \dSet$; the corresponding constraint is $\drug_{d}(t) \leq \rhs_{d, \mathrm{rate}}, \forAll t \in [0,\tymEnd].$ 
\item \textbf{Maximum daily cumulative dose} \citep{Hande1998,Reigner2001,Baker2006,Ershler2006,Palmeri2008}: Clinical studies often seek to determine appropriate thresholds for drug administration within particular time periods. Daily cumulative dose constraints ensure that the administrated drugs in the model are reasonably close to tested protocols. Let $\rhs_{d,\mathrm{cum}}$ denote the maximum cumulative daily dose of a drug $d \in \dSet$; the corresponding constraint is $\int\limits_{t \in D_{m}} \drug_{d}(t) \, dt \leq \rhs_{d, \mathrm{cum}}, \forAll m \in \DAYS.$ 
\item \textbf{Pill administration}:
Certain drugs are available via oral administration and, therefore, must be administered in discrete amounts~\citep{Hande1998,Reigner2001,Ershler2006,Sharma2006}. These drugs are often recommended to be taken with food. For a drug $d \in \dSet$ that is available in an orally administered pill, let $\para_{d,\mathrm{pill}}$ denote the pill's mass and the integer decision variable $\bin_{d,\mathrm{pill}}(t)$ be the number of pills administered at time $t$; we model this constraint as follows: $\drug_{d}(t) = \para_{d,\mathrm{pill}} \, \bin_{d,\mathrm{pill}}(t),~\bin_{d,\mathrm{pill}}(t) \in \mathbb{Z}_{+}, \forAll t \in \MEALS$, and  $\drug_{d}(t) = 0, \forAll t \notin \MEALS$.
\item \textbf{Rest days} (following treatment administration): Rest periods, in which no amount of a particular drug can be administered, may be mandatory clinically; see e.g., \cite{Baker2006}. We introduce binary decision variables $\bin_{d,\mathrm{rest}}^{m}$ to indicate if a drug $d \in \dSet$ is not administered during day $D_{m}$. Given a mandated number of rest days $\para_{d,\mathrm{rest}}$, we enforce this constraint as follows: 
$\int\limits_{t \in D_m} \drug_{d}(t) \, dt \leq \rhs_{d,\mathrm{cum}}(1 - \bin^{m}_{d,\mathrm{rest}}), 
\sum\limits_{l = 0}^{\min\{\para_{d,\mathrm{rest}},\,M - m\}} (1 - \bin^{m+l}_{d,\mathrm{rest}}) \leq 1,~\bin^{m}_{d, \mathrm{rest}} \in \mathbb{B}, \forAll m \in \DAYS.$ 
\item \textbf{Toxicity}: Drug toxicity is a major consideration in chemotherapy planning. 
We use the white blood cell count as a measure of toxicity and distinguish between neutrophil and lymphocyte, two major white blood cell types, which are responsible for different side effects, namely neutropenia and lymphocytopenia, with different thresholds. Let $N_{\text{neu}}(t)$ and $N_{\text{lym}}(t)$ denote the count of neutrophils and lymphocytes, respectively, at time $t \in [0,\tymEnd]$; we assume neutrophils and lymphocytes account for the fractions $\theta_{\text{neu}}$ and $\theta_{\text{lym}}$ of the total white blood cell count. The neutropenia and lymphocytopenia constraints are as follows: $N_{\text{neu}}(t) \geq \beta_{\text{neu}},~N_{\text{neu}}(t) = \theta_{\text{neu}} \, N_{w}(t),  \forAll t \in [0,T]$, and $N_{\text{lym}}(t) \geq \beta_{\text{lym}},~N_{\text{lym}}(t) = \theta_{\text{lym}} \, N_{w}(t),  \forAll t \in [0,T]$, where $\beta_{\text{neu}}$ and $\beta_{\text{lym}}$ are the clinical thresholds for neutropenia and lymphocytopenia, respectively.
\end{enumerate}

In the presentation of our models, we denote the set of treatment solutions satisfying the operational constraints by $\operSet$. Because the neutropenia and lymphocytopenia constraints are operational constraints, we also include the accompanying pharmacodynamics constraints, i.e., Eqs.~\eqref{wbcModel} and~\eqref{wbcModelDrugs}, within this set when we present our models.
\section{Chemotherapy Optimization Models}
\label{sec:models}
The primary goal of chemotherapy is to reduce the number of cancer cells in the body. There are multiple ways to express this goal; one option is to focus on the end-of-treatment cell count. The objective of our model is to minimize $\sum\limits_{q \in \typeSet} \pop_{q}(\finishTym)$, which is equivalent to minimizing the geometric mean of the cancer cell type populations at the end of the treatment period, i.e., $\big(\prod\limits_{q \in \typeSet} \epop_{q} (T) \big)^{1/|\typeSet|}$. We note that one may prioritize cell types according to their levels of malignancy by considering different coefficients for population variables in the objective function.

Based on the model components described in Section~\ref{sec:prelim}, the combination chemotherapy optimization problem can be formulated as follows: 
\begin{subequations}\label{exactForwardSingleStage}
\begin{align}
\min~~~&\sum\limits_{q \in \typeSet}\pop_{q}(\finishTym)\\
\text{s.t.}~~~&\dot{\conc_{d}}(t) = - \concPar_{d} \, \conc_{d}(t) + \drug_{d}(t)/\vol,\forAll d \in \dSet,~t \in [0,\finishTym], \label{eq:OC_PK}\\
&\conc_{d}(0) = 0,\forAll d \in \dSet, \label{eq:OC_PK_boundary}\\
&\dot{\pop_{q}}(t) = \gompOne \, \Big ( \ln(\epop_{q,\infty}) - \pop_{q}(t) \Big ) - \sum\limits_{d \in \dSet} \dCons_{d,q} \, \exp(-\dExp_{d,q} \, t) \, \econ_{d}(t),\forAll q \in \typeSet,~t \in [0,\finishTym], \label{eq:OC_PD}\\
&\pop_{q}(0) = \ln (\epop_{q,0}),\forAll q \in \typeSet, \label{eq:OC_PD_boundary}\\
&\econ_{d}(t) = \max\{0,~\conc_{d}(t) - \rhs_{d,\mathrm{eff}}\},\forAll d \in \dSet,~t \in [0,\finishTym],\label{eq:OC_effectiveConc}\\
&\big ( {\bf \drug}(t), \, {\bf \conc}(t),\,  \epop_{\text{neu}}(t),\, \epop_{\text{lym}}(t) \big ) \in \operSet,\forAll t \in [0,\finishTym],
\end{align}
\end{subequations}
where ${\bf \drug}(t) = \big(\drug_{1}(t), \ldots, \drug_{|\dSet|}(t)\big)$ and ${\bf \conc}(t) = \big(\conc_{1}(t), \ldots, \conc_{|\dSet|}(t)\big)$ are the variable vectors representing drug administration and concentration, respectively. Nonegativity is enforced, except for $\pop_{q}(t),\forAll q \in \typeSet$, which we consider a part of the definition of $\operSet$.

Formulation~\eqref{exactForwardSingleStage} is an OC problem involving discrete and continuous controls and nonlinear functions, which, given the scale of the instances arising in practice, is extremely hard to solve exactly. We use discretization and linearization techniques to approximate this problem by a MILP formulation, which is significantly more tractable.~\cite{Harrold2009} employed MILP to approximate an OC problem for single-drug chemotherapy optimization. While their model does not consider several operational constraints, they show that this framework provides high-quality approximations for the pharmacokinetics and pharmacodynamics ODE models, i.e., Eqs.~\eqref{eq:OC_PK}--\eqref{eq:OC_effectiveConc}, using two case studies. Here, we also provide some analytical results concerning the approximation quality of such a transformation.

In the rest of this section, we first present our MILP model for combination chemotherapy optimization, including details of the discretization and linearization techniques we use. We extend this model to address uncertainty in the heterogeneity of tumors and present a model for neoadjuvant chemotherapy planning under this uncertainty. Finally, we present our analytical results concerning numerical stability and approximation quality of the proposed models.
\subsection{MILP Model}
\label{sec:deterministic}
To approximate the ODEs in~\eqref{exactForwardSingleStage}, one may use Runge-Kutta (RK) methods as approximation schemes; we refer to~\citep{butcher2007runge} for details. We use (forward) Euler's method, i.e., the first-order RK method, which has been previously used in chemotherapy optimization; see e.g.,~\citep{Harrold2009}. Given a (fixed) time-step $h$, consider the discretization of the planning horizon $[0,T]$ by $\lastTym+1$ points with the index set $\timeSet=\{0,\ldots,\lastTym\}$, where $t(0) = 0$ and  $t(\lastTym) = \finishTym = \lastTym h$. The Euler's approximation of the pharmacokinetics model, i.e., Eqs.~\eqref{eq:OC_PK}--\eqref{eq:OC_PK_boundary}, is as follows:
\begin{equation}\label{eq:PK}
\begin{aligned}
&\conc_{d, \tym+1} = \conc_{d, \tym} - h \, \concPar_{d} \, \conc_{d,\tym} + \frac{\drug_{d,\tym}}{\vol},\forAll d \in \dSet,~s \in \{0,\ldots,S-1\},\\
&\conc_{d,0} = 0,\forAll d \in \dSet.
\end{aligned}
\end{equation}
Note that, in Eq.~\eqref{eq:OC_PK}, $\drug_{d}(t)$ represents the {\em flux} of a drug $d \in \dSet$, i.e., dose administered per unit of time; setting the unit of time equal to $h$, $\drug_{d, \tym}$ denotes the dose administered at (discrete) time $t(s)$. In~\eqref{eq:PK}, $\conc_{d, \tym}$ is the concentration of a drug $d \in \dSet$ at time $t(s)$.

The pharmacodynamics model, i.e., Eqs.~\eqref{eq:OC_PD}--\eqref{eq:OC_PD_boundary}, is approximated in a similar manner:
\begin{equation*}\label{eq:PD}
\begin{aligned}
&\pop_{q, \tym+1} = \pop_{q, \tym} + h\left(\gompOne \, \Big ( \ln(\epop_{q,\infty}) - \pop_{q, \tym} \Big ) - \sum\limits_{d \in \dSet} \dCons_{d,q} \, \exp \big (-\dExp_{d,q} \, t(\tym) \big ) \, \econ_{d,\tym}\right), \forAll q \in \typeSet,~s \in \{0,\ldots,S-1\},\\
&\pop_{q,0} = \ln(\epop_{q,0}),\forAll q \in \typeSet,
\end{aligned}
\end{equation*}
where $\pop_{q, \tym}$ denotes the logarithm of the population of cancer cell type $q \in \typeSet$ at time $t(s)$.

Conventionally, the effective concentration constraints, i.e., Eq.~\eqref{eq:OC_effectiveConc}, are linearized by introducing auxiliary binary variables $\bin_{\econ,d,\tym},\forAll d \in \dSet,~s \in \{0,\ldots,S\}$, as follows:
\begin{subequations}
\begin{align}
\forAll d \in \dSet,~s \in \{0&,\ldots,S\}: \nonumber\\
&\econ_{d,s} \geq 0,\label{relULinear_first} \\
&\econ_{d,s} \geq \conc_{d,s} - \rhs_{d,\mathrm{eff}},\\
&\econ_{d,s} \leq \rhs_{d,\mathrm{conc}} \, \bin_{\econ,d,s},\\
&\econ_{d,s} \leq \conc_{d,s} - \rhs_{d,\mathrm{eff}} + \rhs_{d,\mathrm{conc}} \, (1 - \bin_{\econ,d,s}),\\
&\bin_{\econ,d,\tym}  \in \mathbb{B}, \label{relULinear_last}
\end{align}
\end{subequations}
where $\econ_{d,s}$ is the effective concentration of a drug $d \in \dSet$ at time $t(s)$; recall that $\rhs_{d,\mathrm{conc}}$ denotes the maximum permissible concentration for a drug $d \in \dSet$, i.e., an upper bound on $\conc_{d,s},\forAll s \in \{0,\ldots,S\}$.

The discretization and linearization of operational constraints are straightforward, except for the white blood cell population dynamics model, i.e, Eq.~\eqref{wbcModelDrugs}. We present our model for the white blood cell population dynamics here and provide the details of other operational constraints in the e-companion; see Appendix~\ref{app:modeling}.

Applying Euler's method to Eq.~\eqref{wbcModelDrugs} results in the following discretization of the white blood cell population dynamics model:
\begin{equation} \label{wbcModelDrugsDisc}
\epop_{w,s+1} = \epop_{w,s} + h \left ( \upsilon_{w} - \nu_{w} \, \epop_{w,s} - \sum\limits_{d \in \dSet} \dCons_{d, w} \, \epop_{w,s} \, \conc_{d,s - \tau} \right ),\forAll s \in \{\tau,\ldots,S-1\},
\end{equation}
where $\epop_{w,s}$ denotes the total white blood cell count at time $t(s)$, and $\tau$ corresponds to the time delay $t_{w}$ in the ODE model. To address the bilinearity of Eq.~\eqref{wbcModelDrugsDisc}, we consider two approaches: McCormick relaxation and discretization. In the first approach, we replace each bilinear term $\epop_{w,s} \, \conc_{d,s - \tau}$ in~\eqref{wbcModelDrugsDisc} with a new variable $\mC_{d,\tym}$ and add the corresponding McCormick envelopes~\citep{McCormick1976, al1983jointly} to the formulation. Recall that $\conc_{d,s - \tau} \in [0, \rhs_{d,\mathrm{conc}}]$; we also assume $\epop_{w,s} \in [\rhs_{w}, \epop_{w,0}]$, where $\rhs_{w}$ is a lower bound on the white blood cell count that can be easily obtained from the clinical bounds on the neutrophil and lymphocyte counts, i.e., $\rhs_{w} = \min \{\frac{\beta_{\text{neu}}}{\theta_{\text{neu}}},~\frac{\beta_{\text{lym}}}{\theta_{\text{lym}}}\}$, and $\epop_{w,0}$ denotes the initial count of white blood cells that serves as an upper bound on the white blood cell count during the treatment period. The white blood cell count is discrete, but the continuity assumption is not far from the reality given the magnitude of this quantity, that is $O(10^{9})$ cells per liter. The McCormick relaxation of Eq.~\eqref{wbcModelDrugsDisc} is as follows:
\begin{subequations}\label{mcCormickConstraints}
\begin{align}
&\epop_{w,\tym+1} = \epop_{w,\tym} + h \left( \upsilon_{w} - \nu_{w} \, \epop_{w,s} - \sum\limits_{d \in \dSet} \dCons_{d, w} \, \mC_{d,\tym} \right),\forAll s \in \{\tau,\ldots,S-1\},\\
&\mC_{d,\tym} \geq \rhs_{w} \, \conc_{d,\tym - \tau},\forAll d \in \dSet,~s \in \{\tau,\ldots,S-1\},\\
&\mC_{d,\tym} \geq \epop_{w,0} \, \conc_{d,\tym - \tau} + \rhs_{d,\mathrm{conc}} \, \epop_{w,s} - \epop_{w,0} \, \rhs_{d,\mathrm{conc}},\forAll d \in \dSet,~s \in \{\tau,\ldots,S-1\},\\
&\mC_{d,\tym} \leq \epop_{w,0} \, \conc_{d,\tym - \tau},\forAll d \in \dSet,~s \in \{\tau,\ldots,S-1\},\\
&\mC_{d,\tym} \leq \rhs_{w} \, \conc_{d,s-\tau} + \rhs_{d,\mathrm{conc}} \, \epop_{w,s} - \rhs_{w} \, \rhs_{d,\mathrm{conc}},\forAll d \in \dSet,~s \in \{\tau,\ldots,S-1\}.
\end{align}
\end{subequations}

The drawback of the continuous McCormick relaxation is that the approximation quality of the bilinear sum is not controllable. Therefore, we also consider a modified form of the discretization approach proposed by~\cite{Gupte2013}; given the scale of approximated variables, we alter the constraints provided in this work. As a factor in the bilinear terms, the white blood cell count is approximated by discrete variables within some specified value $\Delta$ of maximum error. We introduce auxiliary binary variables $\bin_{w,s,k},\forAll k \in \{0,\dots,K\}$, which select the approximations of $\epop_{w,\tym}$, and continuous variables $\mv_{d,\tym,k},\forAll k \in \{0,\dots,K\}$, which mirror the value of $\conc_{d,s-\tau}$. Following~\cite{Gupte2013}, the approximation of Eq.~\eqref{wbcModelDrugsDisc} is as follows:
\begin{subequations}\label{GupteMcCormick}
\begin{align}
&\epop_{w,\tym+1} = \epop_{w,\tym} + h \left( \upsilon_{w} - \nu_{w} \, \epop_{w,s} - \sum\limits_{d \in \dSet} \dCons_{d, w} \, \mC_{d,\tym} \right),\forAll s \in \{\tau,\ldots,S-1\},\\
&\epop_{w,\tym} - \left(\rhs_{w} + \sum\limits_{k = 1}^{K} k \, \Delta \, \bin_{w,s,k}\right) \leq \frac{\Delta}{2},\forAll s \in \{\tau,\ldots,S-1\}, \\
&-\epop_{w,\tym} + \left(\rhs_{w} + \sum\limits_{k = 1}^{K} k \, \Delta \, \bin_{w,s,k}\right)\leq \frac{\Delta}{2},\forAll s \in \{\tau,\ldots,S-1\},\\
&\sum\limits_{k = 0}^{K} \bin_{w,s,k} = 1,\forAll s \in \{\tau,\ldots,S-1\},\\
&\mC_{d,s} = \sum\limits_{k = 1}^{K} (\rhs_{w} + k \, \Delta ) \mv_{d,s,k},\forAll d \in \dSet,~s \in \{\tau,\ldots,S-1\},\\
&\mv_{d,s,k} \leq \rhs_{d,\mathrm{conc}} \, \bin_{w,s,k},\forAll d \in \dSet,~s \in \{\tau,\ldots,S-1\},~k \in \{1,\dots,K\},\\
&\mv_{d,s,k} \leq \conc_{d,s-\tau},\forAll d \in \dSet,~s \in \{\tau,\ldots,S-1\},~k \in \{1,\dots,K\},\\
&\mv_{d,s,k} \geq \conc_{d,s-\tau} + \rhs_{d,\mathrm{conc}} \, (\bin_{w,s,k} - 1),\forAll d \in \dSet,~s \in \{\tau,\ldots,S-1\},~k \in \{1,\dots,K\},\\
&\mv_{d,s,k} \geq 0,\forAll d \in \dSet,~s \in \{\tau,\ldots,S-1\},~k \in \{1,\dots,K\},\\
&\bin_{w,s,k} \in \mathbb{B},\forAll s \in \{\tau,\ldots,S-1\},~k \in \{1,\dots,K\}.
\end{align}
\end{subequations}
In~\eqref{GupteMcCormick}, the quantity $\rhs_{w} + \sum\limits_{k = 1}^{K} k \, \Delta \, \bin_{w,s,k}$ approximates the value of $\epop_{w,\tym}$; the variable $\mv_{d,s,k}$ equals $\conc_{d,s-\tau}$ if and only if $\bin_{w,s,k} = 1$ (0 otherwise), and $\mC_{d,s}$ approximates the bilinear term $\epop_{w,s}\conc_{d,s-\tau}$.

Before the effect of drugs on white blood cells starts, i.e, $\forAll s \in \{0,\ldots,\tau-1\}$, we have $\epop_{w,\tym+1} = \epop_{w,\tym} + h \left( \upsilon_{w} - \nu_{w} \, \epop_{w,s} \right )$, where $\epop_{w,0}$ equals the count of white blood cells before treatment. The discretized neutropenia and lymphocytopenia constraints are as follows: $N_{\text{neu},s} \geq \beta_{\text{neu}},~N_{\text{neu},s} = \theta_{\text{neu}} \, N_{w,s}, \forAll s \in \{0,\dots,S\}$, and $N_{\text{lym},s} \geq \beta_{\text{lym}},~N_{\text{lym},s} = \theta_{\text{lym}} \, N_{w,s},  \forAll s \in \{0,\dots,S\}$.

Formulation~\eqref{eq:detMILP} presents our MILP model for combination chemotherapy optimization, as the result of described discretization and linearization techniques applied to~\eqref{exactForwardSingleStage}. In this formulation, we use $\doperSet$ to denote the set of treatment solutions satisfying the discretized operational constraints, including the constraint sets~\eqref{mcCormickConstraints} or~\eqref{GupteMcCormick}.
\begin{subequations}\label{eq:detMILP}
\begin{align}
\min~~~&\sum\limits_{q \in \typeSet}\pop_{q,\lastTym}\\
\text{s.t.}~~~&\conc_{d, \tym+1} = \conc_{d, \tym} - h \, \concPar_{d} \, \conc_{d,\tym} + \drug_{d,\tym}/\vol,\forAll d \in \dSet,~s \in \{0,\ldots,S-1\}, \label{eq:MILP_PK}\\
&\conc_{d,0} = 0,\forAll d \in \dSet, \label{eq:MILP_PK_boundary}\\
&\pop_{q, \tym+1} = \pop_{q, \tym} + h\left(\gompOne \, \Big ( \ln(\epop_{q,\infty}) - \pop_{q, \tym} \Big ) - \sum\limits_{d \in \dSet} \dCons_{d,q} \, \exp \big (-\dExp_{d,q} \, t(\tym) \big ) \, \econ_{d,\tym}\right), \nonumber \\
&\qquad\forAll q \in \typeSet,~s \in \{0,\ldots,S-1\}, \label{eq:MILP_PD} \\
&\pop_{q,0} = \ln(\epop_{q,0}),\forAll q \in \typeSet, \label{eq:MILP_PD_boundary}\\
&\eqref{relULinear_first}-\eqref{relULinear_last},\forAll d \in \dSet,~s \in \{0,\ldots,S\}, \label{eq:MILP_effectiveConc}\\
&({\bf \drug},\, {\bf \conc},\, {\bf \epop}_{\text{neu}},\, {\bf \epop}_{\text{lym}}) \in \doperSet,
\end{align}
\end{subequations}
where ${\bf \drug} = \big[\drug_{d,\tym}, d \in \dSet, \tym \in \timeSet\big]$ and ${\bf \conc} = \big[\conc_{d,\tym}, d \in \dSet, \tym \in \timeSet\big]$ are the variable matrices representing drug administration and concentration, respectively, and ${\bf \epop}_{\text{neu}} = (\epop_{\text{neu},0}, \ldots, \epop_{\text{neu},\lastTym})$ and ${\bf \epop}_{\text{lym}} = (\epop_{\text{lym},0}, \ldots, \epop_{\text{lym},\lastTym})$ are the variable vectors of the neutrophil and lymphocyte count, respectively. Eqs.~\eqref{eq:MILP_PK}--\eqref{eq:MILP_PK_boundary} are the discretized pharmacokinetics model \big(see~\eqref{eq:OC_PK}--\eqref{eq:OC_PK_boundary}\big), Eqs.~\eqref{eq:MILP_PD}--\eqref{eq:MILP_PD_boundary} are the discretized pharmacodynamics model for cancer cells \big(see~\eqref{eq:OC_PD}--\eqref{eq:OC_PD_boundary}\big), and the constraint set~\eqref{eq:MILP_effectiveConc} is the linearized model for effective drug concentration \big(see~\eqref{eq:OC_effectiveConc}\big).

The approximation quality and computational burden of formulation \eqref{eq:detMILP} depends on the choice of discretization time-step $h$ and the method of bilinearity approximation, i.e, McCormick relaxation or the discretization technique, as well as the choice of $\Delta$ in the latter method. We present the results of our computational experiments with respect to these factors in Section~\ref{sec:results}. Note that if $h' < h$, then the discretization with $h'$ provides a finer resolution of the time domain. However, models with differing time-steps are not necessarily relaxations of each other. Similarly, the discretized formulation~\eqref{eq:detMILP} approximates the continuous formulation~\eqref{exactForwardSingleStage}, but it is not necessarily a relaxation of the continuous problem.
\subsection{Neoadjuvant Chemotherapy under Uncertain Tumor Heterogeneity}
\label{sec:stochastic}
In advanced stages of cancer, chemotherapy is often used as a form of neoadjuvant therapy to reduce a tumor to an operable size prior to tumor removal surgery~\citep{Senkus2015}. In this section, we extend the proposed MILP model~\eqref{eq:detMILP} to address the uncertainty of tumor heterogeneity in neoadjuvant chemotherapy planning. The cell composition of a tumor can be characterized through multiple biopsies and next-generation sequencing (NGS) technologies~\citep{gerlinger2012,PIRAINO2019}. Multiple biopsies, however, come with a higher risk of needle tract seeding and cancer metastasis~\citep{tyagi2014}. Hence, with a limited number of biopsies, tumor heterogeneity remains uncertain for treatment planning~\citep{abecassis2019}. Under this uncertainty, the objective is to find a drug administration regimen that reduces the tumor to a clinically determined operable size with a high probability, which can be expressed by the following chance constraint:
\begin{align} \label{ogChanceConstraint}
\Prob\left\{\sum\limits_{q \in \typeSet} \epop_{q,\lastTym} \leq \epop_{\text{surg}}\right\} \geq 1 - \epsilon, 
\end{align}
where $\epop_{\text{surg}}$ denotes the clinically determined operable size for the tumor, and $\epsilon$ is the probability of not meeting the target at the end of treatment period.

Let $\ipop \in \mathbb{R}_{+}^{|\typeSet|}$ be a discrete random variable describing tumor heterogeneity, which can take on a value from the finite set $\{\ipop^{(1)},\ldots,\ipop^{(K)}\}$, and denote the probability of a scenario $k \in \{1,\ldots,K\}$ by $\mu^{(k)}=\Prob\{\ipop = \ipop^{(k)}\}$. Given the logarithmic transformation of the cancer cell population variables in our model, we use the following, more conservative constraints to enforce~\eqref{ogChanceConstraint}:
\begin{subequations}\label{chanceSurgery}
\begin{align}
&\pop^{(k)}_{q,\lastTym} \leq \pop_{\text{surg}} + \ln \left ( \dfrac{N^{(k)}_{q,0}}{\sum_{q \in \typeSet} N^{(k)}_{q,0}} \right ) + \pop_{q,\infty} (1 - \bin_{\text{surg}}^{(k)}),\forAll q \in \typeSet,~k \in \{1,\dots,K\}, \label{chanceSurgery.1}\\
&\sum\limits_{k = 1}^{K} \mu^{(k)} \, \bin^{(k)}_{\text{surg}} \geq 1 - \epsilon, \label{chanceSurgery.2}\\
&\bin^{(k)}_{\text{surg}} \in \mathbb{B},\forAll k \in \{1,\dots, K\}, \label{chanceSurgery.3}
\end{align}
\end{subequations}
where $\pop_{\text{surg}} = \ln (\epop_{\text{surg}})$, $\pop_{q,\infty} = \ln (\epop_{q,\infty})$, and the superscript $(k)$ denotes the value of previously defined variables under a realization scenario $k \in \{1,\ldots,K\}$; the binary variable $\bin^{(k)}_{\text{surg}}$ indicates whether the target is met under a scenario $k \in \{1,\dots, K\}$. It can be easily verified that, given $\bin^{(k)}_{\text{surg}} = 1$ for some scenario $k$, a treatment solution satisfying the set of constraints~\eqref{chanceSurgery.1} will also satisfy $\sum\limits_{q \in \typeSet} \epop^{(k)}_{q,\lastTym} \leq \epop_{\text{surg}}$; these constraints become trivial if $\bin^{(k)}_{\text{surg}} = 0$.

A conventional objective for such a chance-constrained optimization model is to maximize the probability of meeting the target, equivalently to minimize $\epsilon$. A drawback of this objective is that if the tumor size is not far from the clinical target, the model may achieve an objective of $\epsilon = 0$, with a solution that may not lead to significant tumor shrinkage. Recall that the main clinical objective of chemotherapy is to reduce cancer cell population as much as possible~\citep{howChemo2021}. In this regard, we also consider a shrinkage-based objective similar to our deterministic model. Note that, given a fixed value for $\epsilon$, obtaining a success probability of $1 - \epsilon$ is guaranteed through the chance constraints~\eqref{chanceSurgery.1}--\eqref{chanceSurgery.3} regardless of the objective function. Section~\ref{sec:results} presents the results of our numerical study with both shrinkage-based and probability-based objectives. Here, we present our stochastic model with an objective that minimizes the cancer cell population at the end of treatment under the most likely scenario. We assume $\mu^{(1)} \geq \mu^{(k)},\forAll k \in \{1,\ldots,K\}$; the neoadjuvant chance-constrained MILP model is as follows:
\begin{subequations}\label{eq:stoMILP}
\begin{align}
\min~~~&\sum\limits_{q \in \typeSet}\pop^{(1)}_{q,\lastTym}\\
\text{s.t.}~~~&\conc_{d, \tym+1} = \conc_{d, \tym} - h \, \concPar_{d} \, \conc_{d,\tym} + \drug_{d,\tym}/\vol,\forAll d \in \dSet,~s \in \{0,\ldots,S-1\}, \label{eq:sto_PK}\\
&\conc_{d,0} = 0,\forAll d \in \dSet, \label{eq:sto_PK_boundary}\\
&\pop^{(k)}_{q, \tym+1} = \pop^{(k)}_{q, \tym} + h\left(\gompOne \, \Big ( \ln(\epop^{(k)}_{q,\infty}) - \pop^{(k)}_{q, \tym} \Big ) - \sum\limits_{d \in \dSet} \dCons_{d,q} \, \exp \big (-\dExp_{d,q} \, t(\tym) \big ) \, \econ_{d,\tym}\right), \nonumber\\
&\qquad\forAll q \in \typeSet,~s \in \{0,\ldots,S-1\},~k \in \{1,\ldots,K\}, \label{eq:sto_PD}\\
&\pop^{(k)}_{q,0} = \ln(\epop^{(k)}_{q,0}),\forAll q \in \typeSet,~k \in \{1,\ldots,K\}, \label{eq:sto_PD_boundary}\\
&\eqref{relULinear_first}-\eqref{relULinear_last},\forAll d \in \dSet,~s \in \{0,\ldots,S\}, \label{eq:sto_effectiveConc}\\
&\pop^{(k)}_{q,\lastTym} \leq \pop_{\text{surg}} - \ln \left ( \dfrac{N^{(k)}_{q,0}}{\sum_{q \in \typeSet} N^{(k)}_{q,0}} \right ) + \pop_{q,\infty} (1 - \bin_{\text{surg}}^{(k)}),\forAll q \in \typeSet,~k \in \{1,\dots,K\}, \label{eq:sto_chanceSurgery.1}\\
&\sum\limits_{k = 1}^{K} \mu^{(k)} \, \bin^{(k)}_{\text{surg}} \geq 1 - \epsilon, \label{eq:sto_chanceSurgery.2}\\
&\bin^{(k)}_{\text{surg}} \in \mathbb{B},\forAll k \in \{1,\dots, K\}, \label{eq:sto_chanceSurgery.3}\\
&({\bf \drug},\, {\bf \conc},\, {\bf \epop}_{\text{neu}},\, {\bf \epop}_{\text{lym}}) \in \doperSet. \label{eq:sto_last}
\end{align}
\end{subequations}
Note that the variables appearing in the operational constraints, i.e., ${\bf \drug}$, ${\bf \conc}$, ${\bf \epop}_{\text{neu}}$, and ${\bf \epop}_{\text{lym}}$, do not depend on the scenarios.
\subsection{Analytical Results}
\label{sec:structural}
Next, we present analytical results for the proposed combination chemotherapy optimization models. A main result is an error bound for the Euler's method approximation of the cell population state variables, which depend on the Euler's method approximation of the drug concentration state variables, presented through Theorem~\ref{nestedEulerApprox} and Corollary~\ref{nestedPopApprox}. All proofs are provided in the e-companion; see Appendix~\ref{app:proofs}.

Theorems~\ref{uniqueCSolution} and~\ref{uniquePSolution} show that if drug administration is continuous, the state variables for drug concentration and cell population are uniquely defined by the control variables governing drug administration in the base formulation~\eqref{exactForwardSingleStage}. Although administration of oral drugs are inherently discontinuous, they can be approximated arbitrarily well by continuous functions.

\begin{theorem}\thlabel{uniqueCSolution}
Suppose that, for a drug $d \in \dSet$, the administration function $\drug_{d}$ is continuous in time. Then, the differential equation
\begin{align*}
\dot{\conc_{d}}(t) = -\concPar_{d} \, \conc_{d}(t) + \drug_{d}(t)/\vol,~t \in [0, \finishTym],
\end{align*}
governing the drug concentration function $\conc_{d}$, has a unique solution.
\end{theorem}

\begin{theorem}\thlabel{uniquePSolution}
Suppose that the administration functions for all drugs, i.e., $\drug_{d},\forAll d \in \dSet$, are continuous in time. Then, for each cancer cell type $q \in \typeSet$, the differential equation
\begin{align*}
\dot{\pop_{q}}(t) = \gompOne\big(\pop_{q,\infty} - \pop_{q}(t)\big) - \sum\limits_{d \in \dSet} \dCons_{d,q}\exp(-\dExp_{d,q}\,t)\,\econ_{d}(t),~t \in [0, \finishTym],
\end{align*}
governing the cell population function $\pop_{q}$, has a unique solution. 
\end{theorem}

The next results concern the stability of Euler's method. Absolutely stable numerical methods produce ``reasonable results" for suitable time-step values \citep{LeVeque2007}.

\begin{theorem}\thlabel{thm:concIsStable}
Let $\{\drug_{s}\}_{s \in \mathbb{Z}_{+}}$ be a bounded sequence, and $\concPar, h, \vol > 0$. Under the stability condition $h < \frac{2}{\concPar}$, the difference equation 
$$
\conc_{s+1} = \conc_{s} - h \, \concPar \, \conc_{s} + \drug_{s}/\vol,
$$
is absolutely stable, for all $s \in \mathbb{Z}_{+}$.  
\end{theorem}

\begin{theorem}\thlabel{thm:popIsStable}
Let $\{F_{s}\}_{s \in \mathbb{Z}_{+}}$ be a bounded sequence and $\gompOne, h > 0$. Under the stability condition $h < \frac{2}{\gompOne}$, the difference equation 
$$
\pop_{s+1} = \pop_{s} + h \, \big(\gompOne \, (\pop_{\infty} - \pop_{s} ) - F_{s}\big), 
$$
is absolutely stable, for all $s \in \mathbb{Z}_{+}$.
\end{theorem}

Concerning the stochastic model~\eqref{eq:stoMILP}, the following result indicates that if a feasible solution's effective concentration is dominated by another feasible solution's effective concentration, the second solution has a smaller end-of-treatment cancer cell population in all scenarios, regardless of the objective function.

\begin{theorem}\thlabel{thm:moreEffectiveConc}
Consider the stochastic model~\eqref{eq:stoMILP}, with $\gompOne h \leq 1$. Let $({\bf \econ}^{[1]}$, ${\bf \pop}^{[1]})$ and $({\bf \econ}^{[2]}$, ${\bf \pop}^{[2]})$ each be components of different feasible solutions. Suppose $\econ^{[1]}_{d,\tym} \geq \econ^{[2]}_{d,\tym}$, for all $d \in \dSet,\tym \in \{0,\dots,\lastTym\}$. Then $\pop_{q,\lastTym}^{[1],(k)} \leq \pop_{q,\lastTym}^{[2],(k)}$, for all $q \in \typeSet, k \in \{1,\dots,K\}$.
\end{theorem}

Finally, Theorem~\ref{nestedEulerApprox} provides a global error bound for Euler's method applied to state variables that depend on another state variable estimated by Euler's method.

\begin{theorem}\thlabel{nestedEulerApprox}
Consider the system of differential equations
\begin{subequations}
\begin{align*}
\dot{y}(t) &= f(t, y, z),~y(0) = y_{0},\\
\dot{z}(t) &= g(t, z),~z(0) = z_{0},
\end{align*}
\end{subequations}
and the Euler's approximation $\{(y_{s}, z_{s})\}_{s = 0}^{S}$ with step size $h$, given by
$y_{s+1} = y_{s} + h \, f\big(t(s), y_{s}, z_{s}\big)$ and $z_{s+1} = z_{s} + h \, g\big(t(s), z_{s}\big)$. 
Let $\lambda_{z} = \max\{|z_{s} - z_{0}|, s \in \{0,\dots,S\}\}$, and suppose $g$ is continuous in both variables and Lipschitz continuous in its second variable, i.e., there exists $L_{g} > 0$ such that for all $t \in [0, T]~\text{and}~u, v \in \mathbb{R}$ with $|u - z_{0}| \leq \lambda_{z}, |v - z_{0}| \leq \lambda_{z}$,
\begin{align*}
|g(t,u) - g(t,v)| \leq L_{g}|u - v|.
\end{align*} 
Similarly, suppose $f$ is continuous in all variables and Lipschitz continuous (with respect to the $\ell_{1}$ norm) in its second and third variables with constant $L_{f}$. Furthermore, suppose $y$ and $z$ are twice continuously differentiable. Then, for all $s \in \{0,\dots,S\}$,
\begin{align*}
|y_{s} - y(t(s))| \leq \frac{h}{2}\left(\frac{\alpha_{z}}{L_{g}}(e^{L_{g}T} - 1) + \frac{\alpha_{y}}{L_{f}}\right)(e^{L_{f}T} - 1), 
\end{align*}
where $\alpha_{z} = \max\limits_{\tau \in [0,T]} |\ddot{z}(\tau)|$ and $\alpha_{y} = \max\limits_{\tau \in [0,T]} |\ddot{y}(\tau)|$.
\end{theorem}

\begin{corollary} \thlabel{nestedPopApprox}
Let $\conc(t)$ and $\pop_{q}(t),\forAll q \in \typeSet$, be the state variable functions for drug concentration and cell population, respectively, in an optimal solution to the (single-drug) chemotherapy optimization problem~\eqref{exactForwardSingleStage} without the effective concentration and operational constraints. Furthermore, suppose that $\conc$ and $\pop_{q},\forAll q \in \typeSet$, are twice continuously differentiable, and let $\tilde{\conc}$ and $\tilde{\pop}_{q},\,\forAll q \in \typeSet$, be the corresponding Euler's approximations with time-step $h$. Then,
\begin{align*}
\bigg\vert \sum\limits_{q \in \typeSet} \tilde{\pop}_{q,S} - \sum\limits_{q \in \typeSet} \pop_{q}(T) \bigg\vert &\leq \sum\limits_{q \in \typeSet} \frac{h}{2}\left(\frac{\alpha_{\conc}}{|\concPar|}(e^{|\concPar|T} - 1) + \frac{\alpha_{q}}{\max\{|\dCons_{q}|, |\gompOne|\}}\right)(e^{\max\{|\dCons_{q}|, |\gompOne|\}T} - 1),
\end{align*}
where $\alpha_{\conc} = \max\limits_{\tau \in [0,T]}|\ddot{\conc}(\tau)|$ and $\alpha_{q} = \max\limits_{\tau \in [0,T]} |\ddot{\pop}_{q}(\tau)|,\forAll q \in \typeSet$.
\end{corollary}
\section{Model Calibration}
\label{sec:calib}
Though the proposed framework applies to many forms of cancer, we specify our numerical study for breast cancer, which kills more than 40,000 American women annually~\citep{CancerFacts2021}. We include three breast cancer drugs in our study: capecitabine, docetaxel, and etoposide; they are labeled 1, 2, and 3, respectively. To account for heterogeneity in the cell population, we include four tumor cell types: 
$0 \equiv ``\text{no resistance, all drugs}"$, 
$1 \equiv ``\text{resistance, capecitabine only}"$,
$2 \equiv ``\text{resistance, docetaxel only}"$, and
$3 \equiv ``\text{resistance, etoposide only}"$.
We do not consider the case in which a cell type is resistant to multiple drugs because the drugs we consider have different mechanisms to attack tumor cells \citep{Luqmani2005}.

The initial tumor size can vary substantially, depending on the progression of the disease, at the start of a treatment.~\cite{Norton1988} estimates the initial population size of untreated breast cancer patients as $\epop_{0} = 4.8\bigcdot 10^{9}$ cells. This estimate is supported by the fact that tumor detection usually does not occur before the tumor has $10^{9}$ cells, 30 generations after the first malignant cell \citep{Asachenkov1994,Cameron1997}. We use $2^{30} \approx 10^{9}$ cells as the initial cancer cell population. We estimate tumor heterogeneity, i.e., $\epop_{q,0},\forAll q \in \{0,1,2,3\}$, through a branching process; e.g.,~\citep{Kimmel2015}. By simulating this process, we generate multiple scenarios concerning tumor heterogeneity and estimate the probability of realization for each scenario. The description of branching process and estimation of the initial cancer cell populations are provided in the e-companion; see Appendix~\ref{app:branching}. We use these scenarios directly in our numerical study with the chance-constrained model~\eqref{eq:stoMILP}; for the deterministic model~\eqref{eq:detMILP}, we use the (empirical) mean cell count for each tumor cell type in these scenarios.~\cite{Norton1988} estimates the steady-state tumor size as $\epop_{\infty} = 3.1\bigcdot 10^{12}$ cells, which is supported by the maximum tumor size detected in mammograms \citep{Cameron1997}; we set this asymptotic limit at approximately $10^{12}$ cells in our models. Following~\cite{Harrold2009}, we estimate the Gompertz shape parameter by $\Lambda = \dfrac{1}{\tau} \ln \left ( \frac{\ln(N_{\infty} / N_0)}{\ln(N_{\infty} / 2N_0)} \right )$, where $\tau$ is the ``doubling time'' of the tumor, set equal to 5 months in their work. For the white blood cell population dynamics, we follow~\cite{Iliadis2000} and use $\epop_{w,0} = 8\bigcdot10^{9}$ cells per liter as the initial population and $\upsilon_{w} = 1.2\bigcdot10^{9}$ cells per liter per day and $\nu_{w} = 0.15$ per day as the production and turnover rates, respectively. Table \ref{table:populationParameters} summarizes the parameters used for population dynamics in our numerical study; see Appendix~\ref{app:estimate}.

We use the values reported in the literature for pharmacokinetics parameters, i.e., elimination rate $\xi_d,\forAll d \in \dSet$, effect compartment $\vol$, and effectiveness threshold $\rhs_{d,\mathrm{eff}},\forAll d \in \dSet$,~\citep{Iliadis2000,Frances2011} and estimate the pharmacodynamics parameters based on published clinical data. In particular, we use the following administration regimens and the response rate observed in the corresponding clinical trial to estimate the effect of cytotoxic drugs on the cancer cell populations:
\begin{itemize}
\item Capecitabine~\citep{o2001randomized}: 
1255 $\text{mg}/\text{m}^2$ twice daily, 6 cycles of a two-week treatment period followed by a one-week rest period, response rate of 30\%.
\item Docetaxel~\citep{chan1999}: 
100 $\text{mg}/\text{m}^2$, 7 cycles of one-hour infusion every three weeks, response rate of 47\%.
\item Etoposide~\citep{yuan2015}: 60 $\text{mg}/\text{m}^2$ daily, 7 cycles of a 10-day treatment period followed by a 11-day rest period, response rate of 9\%.
\end{itemize}
In clinical studies, the dose administration is commonly reported based on body surface area; we use 1.7 $\text{m}^2$ as an average person's body surface area in our study~\citep{bonate2011}. The treatment (partial) response rate is defined as the percentage of the patients participating in a clinical trial who show 50\% or more decrease in tumor size (diameter) as a result of the therapy~\citep{world1979handbook}. Details of pharmacodynamics parameters estimation are provided in the e-companion; see Appendix~\ref{app:estimate}. Given the narrow therapeutic margin of cytotoxic drugs, we make a conservative assumption that these drugs have the same fractional kill effect on the white blood cells as the cancer cells. Table~\ref{table:PkPDParameters} displays the pharmacokinetics and pharmacodynamics parameters in our numerical study; see Appendix~\ref{app:estimate}.

We also use the simulation results of the aforementioned clinical trials to determine the operational constraints parameters concerning maximum drug concentration, maximum infusion rate, and maximum cumulative daily dose. This ensures that an (optimal) treatment solution stays within the common range of drug administration in practice. Capecitabine and etoposide are orally administered via pills of size 500 mg and 50 mg, respectively~\citep{Hande1998,Sharma2006}. Oral drugs are often taken with food, so we designate three time points within each day at which the oral drugs may be taken. Docetaxel use is constrained by a week-long rest period after each administration day. Finally, we use $\beta_{\text{neu}} = 2.5\bigcdot 10^{9}$ and $\beta_{\text{lym}} = 1\bigcdot 10^{9}$ cells per liter as the neutropenia  and lymphocytopenia thresholds, respectively~\citep{rosado2011hyper,mitrovic2012prognostic}. Table \ref{table:oprParameters} in Appendix~\ref{app:estimate} summarizes the operational parameters used in our numerical study.
\section{Computational Results}
\label{sec:results}
In this section, we present the results of our numerical study, given the model specification and calibration details provided in Section~\ref{sec:calib}. To solve the proposed mixed-integer linear programs, we used Gurobi 9.1.2 with default parameters on a machine with Intel(R) Core(TM) i7-3520M CPU @ 2.90 GHz and restricted the solution time for each instance to a maximum of two hours. Unless otherwise stated, we set the time-step parameter $h$ equal to one hour and use the constraint set~\eqref{GupteMcCormick} with $\Delta = \frac{1}{20}(\epop_{w,0} - \rhs_{w})$ to approximate the bilinear terms; to avoid explosion of the binary variables introduced in~\eqref{GupteMcCormick}, we set the time-step length of one day for the white blood cell count model. We present our results on the computational performance of the models under different discretization resolutions and bilinear approximation methods separately. Finally, we consider a treatment period of 21 days in our numerical study, which is the common length for a chemotherapy cycle~\citep{Ershler2006}.

\begin{figure}
\centering
\begin{subfloat}[Administration]{
\centering
\includegraphics[width = .45\textwidth]{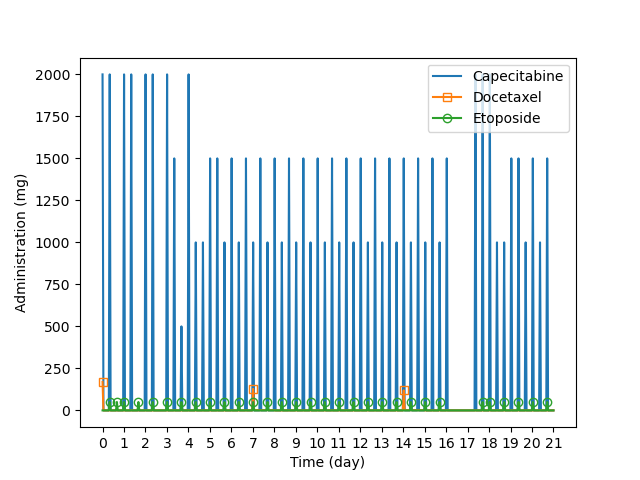}
\label{fig:deterministicDrug}}
\end{subfloat}
\begin{subfloat}[Concentration]{
\centering
\includegraphics[width = .45\textwidth]{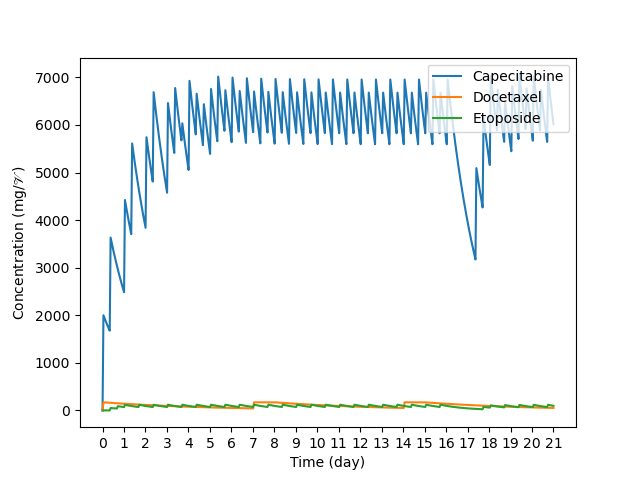}
\label{fig:deterministicConc}}
\end{subfloat}
\caption{Optimal drug administration and the corresponding concentration, given by the deterministic model~\eqref{eq:detMILP}}
\label{fig:deterministicDrugAndConc}
\end{figure}

Figure~\ref{fig:deterministicDrugAndConc} displays the optimal drug administration and concentration over the treatment period using the deterministic model~\eqref{eq:detMILP}. More details are provided in the e-companion; see Figure~\ref{fig:optAdmin} in Appendix~\ref{app:figs}. From the optimal solution, the highest doses occur at the beginning of the treatment period for all drugs. For capecitabine and etoposide that can be administered frequently, the initial high doses make drug concentrations reach the maximum permissible levels, and after that, the administration regimens force relatively constant concentrations until an induced rest period, on day 17; it is important to note that this rest period is not a modeling mandate. In fact, without an explicit rest constraint on the oral drug administration, the optimal solution shows a necessary rest period for these drugs to avoid the violation of toxicity (neutropenia) constraints, akin to mandated policies in clinical practice. Docetaxel, the intravenous drug, is administered weekly starting on the first day, due to the mandated rest period constraints. 

\begin{figure}
\centering
\begin{subfloat}[Cancer cell population (logarithmic scale)]{
\centering	
\includegraphics[width = .45\textwidth]{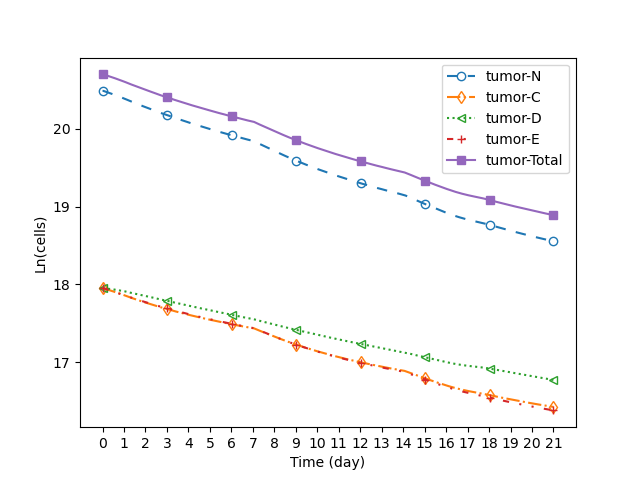}
\label{fig:deterministicPop}}
\end{subfloat}
\begin{subfloat}[White blood cell population]{
\centering
\includegraphics[width = .45\textwidth]{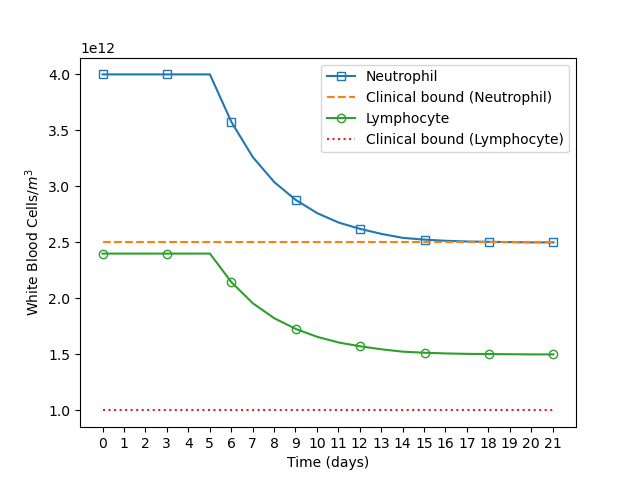}
\label{fig:deterministicWBC}}
\end{subfloat}
\caption{Effect of the optimal drug administration on tumor and white blood cell populations, given by the deterministic model~\eqref{eq:detMILP}}
\label{fig:deterministicPoplogAndWBC}
\end{figure}

The treatment effects on the cancer and white blood cell populations are illustrated in Figure~\ref{fig:deterministicPoplogAndWBC}. In this figure, $N \equiv$ non-resistant, $C \equiv$ capecitabine-resistant, $D \equiv$ docetaxel-resistant, and $E \equiv$ etoposide-resistant. All cancer cell types decrease fairly consistently (with respect to the logarithm) over time. The capecitabine- and etoposide-resistant cell types have similar outcomes; the docetaxel-resistant cell type shows a lower level of response to the treatment. Patently, the steepest descent belongs to the non-resistant cancer cell type. Regarding the white blood cell population, the neutropenia constraint is tight at the optimal solution, which is consistent with the clinical observation that neutropenia is often a toxicity of concern in chemotherapy~\citep{Pizzo1993,Kosaka2015,kasi2018}.

\begin{figure}[t]
\centering
\begin{subfloat}[Capecitabine administration and concentration]{
\centering	
\includegraphics[width = .45\textwidth]{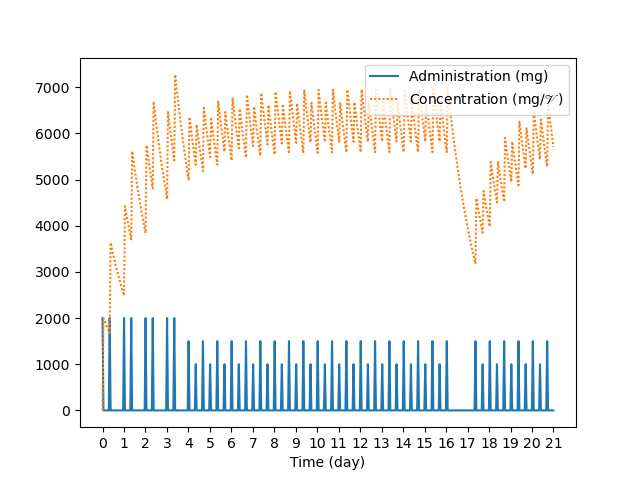}}
\end{subfloat}
\begin{subfloat}[Docetaxel administration and concentration]{
\centering
\includegraphics[width = .45\textwidth]{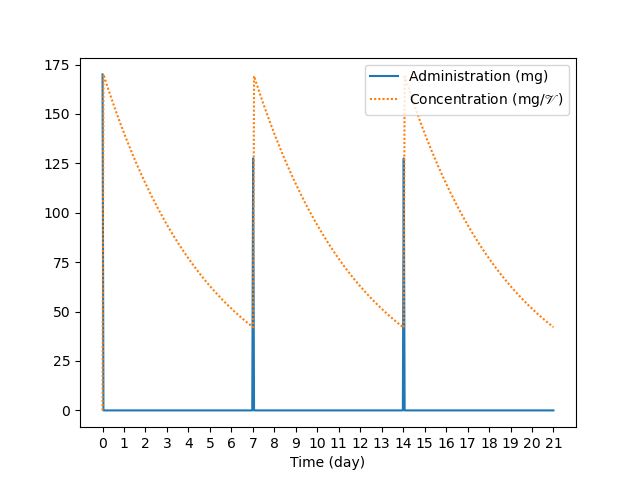}}
\end{subfloat}
\begin{subfloat}[Etoposide administration and concentration]{
\centering	
\includegraphics[width = .45\textwidth]{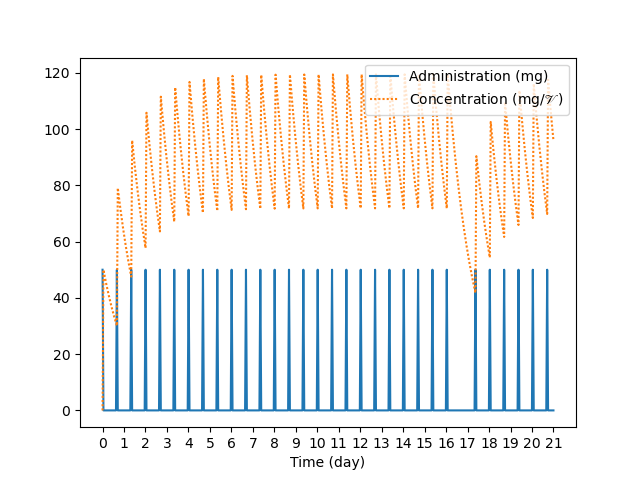}}
\end{subfloat}
\begin{subfloat}[Cancer cell population (logarithmic scale)]{
\centering
\includegraphics[width = .45\textwidth]{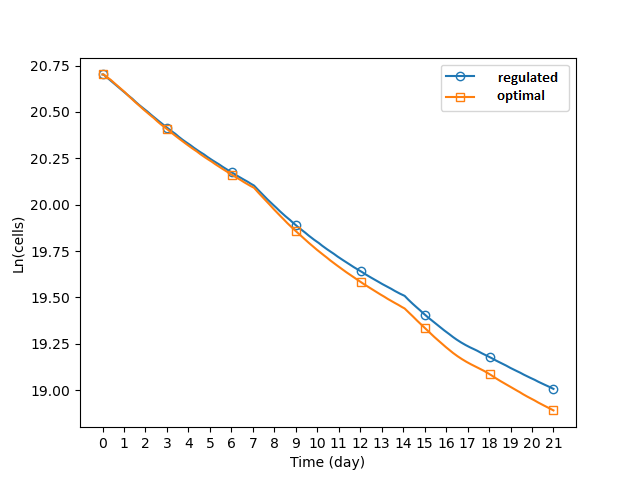}}
\end{subfloat}
\caption{Regulated drugs administration (and concentration) based on the optimal solution of the deterministic model~\eqref{eq:detMILP} and the corresponding tumor-shrinkage outcome compared with that of the optimal solution}
\label{fig:reg}
\end{figure}

Regularity of drug administration is a logistic consideration for oral drugs, as they are usually taken without direct medical supervision~\citep{urquhart1998contending}. In this regard, we have considered converting the optimal solution presented in Figure~\ref{fig:deterministicDrugAndConc} to a regulated administration plan. Figure~\ref{fig:reg} shows the result along with the corresponding tumor-shrinkage outcome; the difference in the cancer cell population between the optimal and regulated plans translates to less than 1 mm change in the tumor diameter. Note that, in construction of the regulated administration, we keep the induced rest period on day 17 because ignoring this period leads to a violation of the neutropenia constraint. We acknowledge that regularity of the administration regimen for oral drugs can be enforced by additional constraints; however, such patterns must be devised carefully, as they can be too restrictive on the outcome. Finally, we point out that, given the pill sizes and maximum dose and concentration constraints for the oral drugs that we consider, the difference between the optimal and regulated plans mainly concerns the administration of capecitabine, which is weaker than etoposide based on their fractional kill effect parameter values. 

In the next phase of our experiments, we investigated the impact of different discretization resolutions on the computational performance of~\eqref{eq:detMILP} by varying the time-step parameter $h$ from four hours (240 minutes) to 15 minutes. Table~\ref{tab:byH_Deterministic} summarizes the results; the gap of $0.01\%$ is the solver's default value within which it considers an incumbent solution optimal. Given the initial objective value of $\sum\limits_{q =0}^3\pop_{q,0} = 74.34$ at the start of the treatment, there is little evidence that the discretization resolution (within the tested range) substantially impacts the optimal objective value, in accordance with our stability results; all tested time-step values satisfy the stability conditions of~\thref{thm:concIsStable,thm:popIsStable}. However, it does impact the model's solvability because as $h$ decreases, the number of variables and constraints increase. The gap for the model with $h=15$ minutes was $0.09\%$, the highest among the tested values.
\begin{table}
\centering
\caption{\label{tab:byH_Deterministic} Solver statistics for different time-step values ($h$)}
\begin{tabular}{ccccccccc}
\toprule
$h$\,(minutes) & Obj Val & Run-time\,(s) & Cons  & \phantom{\large -h}Vars\phantom{\huge $._g$} & IVars &BVars &Gap\\
\hline
240 & 68.09 & 1163 & 8547 & 3759 & 714 & 588 & $<$ 0.01\%\\
120 & 68.12 & 7200 & 11571 & 5523 & 840 & 714 & \quad 0.02\%\\
60 & 68.13 & 2029 & 17619 & 9051 & 1092 & 966 & $<$ 0.01\%\\
30 & 68.21 & 7200 & 29715 & 16107 & 1596 & 1470 & \quad 0.02\%\\
15 & 68.24 & 7200 & 53907 & 30219 & 2604 & 2478 & \quad 0.09\%\\
\bottomrule
\end{tabular}
\centering
\end{table}

We also examined the impact of different approximation methods for bilinear terms, i.e., constraints~\eqref{mcCormickConstraints} and~\eqref{GupteMcCormick}; for the latter, we considered two different discretization intervals for the white blood cell count, i.e., $\Delta = \frac{1}{20}(\epop_{w,0} - \rhs_{w})$ and $\Delta = \frac{1}{40}(\epop_{w,0} - \rhs_{w})$. Recall that $\epop_{w,0}$ and $\rhs_{w}$ denote upper and lower bounds on the white blood cell count, respectively. The results are presented in Table~\ref{tab:byMC_Deterministic}, in which ``Continuous'' refers to the McCormick relaxation method, i.e., the constraint set~\eqref{mcCormickConstraints}, and ``Discrete'' to the (modified) method of~\cite{Gupte2013}, i.e., the constraint set~\eqref{GupteMcCormick}. As expected, the McCormick relaxation method leads to a lower optimal objective value, due to its flexibility. This method, however, does not provide a means to control the approximation quality of the bilinear terms. The solvability of the model decreases as more integer variables are introduced through refining the discretization interval in the method of~\cite{Gupte2013}.   
\begin{table}
\centering
\caption{\label{tab:byMC_Deterministic} Solver statistics for different bilinearity approximation methods}
\begin{tabular}{lcccccccc}
\toprule
Method & $\Delta$ & Obj Val & Run-time\,(s) & Cons  & \phantom{\large -h}Vars\phantom{\huge $._g$} & IVars &BVars &Gap\\
\hline
Continuous & -- & 68.01 & 387 & 12579 & 7350 & 651 & 525 & $<$ 0.01\%\\
Discrete & ${1}/{20}$ &  68.13 & 2029 & 17619 & 9051 & 1092 & 966 & $<$ 0.01\%\\
Discrete & ${1}/{40}$ & 68.14 & 7200 & 22659 & 10731 & 1512 & 1386 & \quad 0.13\%\\
\bottomrule
\end{tabular}
\centering
\end{table}

Next, we present the results of our numerical study with the chance-constrained optimization model~\eqref{eq:stoMILP}. We simulated a branching process to generate a set of scenarios describing the heterogeneity of tumor. Details of the branching process and scenario generation are provided in Appendix~\ref{app:branching} of the e-companion.  Table~\ref{tab:simulatedScenarioGeneration} displays the (logarithm of) cancer cell populations and realization probability for each scenario. Based on the branching process, there is a single dominant scenario (Scenario 1) with the associated probability of over 0.77, in which the non-resistant cell type has the largest cell count. The probability of a scenario in which a drug-resistant cell type has the largest cell count is less than 0.02. Although the size of the drug-resistant part of a tumor is expected to be smaller that the non-resistant part, it can become dominant without treatment. 
\begin{table}
\centering
\caption{\label{tab:simulatedScenarioGeneration} Simulated scenarios generated by a branching process}
\begin{tabular}{cccccc}
\toprule
Scenario & Non-resist.(0) & Capec.-resist.(1) & Docet.-resist.(2) & Etopo.-resist.(3) & \phantom{\large -h}Prob.\phantom{\huge $._g$}\\
\hline
1 & 20.53 & 17.89 & 17.89 & 17.89 & 0.7705\\
2 & 20.44 & 17.85 & 17.83 & 18.69 & 0.0619\\
3 & 20.44 & 18.70 & 17.86 & 17.86 & 0.0603\\
4 & 20.44 & 17.84 & 18.74 & 17.82 & 0.0579\\
5 & 20.22 & 19.50 & 17.74 & 17.72 & 0.0109\\
6 & 20.23 & 17.71 & 19.49 & 17.71 & 0.0109\\
7 & 20.25 & 17.72 & 17.67 & 19.46 & 0.0103\\
8 & 19.80 & 17.35 & 17.45 & 20.09 & 0.0064\\
9 & 19.81 & 17.20 & 20.10 & 17.28 & 0.0059\\
10 & 19.80 & 20.11 & 17.18 & 17.39 & 0.0050\\
\bottomrule
\end{tabular}
\centering
\end{table}

\begin{figure}[t]
\centering
\begin{subfloat}{
\centering	
\includegraphics[width = .45\textwidth]{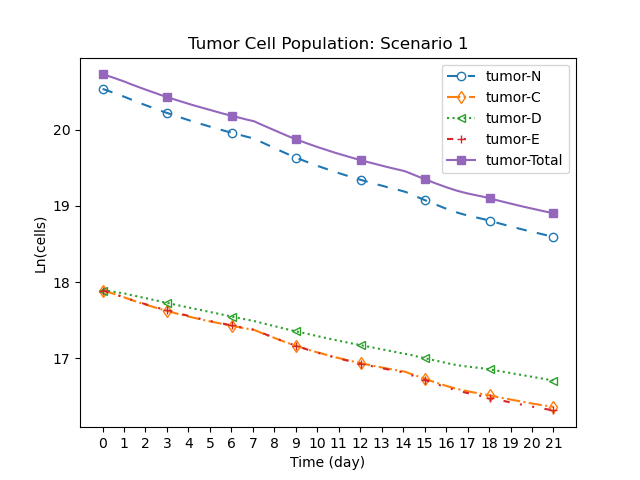}}
\end{subfloat}
\begin{subfloat}{
\centering
\includegraphics[width = .45\textwidth]{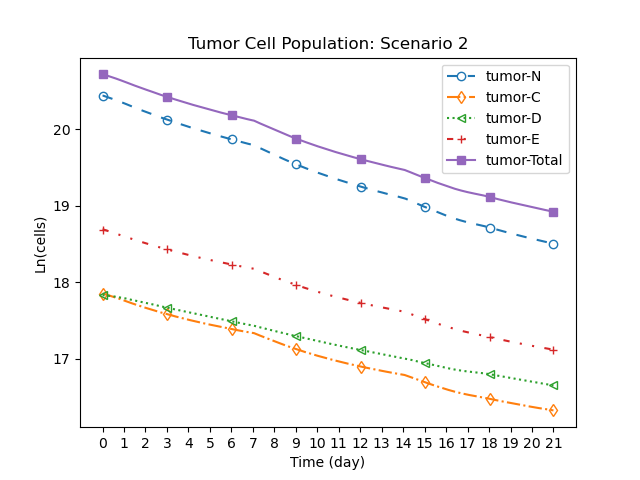}}
\end{subfloat}
\begin{subfloat}{
\centering	
\includegraphics[width = .45\textwidth]{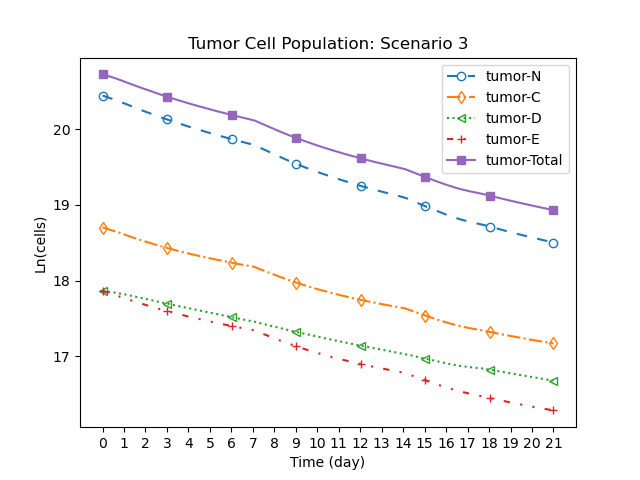}}
\end{subfloat}
\begin{subfloat}{
\centering
\includegraphics[width = .45\textwidth]{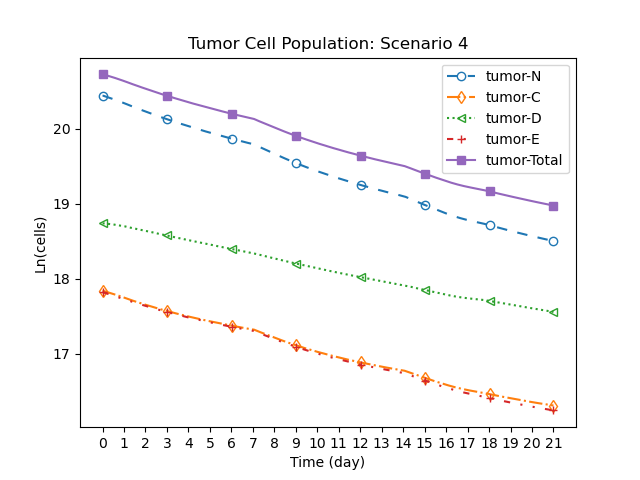}}
\end{subfloat}
\caption{Treatment effect on tumor cell populations under Scenarios 1--4, given by the neoadjuvant chance-constrained model~\eqref{eq:stoMILP}}
\label{fig:sto}
\end{figure}

The initial cancer cell population in all scenarios is $10^9$ cells, which translates to a tumor with a diameter of about 25 mm~\citep{del2009does}. While the operable size of a tumor must be determined clinically based on the tumor location and patient's health conditions, the diameter of 20 mm is commonly considered the border of stage II and stage III breast cancer~\citep{narod2013two,Senkus2015}. Thus, in our study, we assume $\epop_{\text{surg}} = 0.4\bigcdot10^{9}$ cells, which translates to a diameter of less than 20 mm. We also set $\epsilon = 0.05$ indicating that the desired probability of reaching an operable tumor size at the end of treatment period is at least 0.95. Using the discrete approximation of the bilinear terms with $\Delta = \frac{1}{20}(\epop_{w,0} - \rhs_{w})$ and setting the time-step $h$ equal to one hour, the neoadjuvant chance-constrained optimization model~\eqref{eq:stoMILP} contained 27,205 variables (1,102 integer, 976 binary) and 35,804 constraints; the solver found an optimal solution in about 17 minutes with the optimal objective value of 67.99. Figure~\ref{fig:sto} illustrates the treatment effect on cancer cell populations under Scenarios 1 to 4; these are the scenarios with the realization probability of at least 0.05. More details on the output of the chanced-constrained model~\eqref{eq:stoMILP} as well as its counterpart with a probability-based objective, i.e., minimize $\epsilon$ subject to~\eqref{eq:sto_PK}--\eqref{eq:sto_last}, are provided in the e-companion; see Appendix~\ref{app:figs}. As stated earlier and illustrated in Figures~\ref{fig:prob_optAdmin} and~\ref{fig:pro_remaining}, the probability-based objective leads to inferior tumor-shrinkage compared with formulation~\eqref{eq:stoMILP}.
\subsection{Sensitivity Analysis}
\label{sec:sensitivity}
The proposed optimization models involve several parameters that need to be estimated based on clinical data and a patient's biological characteristics. In this section, we present the results of our analysis to determine the sensitivity of an optimal solution and objective value of the combination chemotherapy optimization problem~\eqref{eq:detMILP} to these parameters.

\begin{figure}
\centering
\begin{subfloat}[Elimination rate]{
\centering	
\includegraphics[width = .45\textwidth]{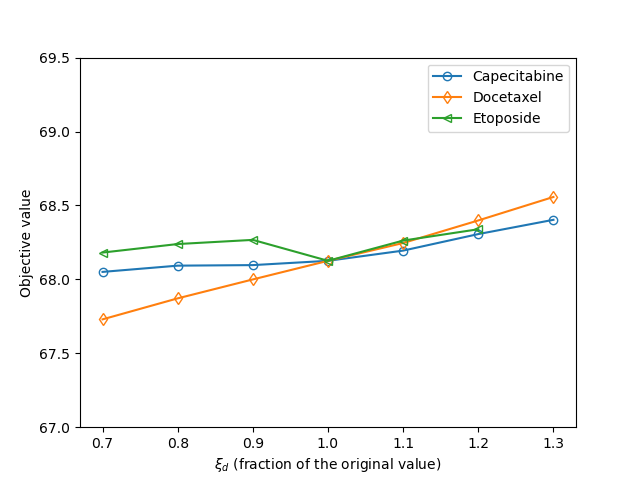}
\label{fig:sen_elimination}}
\end{subfloat}
\begin{subfloat}[Kill effect on cancer cells]{
\centering	
\includegraphics[width = .45\textwidth]{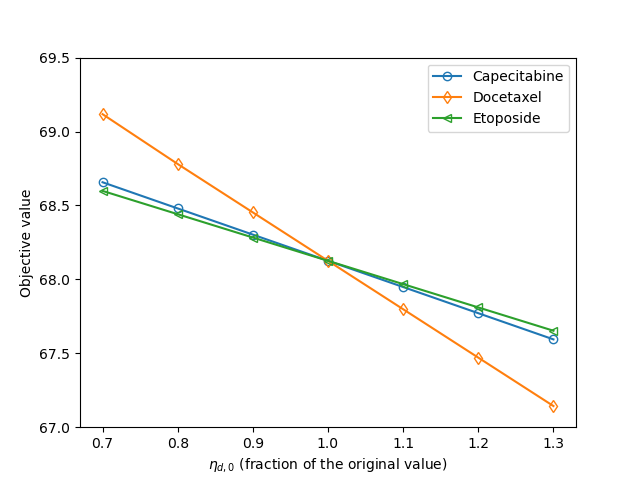}
\label{fig:sen_kill_cancer}}
\end{subfloat}
\begin{subfloat}[Kill effect on white blood cells]{
\centering
\includegraphics[width = .45\textwidth]{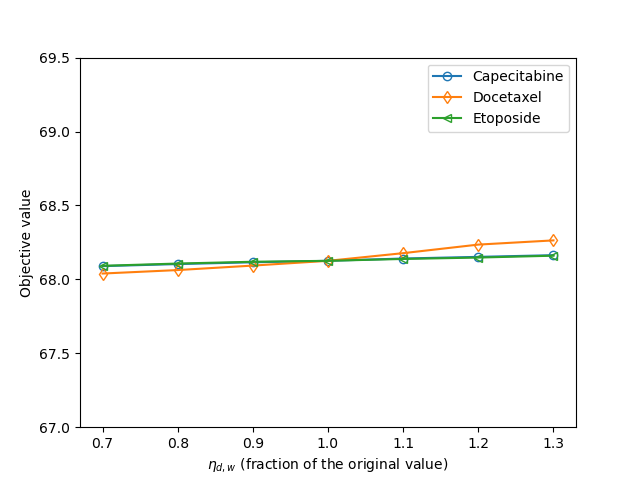}
\label{fig:sen_kill_WBC}}
\end{subfloat}
\begin{subfloat}[Temporal resistance]{
\centering
\includegraphics[width = .45\textwidth]{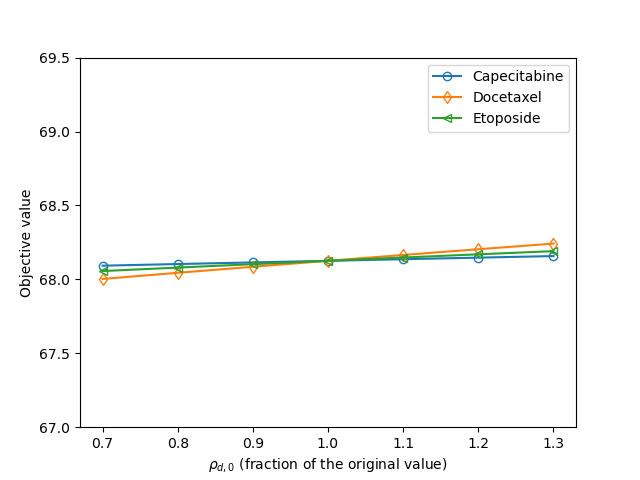}
\label{fig:sen_temp_res}}
\end{subfloat}
\caption{Sensitivity analysis results for the deterministic model~\eqref{eq:detMILP} with respect to pharmacokinetics and pharmacodynamics parameters}
\label{fig:sen_Pk_PD}
\end{figure}

Figure~\ref{fig:sen_Pk_PD} illustrates the results of the sensitivity analysis with respect to the pharmacokinetics and pharmacodynamics parameters. We used clinical data to estimate the fractional kill effect parameter of each drug for the non-resistant cell type, i.e., $\dCons_{d, 0}$, and set $\dCons_{d, q} = 0.25 \, \dCons_{d, 0},\forAll q \in \{1,2,3\}$, to account for drug-resistance in our numerical study; see Appendix~\ref{app:estimate} for more information. We also assumed that temporal resistance is constant across the cancer cell types, i.e., $\dExp_{d, q} = \dExp_{d, 0},\forAll q \in \{1,2,3\}$. Thus, for these parameters, we focused on $\dCons_{d, 0}$ and $\dExp_{d, 0}$. For each drug, we varied the corresponding parameters in 10\% increments, while the other parameters of the model remained constant, and measured the corresponding impact on the optimal objective value. To provide a convenient comparison, in Figure~\ref{fig:sen_Pk_PD}, the horizontal axes demonstrate (changed) parameter values as a fraction of the original value; the vertical axes display the optimal objective value, i.e., $\sum\limits_{q \in \typeSet}\pop^*_{q,\lastTym}$. 

In general, the optimal objective value of the combination chemotherapy optimization problem~\eqref{eq:detMILP} is more sensitive to the fractional kill effect of a drug on cancer cells than any other pharmacokinetics or pharmacodynamics parameter. Based on Figure~\ref{fig:sen_kill_cancer}, docetaxel, the intravenous drug, is the most influential on the optimal objective value among the considered drugs; the oral drugs, capecitabine and etoposide, are less impactful and show very similar patterns. Apart from the fractional kill effect on cancer cells, the objective is most sensitive to the drug elimination rate parameters. Recall that, for a drug $d$, the elimination rate $\xi_d$ determines how fast the drug concentration in the body declines. Figure~\ref{fig:sen_elimination} suggests docetaxel has the most influential elimination rate on the objective value as well. It is important to note that docetaxel has a mandated one-week rest period; hence, the model has very little flexibility for administration of this drug. In other words, as opposed to oral drugs that can be administrated frequently to keep their concentrations constantly high during the treatment period, the docetaxel concentration is much lower than its maximum level for most of the treatment period. Due to this operational constraint, changes to the elimination rate of this drug directly influence its concentration profile and lead to relatively high impacts on the optimal objective value. As shown in Figures~\ref{fig:sen_kill_WBC} and~\ref{fig:sen_temp_res}, the model shows much less sensitivity to the parameters representing fractional kill effect on white blood cells and temporal resistance, i.e., $\dCons_{d, w}$ and $\dExp_{d, 0}$. The low impact of $\dCons_{d, w}$ comparative to $\dCons_{d, 0}$ is justified by the fact that the model has a tumor shrinkage-based objective.

Another important observation regarding the impact of operational constraints on the optimal objective value of model~\eqref{eq:detMILP} concerns the elimination rate of etoposide; see Figure~\ref{fig:sen_elimination}. The increasing patterns observed from 0.7 to 0.9 (as fractions of the original value) and from 1.0 to 1.2 are justified by the fact that higher elimination rates lead to lower drug concentrations and smaller effects on tumor cells. However, a better objective value is obtained when the elimination rate changes from 0.9 to 1.0. This is due to the fact that the optimal administration regimen for etoposide changes in this interval. In fact, a higher elimination rate allows the model to administer etoposide more frequently without violation of the maximum permissible concentration, given the discrete administration times. A similar phenomenon underlies the change of slope observed in this figure for capecitabine. Finally, we note that the value of the elimination rate parameter of a drug is physically restricted to $\xi = 1.0~\text{day}^{-1}$; hence, the etoposide curve in Figure~\ref{fig:sen_elimination} terminates at $\xi = 1.2 \times 0.8 = 0.96~\text{day}^{-1}$.

\begin{figure}
\centering
\begin{subfloat}[Maximum permissible dose (and concentration)]{
\centering
\includegraphics[width = .45\textwidth]{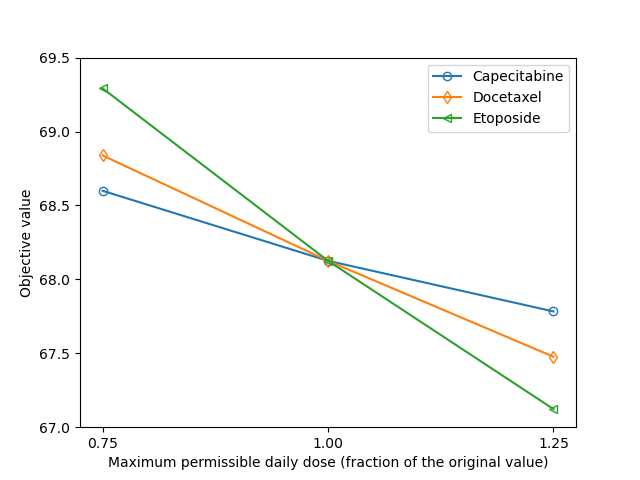}
\label{fig:sen_dose}}
\end{subfloat}
\begin{subfloat}[Neutropenia threshold]{
\centering
\includegraphics[width = .45\textwidth]{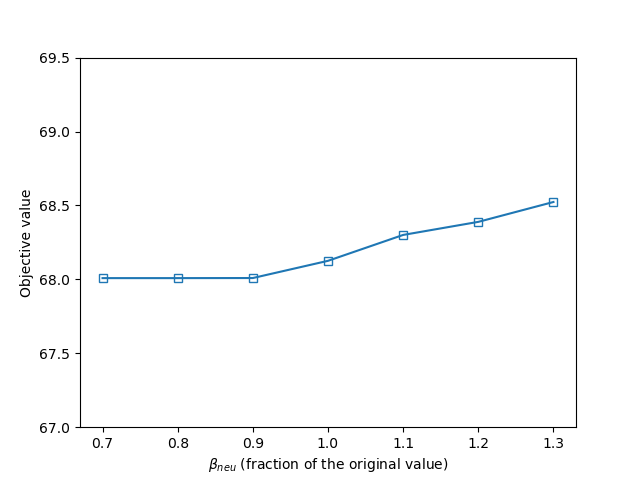}
\label{fig:sen_Neu}}
\end{subfloat}
\caption{Sensitivity analysis results for the deterministic model~\eqref{eq:detMILP} with respect to operational parameters}
\label{fig:sen_opr}
\end{figure}

Aside from the pharmacokinetics and pharmacodynamics parameters, we evaluated the impact of changes to operational parameters governing the combination chemotherapy optimization model~\eqref{eq:detMILP}. In particular, we investigated changes to the maximum permissible dose (and concentration) for each drug and the neutropenia threshold. Recall that the maximum permissible doses (and concentrations) used in our numerical study are based on common administration regimens in clinical trials. For capecitabine, the clinical administration dose translates to 8 pills per day, restricted to 4 pills per meal; we denote this regimen by 8/4 hereafter. We considered increasing these limits to 10/5 and decreasing them to 6/3. Based on the new values, we calculated the maximum permissible doses and drug concentration through simulation of the corresponding clinical trial. Similarly, for etoposide, we considered changing the current administration regimen of 2/1 to 3/2 and 1/1. In the presentation of the result to follow,  with a slight abuse of the term dose, we demonstrate these changes as 25\% increase and decrease to the maximum permissible dose. For the intravenous drug, docetaxel, we directly applied the 25\% increase and decrease to the current maximum permissible daily dose and infusion rate; similar to oral drugs, we obtained the corresponding maximum permissible concentration from simulation of the corresponding clinical trial with the new administration values. Figure~\ref{fig:sen_dose} illustrates the results. According to this figure, the maximum permissible dose (and concentration) constraints are the most restrictive for etoposide, and their relative relaxation can significantly impact the treatment outcome. Note that the new optimal solutions satisfy the neutropenia constraint; under higher administration doses, the count of neutrophils never falls below the neutropenia threshold, but it reaches this level much earlier and stays there for the rest of the treatment period. Such an increased dose may lead to side effects other than neutropenia, which are not captured by our model. 

Normal white blood cell counts vary across patients; different individuals may have different white blood cell thresholds for the associated side effects. In our sensitivity analysis, we changed the neutropenia threshold in 10\% increments and  measured its impact on the optimal objective value; the results are illustrated in Figure~\ref{fig:sen_Neu}. Given the tightness of the neutropenia constraint, as mentioned before, deterioration of the optimal objective value due to the increase of neutropenia threshold is well justified. This figure also shows that decreasing the neutropenia threshold below 90\% of the current level will not affect the optimal objective value, given the specified values for other parameters in our model. This observation demonstrates that, below this level, other operational constraints, i.e., maximum permissible doses and concentrations, are the driving factors. 
\section{Conclusion}
\label{sec:conc}
In this paper, we present a mixed-integer linear programming model for combination chemotherapy optimization, which seeks to find optimal administration dose and schedule for cytotoxic drugs by minimizing cancer cell population at the end of a treatment period. As opposed to previous works that often ignore operational considerations or low white blood cell counts as a toxic effect, we incorporate these constraints. We also extend this model to account for the uncertainty of tumor heterogeneity and present a chanced-constrained model for neoadjuvant chemotherapy. We use the literature and published clinical data to calibrate our model parameters for a case of breast cancer and present the results of our numerical study. We perform sensitivity analyses to identify the most influential parameters on the model outcomes. Our models provide a framework for the exploration of new, individualized dose guidelines. Future directions stemming from this work include improving estimates of model parameters, considering other drugs or types of cancer in the numerical study, and factoring additional toxicities. 
\ACKNOWLEDGMENT{%
The authors would like to thank Dr. Jeffrey Myers 
and Cem Dede of The University of Texas MD Anderson Cancer Center, Adam Palmer of the University of North Carolina, and David Mildebrath, Soheil Hemmati, Saumya Sinha, and M. Can Camur of Rice University for their helpful comments. This research was supported in part by National Science Foundation grants CMMI-1933369 and CMMI-1933373.
}


\begin{thebibliography}{90}
\providecommand{\natexlab}[1]{#1}
\providecommand{\url}[1]{\texttt{#1}}
\providecommand{\urlprefix}{URL }

\bibitem[{Ab{\'e}cassis et~al.(2019)Ab{\'e}cassis, Hamy, Laurent
  et~al.}]{abecassis2019}
Ab{\'e}cassis J, Hamy A, Laurent C, et~al. (2019) Assessing reliability of
  intratumor heterogeneity estimates from single sample whole exome sequencing
  data. \emph{PLoS ONE} 14(11):e0224143.

\bibitem[{Al-Khayyal \protect\BIBand{} Falk(1983)}]{al1983jointly}
Al-Khayyal F, Falk J (1983) Jointly constrained biconvex programming.
  \emph{Mathematics of Operations Research} 8(2):273--286.

\bibitem[{Alam et~al.(2013)Alam, Hossain, Algoul et~al.}]{alam2013}
Alam M, Hossain M, Algoul S, et~al. (2013) Multi-objective multi-drug
  scheduling schemes for cell cycle specific cancer treatment. \emph{Computers
  \& Chemical Engineering} 58:14--32.

\bibitem[{American {C}ancer {S}ociety({\natexlab{a}})}]{CancerFacts2021}
American {C}ancer {S}ociety (2021{\natexlab{a}}) Cancer {F}acts \& {F}igures
  2021.
  \url{https://www.cancer.org/research/cancer-facts-statistics/all-cancer-facts-figures/cancer-facts-figures-2021.html}.

\bibitem[{American {C}ancer {S}ociety({\natexlab{b}})}]{CancerSurvivorship2021}
American {C}ancer {S}ociety (2021{\natexlab{b}}) Cancer {T}reatment \&
  {S}urvivorship {F}acts \& {F}igures.
  \url{https://www.cancer.org/research/cancer-facts-statistics/survivor-facts-figures.html}.

\bibitem[{American {C}ancer {S}ociety({\natexlab{c}})}]{sideEffects}
American {C}ancer {S}ociety (2021{\natexlab{c}}) {C}hemotherapy {S}ide
  {E}ffects.
  \url{https://www.cancer.org/treatment/treatments-and-side-effects/treatment-types/chemotherapy/chemotherapy-side-effects.html}.

\bibitem[{American {C}ancer {S}ociety({\natexlab{d}})}]{howChemo2021}
American {C}ancer {S}ociety (2021{\natexlab{d}}) How {I}s {C}hemotherapy {U}sed
  to {T}reat {C}ancer?
  \url{https://www.cancer.org/treatment/treatments-and-side-effects/treatment-types/chemotherapy/how-is-chemotherapy-used-to-treat-cancer.html}.

\bibitem[{Asachenkov et~al.(1994)Asachenkov, Marchuk, Mohler, \protect\BIBand{}
  Zuev}]{Asachenkov1994}
Asachenkov A, Marchuk G, Mohler R, Zuev S (1994) \emph{Disease Dynamics}
  (Basel, Basel-Stadt, Switzerland: Birkh\"{a}user Basel), 1st edition.

\bibitem[{Baker et~al.(2006)Baker, Sparreboom, \protect\BIBand{}
  Verweij}]{Baker2006}
Baker S, Sparreboom A, Verweij J (2006) Clinical pharmacokinetics of docetaxel.
  \emph{Clinical Pharmacokinetics} 45(3):235--252.

\bibitem[{Beraldi \protect\BIBand{} Bruni(2014)}]{Beraldi2014}
Beraldi P, Bruni M (2014) A clustering approach for scenario tree reduction:
  {A}n application to a stochastic programming portfolio optimization problem.
  \emph{TOP} 22:934--949.

\bibitem[{Bonate(2011)}]{bonate2011}
Bonate P (2011) \emph{Pharmacokinetic-Pharmacodynamic Modeling and Simulation},
  volume~20 (Springer).

\bibitem[{Butcher(2007)}]{butcher2007runge}
Butcher J (2007) Runge-{K}utta methods. \emph{Scholarpedia} 2(9):3147.

\bibitem[{Butcher(2008)}]{Butcher2008}
Butcher J (2008) \emph{Numerical Methods for Ordinary Differential Equations}
  (John Wiley \& Sons Ltd.), 2nd edition.

\bibitem[{Cajal et~al.(2020)Cajal, Ses{\'e}, Capdevila et~al.}]{y2020}
Cajal S, Ses{\'e} M, Capdevila C, et~al. (2020) Clinical implications of
  intratumor heterogeneity: Challenges and opportunities. \emph{Journal of
  Molecular Medicine} 98(2):161--177.

\bibitem[{Cameron(1997)}]{Cameron1997}
Cameron D (1997) Mathematical modelling of the response of breast cancer to
  drug therapy. \emph{Journal of Theoretical Medicine} 2:137--151.

\bibitem[{Chan et~al.(1999)Chan, Friedrichs, Noel et~al.}]{chan1999}
Chan S, Friedrichs K, Noel D, et~al. (1999) Prospective randomized trial of
  docetaxel versus doxorubicin in patients with metastatic breast cancer.
  \emph{Journal of Clinical Oncology} 17(8):2341--2341.

\bibitem[{Coldman \protect\BIBand{} Goldie(1983)}]{Coldman1983}
Coldman A, Goldie J (1983) A model for the resistance of tumor cells to cancer
  chemotherapeutic agents. \emph{Mathematical Biosciences} 65(2):291--307.

\bibitem[{Coldman \protect\BIBand{} Murray(2000)}]{coldman2000}
Coldman A, Murray J (2000) Optimal control for a stochastic model of cancer
  chemotherapy. \emph{Mathematical Biosciences} 168(2):187--200.

\bibitem[{Costa \protect\BIBand{} Boldrini(1997)}]{costa1997}
Costa M, Boldrini J (1997) Chemotherapeutic treatments: A study of the
  interplay among drug resistance, toxicity and recuperation from side effects.
  \emph{Bulletin of Mathematical Biology} 59(2):205--232.

\bibitem[{Day(1986)}]{Day1986TreatmentChemotherapy}
Day R (1986) Treatment sequencing, asymmetry, and uncertainty: {P}rotocol
  strategies for combination chemotherapy. \emph{Cancer Research}
  46(8):3876--3885.

\bibitem[{{de}~Pillis et~al.(2007){de}~Pillis, Gu, Fister
  et~al.}]{dePillis2007}
{de}~Pillis L, Gu W, Fister K, et~al. (2007) Chemotherapy for tumors: An
  analysis of the dynamics and a study of quadratic and linear optimal
  controls. \emph{Mathematical Biosciences} 209(1):292--315.

\bibitem[{Del~Monte(2009)}]{del2009does}
Del~Monte U (2009) Does the cell number 10$^9$ still really fit one gram of
  tumor tissue? \emph{Cell Cycle} 8(3):505--506.

\bibitem[{d'Onofrio et~al.(2009)d'Onofrio, Ledzewicz, Maurer, \protect\BIBand{}
  Sch{\"{a}}ttler}]{dOnofrio2009OnTumors}
d'Onofrio A, Ledzewicz U, Maurer H, Sch{\"{a}}ttler H (2009) {On optimal
  delivery of combination therapy for tumors}. \emph{Mathematical Biosciences}
  222(1):13--26.

\bibitem[{Ebata et~al.(2018)Ebata, Hirano, Konishi
  et~al.}]{ebata2018randomized}
Ebata T, Hirano S, Konishi M, et~al. (2018) Randomized clinical trial of
  adjuvant gemcitabine chemotherapy versus observation in resected bile duct
  cancer. \emph{Journal of British Surgery} 105(3):192--202.

\bibitem[{Ershler(2006)}]{Ershler2006}
Ershler WB (2006) Capecitabine monotherapy: Safe and effective treatment for
  metastatic breast cancer. \emph{The Oncologist} 11(4):325--335.

\bibitem[{Floares et~al.(2003)Floares, Floares, Cucu, \protect\BIBand{}
  Lazar}]{Floares2003}
Floares A, Floares C, Cucu M, Lazar L (2003) Adaptive neural networks control
  of drug dosage regimens in cancer chemotherapy. \emph{Proceedings of the
  International Joint Conference on Neural Networks}, volume~1, 154--159
  (IEEE).

\bibitem[{Frances et~al.(2011)Frances, Claret, Bruno, \protect\BIBand{}
  Iliadis}]{Frances2011}
Frances N, Claret L, Bruno R, Iliadis A (2011) Tumor growth modeling from
  clinical trials reveals synergistic anticancer effect of the capecitabine and
  docetaxel combination in metastatic breast cancer. \emph{Cancer Chemotherapy
  and Pharmacology} 68(6):1413--1419.

\bibitem[{Gerlinger et~al.(2012)Gerlinger, Rowan, Horswell
  et~al.}]{gerlinger2012}
Gerlinger M, Rowan A, Horswell S, et~al. (2012) Intratumor heterogeneity and
  branched evolution revealed by multiregion sequencing. \emph{The New England
  Journal of Medicine} 366:883--892.

\bibitem[{G\"{u}lpinar et~al.(2004)G\"{u}lpinar, Rustem, \protect\BIBand{}
  Settergren}]{Gulpinar2004}
G\"{u}lpinar N, Rustem B, Settergren R (2004) Simulation and optimization
  approaches to scenario tree generation. \emph{Journal of Economic Dynamics
  and Control} 28(7):1291--1315.

\bibitem[{Gupte et~al.(2013)Gupte, Ahmed, Cheon, \protect\BIBand{}
  Dey}]{Gupte2013}
Gupte A, Ahmed S, Cheon M, Dey S (2013) Solving mixed integer bilinear problems
  using {MILP} formulations. \emph{SIAM Journal on Optimization}
  23(2):721--744.

\bibitem[{Hande(1998)}]{Hande1998}
Hande K (1998) Etoposide: {F}our decades of development of a topoisomerase {II}
  inhibitor. \emph{European Journal of Cancer} 34(10):1514--1521.

\bibitem[{Harrold \protect\BIBand{} Parker(2009)}]{Harrold2009}
Harrold J, Parker R (2009) {Clinically relevant cancer chemotherapy dose
  scheduling via mixed-integer optimization}. \emph{Computers and Chemical
  Engineering} 33(12):2042--2054.

\bibitem[{Hu et~al.(2016)Hu, Sun, Wang, \protect\BIBand{} Gu}]{hu2016recent}
Hu Q, Sun W, Wang C, Gu Z (2016) Recent advances of cocktail chemotherapy by
  combination drug delivery systems. \emph{Advanced Drug Delivery Reviews}
  98:19--34.

\bibitem[{Hu et~al.(2017)Hu, Sun, \protect\BIBand{} Curtis}]{hu2017}
Hu Z, Sun R, Curtis C (2017) A population genetics perspective on the
  determinants of intratumor heterogeneity. \emph{Biochimica et Biophysica Acta
  (BBA)-Reviews on Cancer} 1867(2):109--126.

\bibitem[{Iliadis \protect\BIBand{} Barbolosi(2000)}]{Iliadis2000}
Iliadis A, Barbolosi D (2000) Optimizing drug regimens in cancer chemotherapy
  by an efficacy-toxicity mathematical model. \emph{Computers and Biomedical
  Research} 33:211--226.

\bibitem[{Itik et~al.(2009)Itik, Salamci, \protect\BIBand{} Banks}]{itik2009}
Itik M, Salamci M, Banks S (2009) Optimal control of drug therapy in cancer
  treatment. \emph{Nonlinear Analysis: Theory, Methods \& Applications}
  71(12):e1473--e1486.

\bibitem[{Jacqmin et~al.(2007)Jacqmin, Snoeck, Van~Schaick
  et~al.}]{jacqmin2007}
Jacqmin P, Snoeck E, Van~Schaick E, et~al. (2007) Modelling response time
  profiles in the absence of drug concentrations: {D}efinition and performance
  evaluation of the {K}--{PD} model. \emph{Journal of Pharmacokinetics and
  Pharmacodynamics} 34(1):57--85.

\bibitem[{Kasi \protect\BIBand{} Grothey(2018)}]{kasi2018}
Kasi P, Grothey A (2018) Chemotherapy-induced neutropenia as a prognostic and
  predictive marker of outcomes in solid-tumor patients. \emph{Drugs}
  78(7):737--745.

\bibitem[{Kimmel \protect\BIBand{} Axelrod(2015)}]{Kimmel2015}
Kimmel M, Axelrod D (2015) \emph{Branching Processes in Biology}
  (Springer-Verlag), 2nd edition.

\bibitem[{Kosaka et~al.(2015)Kosaka, Rai, Masuda et~al.}]{Kosaka2015}
Kosaka Y, Rai Y, Masuda N, et~al. (2015) Phase {III} placebo-controlled,
  double-blind, randomized trial of pegfilgrastim to reduce the risk of febrile
  neutropenia in breast cancer patients receiving docetaxel/dyclophosphamide
  chemotherapy. \emph{Supportive Care in Cancer} 23:1137--1145.

\bibitem[{Laird(1964)}]{Laird1964}
Laird A (1964) Dynamics of tumour growth. \emph{British Journal of Cancer}
  18(3):490--502.

\bibitem[{Laird et~al.(1965)Laird, Tyler, Barton et~al.}]{laird1965dynamics}
Laird A, Tyler SA, Barton A, et~al. (1965) Dynamics of normal growth.
  \emph{Growth} 29:233--248.

\bibitem[{LeVeque(2007)}]{LeVeque2007}
LeVeque R (2007) \emph{Finite Difference Methods for Ordinary and Partial
  Differential Equations} (Philadelphia: Society for Industrial and Applied
  Mathematics).

\bibitem[{Liang et~al.(2006)Liang, Leung, \protect\BIBand{} Mok}]{Liang2006}
Liang Y, Leung K, Mok T (2006) A novel evolutionary drug scheduling model in
  cancer chemotherapy. \emph{IEEE Transactions on Information Technology in
  Biomedicine} 10(2):237--245.

\bibitem[{Luqmani(2005)}]{Luqmani2005}
Luqmani Y (2005) Mechanisms of drug resistance in cancer chemotherapy.
  \emph{Medical Principles and Practice} 14(suppl.1):35--48.

\bibitem[{Mariotti et~al.(2021)Mariotti, Han, Ismail-Khan
  et~al.}]{mariotti2021effect}
Mariotti V, Han H, Ismail-Khan R, et~al. (2021) Effect of taxane chemotherapy
  with or without indoximod in metastatic breast cancer: A randomized clinical
  trial. \emph{JAMA Oncology} 7(1):61--69.

\bibitem[{Martin(1992)}]{martin1992}
Martin R (1992) Optimal control drug scheduling of cancer chemotherapy.
  \emph{Automatica} 28(6):1113--1123.

\bibitem[{Martin et~al.(1990)Martin, Fisher, Minchin, \protect\BIBand{}
  Teo}]{martin1990}
Martin R, Fisher M, Minchin R, Teo K (1990) A mathematical model of cancer
  chemotherapy with an optimal selection of parameters. \emph{Mathematical
  Biosciences} 99(2):205--230.

\bibitem[{Martin et~al.(1992{\natexlab{a}})Martin, Fisher, Minchin,
  \protect\BIBand{} Teo}]{Martin1992Low}
Martin R, Fisher M, Minchin R, Teo K (1992{\natexlab{a}}) Low-intensity
  combination chemotherapy maximizes host survival time for tumors containing
  drug-resistant cells. \emph{Mathematical Biosciences} 110(2):221--252.

\bibitem[{Martin et~al.(1992{\natexlab{b}})Martin, Fisher, Minchin,
  \protect\BIBand{} Teo}]{Martin1992Resistant}
Martin R, Fisher M, Minchin R, Teo K (1992{\natexlab{b}}) Optimal control of
  tumor size used to maximize survival time when cells are resistant to
  chemotherapy. \emph{Mathematical Biosciences} 110(2):201--219.

\bibitem[{Martin \protect\BIBand{} Teo(1994)}]{Martin1994}
Martin R, Teo K (1994) \emph{Optimal Control of Drug Administration in Cancer
  Chemotherapy} (World Scientific).

\bibitem[{McCormick(1976)}]{McCormick1976}
McCormick G (1976) Computability of global solutions to factorable nonconvex
  programs: {P}art {I} -- {C}onvex underestimating problems. \emph{Mathematical
  Programming} 10:147--175.

\bibitem[{Mitrovic et~al.(2012)Mitrovic, Perry, Suzumiya
  et~al.}]{mitrovic2012prognostic}
Mitrovic Z, Perry A, Suzumiya J, et~al. (2012) The prognostic significance of
  lymphopenia in peripheral {T}-cell and natural killer/{T}-cell lymphomas: A
  study of 826 cases from the {I}nternational {P}eripheral {T}-cell {L}ymphoma
  {P}roject. \emph{American Journal of Hematology} 87(8):790--794.

\bibitem[{Murray(1990)}]{Murray1990}
Murray J (1990) Some optimal control problems in cancer chemotherapy with a
  toxicity limit. \emph{Mathematical Biosciences} 100(1):49--67.

\bibitem[{Murray(1994)}]{Murray1994}
Murray J (1994) Optimal drug regimens in cancer chemotherapy for single drugs
  that block progression through the cell cycle. \emph{Mathematical
  Biosciences} 123(2):183--213.

\bibitem[{Murray(1997)}]{murray1997}
Murray J (1997) The optimal scheduling of two drugs with simple resistance for
  a problem in cancer chemotherapy. \emph{Mathematical Medicine and Biology: A
  Journal of the IMA} 14(4):283--303.

\bibitem[{Nanda et~al.(2007)Nanda, Moore, \protect\BIBand{}
  Lenhart}]{nanda2007}
Nanda S, Moore H, Lenhart S (2007) Optimal control of treatment in a
  mathematical model of chronic myelogenous leukemia. \emph{Mathematical
  Biosciences} 210(1):143--156.

\bibitem[{Narod et~al.(2013)Narod, Iqbal, Jakubowska et~al.}]{narod2013two}
Narod S, Iqbal J, Jakubowska A, et~al. (2013) Are two-centimeter breast cancers
  large or small? \emph{Current Oncology} 20(4):205--211.

\bibitem[{Norton(1988)}]{Norton1988}
Norton L (1988) A {G}ompertzian model of human breast cancer growth.
  \emph{Cancer Research} 48(24 Part 1):7067--7071.

\bibitem[{O’{S}haughnessy et~al.(2001)O’{S}haughnessy, Blum, Moiseyenko
  et~al.}]{o2001randomized}
O’{S}haughnessy J, Blum J, Moiseyenko V, et~al. (2001) Randomized,
  open-label, {P}hase {II} trial of oral capecitabine
  ({X}eloda{\textregistered}) vs. a reference arm of intravenous {CMF}
  (cyclophosphamide, cethotrexate and 5-fluorouracil) as first-line therapy for
  advanced/metastatic breast cancer. \emph{Annals of Oncology}
  12(9):1247--1254.

\bibitem[{Palmeri et~al.(2008)Palmeri, Vaglica, \protect\BIBand{}
  Palmeri}]{Palmeri2008}
Palmeri L, Vaglica M, Palmeri S (2008) Weekly docetaxel in the treatment of
  metastatic breast cancer. \emph{Therapeutics and Clinical Risk Management}
  4(5):1047--1059.

\bibitem[{Panetta \protect\BIBand{} Adam(1995)}]{Panetta1995}
Panetta J, Adam J (1995) A mathematical model of cycle-specific chemotherapy.
  \emph{Mathematical and Computer Modelling} 22(2):67--82.

\bibitem[{Pereira et~al.(1995)Pereira, Pedreira, \protect\BIBand{}
  De~Sousa}]{pereira1995}
Pereira F, Pedreira C, De~Sousa J (1995) A new optimization based approach to
  experimental combination chemotherapy. \emph{Frontiers of Medical and
  Biological Engineering} 64(4):257--268.

\bibitem[{Petrovski et~al.(2004)Petrovski, Sudha, \protect\BIBand{}
  McCall}]{Petrovski2004}
Petrovski A, Sudha B, McCall J (2004) Optimising cancer chemotherapy using
  particle swarm optimisation and genetic algorithms. Yao X, Burke E, Lozano J,
  et~al., eds., \emph{Parallel Problem Solving from Nature - PPSN VIII},
  633--641 (Berlin, Heidelberg: Springer Berlin Heidelberg).

\bibitem[{Piraino et~al.(2019)Piraino, Thomas, O’Donovan, \protect\BIBand{}
  Furney}]{PIRAINO2019}
Piraino S, Thomas V, O’Donovan P, Furney S (2019) Mutations: Driver versus
  passenger. Boffetta P, Hainaut P, eds., \emph{Encyclopedia of Cancer (Third
  Edition)}, 551--562 (Oxford: Academic Press), third edition, ISBN
  978-0-12-812485-7.

\bibitem[{Pizzo(1993)}]{Pizzo1993}
Pizzo P (1993) Management of fever in patients with cancer and
  treatment-induced neutropenia. \emph{New England Journal of Medicine}
  328:1323--1332.

\bibitem[{Polyak(2011)}]{Polyak2011}
Polyak K (2011) Heterogeneity in breast cancer. \emph{The Journal of Clinical
  Investigation} 121(10):3786--3788.

\bibitem[{Reigner et~al.(2001)Reigner, Blesch, \protect\BIBand{}
  Weidekamm}]{Reigner2001}
Reigner B, Blesch K, Weidekamm E (2001) Clinical pharmacokinetics of
  capecitabine. \emph{Clinical Pharmacokinetics} 40(2):85--104.

\bibitem[{Rosado et~al.(2011)Rosado, Diamanti, Cascioli, Ceccarelli,
  Caporuscio, D'Amelio, Carsetti, \protect\BIBand{} Lagana}]{rosado2011hyper}
Rosado M, Diamanti A, Cascioli S, Ceccarelli S, Caporuscio S, D'Amelio R,
  Carsetti R, Lagana B (2011) Hyper-{I}g{M}, neutropenia, mild infections and
  low response to polyclonal stimulation: {H}yper-{I}g{M} syndrome or common
  variable immunodeficiency? \emph{International Journal of Immunopathology and
  Pharmacology} 24(4):983--991.

\bibitem[{Sager(2005)}]{sager2005numerical}
Sager S (2005) \emph{Numerical Methods for Mixed-Integer Optimal Control
  Problems} (Der Andere Verlag T{\"o}nning), ISBN 3-89959-416-9.

\bibitem[{Saville et~al.(2019)Saville, Smith, \protect\BIBand{}
  Bijak}]{saville2019}
Saville C, Smith H, Bijak K (2019) Operational research techniques applied
  throughout cancer care services: a review. \emph{Health Systems} 8(1):52--73.

\bibitem[{Segal et~al.(2014)Segal, Flood, Mancini et~al.}]{Segal2014}
Segal E, Flood M, Mancini R, et~al. (2014) Oral chemotherapy food and drug
  interactions: {A} comprehensive review of the literature. \emph{Journal of
  Oncology Practice} 10(4):e255--e268.

\bibitem[{Senkus et~al.(2015)Senkus, Kyriakides, Ohno et~al.}]{Senkus2015}
Senkus E, Kyriakides S, Ohno S, et~al. (2015) Primary breast cancer: {ESMO}
  clinical practice guidelines for diagnosis, treatment, and follow-up.
  \emph{Annals of Oncology} Supplement 5:v8--30.

\bibitem[{Sharma et~al.(2006)Sharma, Rivory, Beale et~al.}]{Sharma2006}
Sharma R, Rivory L, Beale P, et~al. (2006) A {P}hase {II} study of fixed-dose
  capecitabine and assessment of predictors of toxicity in patients with
  advanced/metastatic colorectal cancer. \emph{British Journal of Cancer}
  94:964--968.

\bibitem[{Shi et~al.(2014)Shi, Alagoz, Erenay, \protect\BIBand{} Su}]{Shi2014}
Shi J, Alagoz O, Erenay F, Su Q (2014) A survey of optimization models on
  cancer chemotherapy treatment planning. \emph{Annals of Operations Research}
  221(1):331--356.

\bibitem[{Skipper et~al.(1964)Skipper, Schabel~Jr., \protect\BIBand{}
  Wilcox}]{Skipper1964XIII}
Skipper H, Schabel~Jr F, Wilcox W (1964) Experimental evaluation of potential
  anticancer agents. {XIII}. {O}n the criteria and kinetics associated with
  ``curability" of experimental leukemia. \emph{Cancer Chemotherapy Reports}
  35:1--111.

\bibitem[{Skipper et~al.(1967)Skipper, Schabel~Jr, \protect\BIBand{}
  Wilcox}]{skipper1967XXI}
Skipper H, Schabel~Jr F, Wilcox W (1967) Experimental evaluation of potential
  anticancer agents. {XXI}. {S}cheduling of arabinosylcytosine to take
  advantage of its {S}-phase specificity against leukemia cells. \emph{Cancer
  Chemotherapy Reports} 51(3):125--165.

\bibitem[{S\"{u}li \protect\BIBand{} Mayers(2003)}]{Suli2003}
S\"{u}li E, Mayers D (2003) \emph{An Introduction to Numerical Analysis} (New
  York, United States of America: Cambridge University Press).

\bibitem[{Swan \protect\BIBand{} Vincent(1977)}]{Swan1977}
Swan G, Vincent T (1977) Optimal control analysis in the chemotherapy of {IgG}
  multiple myeloma. \emph{Bulletin of Mathematical Biology} 39(3):317--337.

\bibitem[{Swierniak et~al.(2009)Swierniak, Kimmel, \protect\BIBand{}
  Smieja}]{Swierniak2009}
Swierniak A, Kimmel M, Smieja J (2009) Mathematical modeling as a tool for
  planning anticancer therapy. \emph{European Journal of Pharmacology}
  625(1):108 -- 121.

\bibitem[{Tan et~al.(2002)Tan, Khor, Cai, Heng, \protect\BIBand{}
  Lee}]{Tan2002}
Tan K, Khor E, Cai J, Heng C, Lee T (2002) Automating the drug scheduling of
  cancer chemotherapy via evolutionary computation. \emph{Artificial
  Intelligence in Medicine} 25(2):169--185.

\bibitem[{Tj{\o}rve \protect\BIBand{} Tj{\o}rve(2017)}]{Tjorve2017}
Tj{\o}rve K, Tj{\o}rve E (2017) The use of {G}ompertz models in growth
  analyses, and new {G}ompertz-model approach: {A}n addition to the
  {U}nified-{R}ichards family. \emph{PLoS ONE} 12(6):1--17.

\bibitem[{Tse et~al.(2007)Tse, Liang, Leung, Lee, \protect\BIBand{}
  Mok}]{tse2007}
Tse S, Liang Y, Leung K, Lee K, Mok T (2007) A memetic algorithm for
  multiple-drug cancer chemotherapy schedule optimization. \emph{IEEE
  Transactions on Systems, Man, and Cybernetics, Part B (Cybernetics)}
  37(1):84--91.

\bibitem[{Tyagi \protect\BIBand{} Dey(2014)}]{tyagi2014}
Tyagi R, Dey P (2014) Needle tract seeding: An avoidable complication.
  \emph{Diagnostic Cytopathology} 42(7):636--640.

\bibitem[{Urquhart \protect\BIBand{} De~Klerk(1998)}]{urquhart1998contending}
Urquhart J, De~Klerk E (1998) Contending paradigms for the interpretation of
  data on patient compliance with therapeutic drug regimens. \emph{Statistics
  in Medicine} 17(3):251--267.

\bibitem[{Villasana \protect\BIBand{} Ochoa(2004)}]{villasana2004}
Villasana M, Ochoa G (2004) Heuristic design of cancer chemotherapies.
  \emph{IEEE Transactions on Evolutionary Computation} 8(6):513--521.

\bibitem[{{W}orld~{H}ealth {O}rganization(1979)}]{world1979handbook}
{W}orld~{H}ealth {O}rganization (1979) \emph{WHO handbook for reporting results
  of cancer treatment} (World Health Organization).

\bibitem[{Yuan et~al.(2015)Yuan, Di, Zhang et~al.}]{yuan2015}
Yuan P, Di L, Zhang X, et~al. (2015) Efficacy of oral etoposide in pretreated
  metastatic breast cancer: A multicenter {P}hase 2 study. \emph{Medicine}
  94(17).

\bibitem[{Zhao et~al.(2020)Zhao, Lian, Huo et~al.}]{zhao2020efficacy}
Zhao Q, Lian C, Huo Z, et~al. (2020) The efficacy and safety of neoadjuvant
  chemotherapy on patients with advanced gastric cancer: A multicenter
  randomized clinical trial. \emph{Cancer Medicine} 9(16):5731--5745.

\bibitem[{Zietz \protect\BIBand{} Nicolini(1979)}]{zietz1979}
Zietz S, Nicolini C (1979) Mathematical approaches to optimization of cancer
  chemotherapy. \emph{Bulletin of Mathematical Biology} 41(3):305--324.

\end{thebibliography}



\newpage
\setcounter{page}{1}
\RUNTITLE{Combination Chemotherapy Optimization with Discrete Dosing {\bf (e-companion)}}
%
%
%
\begin{APPENDICES}
\section{Additional Modeling Details}
\label{app:modeling}
%

This section provides details of the operational constraints concerning drug concentration, infusion rate, daily cumulative dose, pill administration, and rest days in our models, described in Section~\ref{sec:operational}. 

\subsection{Capecitabine}
\begin{subequations}\label{eq:capec}
\begin{align}
&\conc_{1, s} \leq \rhs_{1, \mathrm{conc}}, \forAll s \in \{0,\dots,S\}, \tag{Capec.Concentration.Max}\label{capec:conc}\\
&U_{1,s} \leq \rhs_{1,\mathrm{rate}}, \forAll s \in \{0,\dots,S\}, \tag{Capec.InfusionRate}\label{capec:rate}\\
&\sum\limits_{s \in D_m} \drug_{1, s} \leq \rhs_{1, \mathrm{cum}}, \forAll m \in \DAYS, \tag{Capec.Daily.Max}\label{capec:dayMax}\\
& \drug_{1, s} = 0, \forAll s \notin \MEALS, \tag{Capec.PillAdmin1}\label{capec:nopill}\\
&\drug_{1, s} = \para_{1, \mathrm{pill}}\,\bin_{1, \mathrm{pill}, s}, \forAll s \in \MEALS,  \tag{Capec.PillAdmin2}\label{capec:pill}\\
&\bin_{1,\mathrm{pill},\tym} \in \mathbb{Z}_{+}, \forAll \tym \in \MEALS. \tag{Capec.IntegerPills}\label{capec:intVars}
\end{align}
\end{subequations}
The first three constraints enforce the maximum concentration, infusion rate, and daily dose requirements. The constraints \eqref{capec:nopill}, \eqref{capec:pill}, and \eqref{capec:intVars} enforce that capecitabine doses are administered via pills (in discrete amounts) only in meal times.

\subsection{Docetaxel}
\begin{subequations}\label{eq:docet}
\begin{align}
&\conc_{2, s} \leq \rhs_{2, \mathrm{conc}}, \forAll s \in \{0,\dots,S\}, \tag{Docet.Concentration.Max}\label{docet:conc}\\
&\drug_{2, s} \leq \rhs_{2, \mathrm{rate}}, \forAll s \in \{0,\dots,\lastTym\}, \tag{Docet.InfusionRate} \label{docet:rate}\\
&\sum\limits_{s \in D_m} \drug_{2, s} \leq \rhs_{2, \mathrm{cum}}(1 - \bin^{m}_{2,\mathrm{rest}}), \forAll m \in \DAYS, \tag{Docet.Daily.Max}\label{docet:dayMax}\\
&\sum\limits_{l = 0}^{\min\{\para_{2,\mathrm{rest}}, M - m\}} (1 - \bin^{m+l}_{2,\mathrm{rest}}) \leq 1, \forAll m \in \DAYS, \tag{Docet.Rest}\label{docet:rest}\\
&\bin^{m}_{2, \mathrm{rest}} \in \mathbb{B}, \forAll m \in \DAYS. \tag{Doc.IntegerRest}\label{docet:intVars}
\end{align}
\end{subequations}
The constraint \eqref{docet:dayMax} enforces a capacity on the maximum administered dose per day depending on the rest mandate. The constraint \eqref{docet:rest} controls the selection of infusion sessions along with rest days. The other constraints are analogous to those of capecitabine.

\subsection{Etoposide}
\begin{subequations}\label{eq:etopo}
\begin{align}
&\conc_{3, s} \leq \rhs_{3, \mathrm{conc}}, \forAll s \in \{0,\dots,S\}, \tag{Etopo.Concentration.Max}\label{etopo:conc}\\
&\drug_{3,s} \leq \rhs_{3,\mathrm{rate}},\forAll s \in \{0,\dots,S\}, \tag{Etopo.InfusionRate}\label{etopo.rate}\\
&\sum\limits_{s \in D_m} \drug_{3, s} \leq \rhs_{3, \mathrm{cum}}, \forAll m \in \DAYS, \tag{Etopo.Daily.Max}\label{etopo:dayMax}\\
& \drug_{3, s} = 0, \forAll s \notin \MEALS, \tag{Etopo.PillAdmin1}\label{etopo:nopill}\\
&\drug_{3, s} = \para_{3, \mathrm{pill}} \, \bin_{3, \mathrm{pill}, s}, \forAll s \in \MEALS, \tag{Etopo.PillAdmin}\label{etopo:pill}\\
&\bin_{3,\mathrm{pill},\tym} \in \mathbb{Z}_{+}, \forAll \tym \in \{0,\dots,\lastTym\}. \tag{Etopo.IntegerPills}\label{etopo:intVars}
\end{align}
\end{subequations}
The constraints for etoposide are analogous to those of capecitabine.
\newpage
\section{Proofs}
\label{app:proofs}
\counterwithin{equation}{section}
\setcounter{equation}{0}
The forward Euler's method aims to approximate
\begin{equation} \label{prototypeODE} 
\dot{y}(t) = f(t,y),~y(0) = y_{0}.
\end{equation}

\begin{lemma}\thlabel{ExistenceUniqueness} \citep{Butcher2008}
Consider \eqref{prototypeODE} in which $f$ is continuous in its first variable and Lipschitz continuous in its second variable. Then,~\eqref{prototypeODE} has a unique solution.
\end{lemma}

\begin{repeattheorem}[\thref{uniqueCSolution}]
Suppose that, for a drug $d \in \dSet$, the administration function $\drug_{d}$ is continuous in time. Then, the differential equation
\begin{align} \label{eq:diffTHE1}
\dot{\conc_{d}}(t) = -\concPar_{d} \, \conc_{d}(t) + \drug_{d}(t)/\vol,~t \in [0, \finishTym],
\end{align}
governing the drug concentration function $\conc_{d}$, has a unique solution.
\end{repeattheorem}

\myProof
For a $d \in \dSet$, let $f_{\conc_{d}}(t, \conc_{d}(t)) = -\concPar_{d}\, \conc_{d}(t) + \drug_{d}(t)/\vol$. Observe that $f_{C_{d}}$ is linear in its second variable, hence Lipschitz continuous. Also, the continuity of $\drug_{d}(t)$ is sufficient to guarantee that $f_{\conc_{d}}$ is continuous in its first variable. Therefore, the differential equation~\eqref{eq:diffTHE1} has a unique solution by \thref{ExistenceUniqueness}.
$\myQED$

\begin{repeattheorem}[\thref{uniquePSolution}]
%
%
Suppose that the administration functions for all drugs, i.e., $\drug_{d},\forAll d \in \dSet$, are continuous in time. Then, for each cancer cell type $q \in \typeSet$, the differential equation
\begin{align} \label{eq:diffTHE2}
\dot{\pop_{q}}(t) = \gompOne\big(\pop_{q,\infty} - \pop_{q}(t)\big) - \sum\limits_{d \in \dSet} \dCons_{d,q}\exp(-\dExp_{d,q}\,t)\,\econ_{d}(t),~t \in [0, \finishTym],
\end{align}
governing the cell population function $\pop_{q}$, has a unique solution. 
\end{repeattheorem}

\myProof
%
%
By \thref{uniqueCSolution}, the functions $\conc_{d}(t),\forAll d \in \dSet$, are defined uniquely, and they are differentiable, hence continuous. This guarantees continuity of $\econ_{d}(t) = \max\{0,~\conc_{d}(t) - \rhs_{d,\mathrm{eff}}\},\forAll d \in \dSet,$ and implies the continuity of $f_{\pop_q}(t, \pop_q(t)) = \gompOne\big(\pop_{q,\infty} - \pop_{q}(t)\big) - \sum\limits_{d \in \dSet} \dCons_{d,q}\exp(-\dExp_{d,q}t)\,\econ_{d}(t)$ in its first variable. In addition, $f_{\pop_q}$ is linear, hence Lipschitz continuous, in its second variable. Therefore, by~\thref{ExistenceUniqueness}, the differential equation~\eqref{eq:diffTHE2} has a unique solution.
$\myQED$

\begin{repeattheorem}[\thref{thm:concIsStable}]
Let $\{\drug_{s}\}_{s \in \mathbb{Z}_{+}}$ be a bounded sequence, and $\concPar, h, \vol > 0$. Under the stability condition $h < \frac{2}{\concPar}$, the difference equation 
\begin{equation} \label{eq:diffTHE3}
\conc_{s+1} = \conc_{s} - h \, \concPar \, \conc_{s} + \drug_{s}/\vol,    
\end{equation}
is absolutely stable, for all $s \in \mathbb{Z}_{+}$.
\end{repeattheorem}

\myProof
We proceed by proving that $\conc_{s+1} = (1 - h\,\concPar)^{s+1}\conc_{0} + \frac{1}{\vol}\sum\limits_{k = 0}^{s}(1 - h\,\concPar)^{s - k}\drug_{k}$ using induction. The base case, $s = 0$, is immediate. Assume that, for all $s \leq s' \in \mathbb{Z}_{+}$, the claim holds. Then, 
\begin{align*}
\conc_{s'+1} &= \conc_{s'} - h\,\concPar\,\conc_{s'} + \drug_{s'}/\vol\\
&= (1 - h\,\concPar)\conc_{s'} + \drug_{s'}/\vol\\
&= (1 - h\,\concPar)\left((1 - h\,\concPar)^{s'}\conc_{0} + \frac{1}{\vol}\sum\limits_{k = 0}^{s'-1}(1 - h\,\concPar)^{s' - 1 - k}\drug_{k}\right )  + \drug_{s'}/\vol \\
&= (1 - h\,\concPar)^{s' + 1}\conc_{0} + \frac{1}{\vol}\sum\limits_{k = 0}^{s'}(1 - h\,\concPar)^{s' - k},
\end{align*}
where we use the induction hypothesis in the penultimate line. By induction, the recurrence relation holds for all $s \in \mathbb{Z}_{+}$. Because $\{\drug_{s}\}_{s \in \mathbb{Z}_{+}}$ is a bounded sequence, $\concPar, h$ are strictly positive; the difference equation~\eqref{eq:diffTHE3} is stable if the condition $h < \frac{2}{\concPar}$ is satisfied because it implies $|1 - h\, \concPar|  < 1$.
$\myQED$

\begin{repeattheorem}[\thref{thm:popIsStable}]
Let $\{F_{s}\}_{s \in \mathbb{Z}_{+}}$ be a bounded sequence and $\gompOne, h > 0$. Under the stability condition $h < \frac{2}{\gompOne}$, the difference equation 
\begin{equation} \label{eq:diffTHE4}
\pop_{s+1} = \pop_{s} + h\big(\gompOne \,  (\pop_{\infty} - \pop_{s}) - F_{s}\big),   
\end{equation}
is absolutely stable, for all $s \in \mathbb{Z}_{+}$.
\end{repeattheorem}

\myProof
We first prove the recurrence relation 
$\pop_{s+1} = (1 - \gompOne h)^{s + 1}\pop_{0} + \sum\limits_{k = 0}^{s}\big((1 - \gompOne h)^{s - k}\,h\,(\gompOne\pop_{\infty} - F_{k})\big)$ using induction. The base case, $s = 0$, is immediate. 
Assume that, for all $s \leq s' \in \mathbb{Z}_{+}$, the claim holds. Then, we have 
\begin{align*}
\pop_{s'+1} &= \pop_{s'} + h\,\big(\gompOne(\pop_{\infty} - \pop_{s'}) - F_{s'}\big)\\
&= (1 - \gompOne h)\pop_{s'} + h(\gompOne \pop_{\infty} - F_{s'})\\
&= (1 - \gompOne h)\left((1 - \gompOne h)^{s'}\pop_{0} + \sum\limits_{k = 0}^{s' - 1}(1 - \gompOne h)^{s' - 1- k}\,h\,(\gompOne \pop_{\infty} - F_{k})\right) + h\,(\gompOne \pop_{\infty} - F_{s'})\\
&= (1 - \gompOne h)^{s' + 1}\pop_{0} + \sum\limits_{k = 0}^{s'}(1 - \gompOne h)^{s' - k}\,h\,(\gompOne \pop_{\infty} - F_{k}),
\end{align*}
where we apply the induction hypothesis in the second-to-last line. Because $\{F_{s}\}_{s \in \mathbb{Z}_{+}}$ is a bounded sequence, $h, \gompOne$ are strictly positive, which implies the difference equation~\eqref{eq:diffTHE4} is stable, if the condition $h < \frac{2}{\gompOne}$ is satisfied.
$\myQED$

\begin{repeattheorem}[\thref{thm:moreEffectiveConc}]
Consider the stochastic model~\eqref{eq:stoMILP}, with $\gompOne h \leq 1$. Let $({\bf \econ}^{[1]}$, ${\bf \pop}^{[1]})$ and $({\bf \econ}^{[2]}$, ${\bf \pop}^{[2]})$ each be components of different feasible solutions. Suppose $\econ^{[1]}_{d,\tym} \geq \econ^{[2]}_{d,\tym}$, for all $d \in \dSet,\tym \in \{0,\dots,\lastTym\}$. Then $\pop_{q,\lastTym}^{[1],(k)} \leq \pop_{q,\lastTym}^{[2],(k)}$, for all $q \in \typeSet, k \in \{1,\dots,K\}$.
\end{repeattheorem}

\myProof
Observe that $\dCons_{d,q} \geq 0$, for each $d \in \dSet, q \in \typeSet$, as it is the fractional kill effect parameter, and additionally, for all $s \in \{0,\dots,\lastTym\},~\exp(-\dExp\, t(\tym)) \geq 0$. Moreover, $h > 0$ because it is a positive unit of time. 
Choose $q \in \typeSet, k \in \{0,\dots,K\}$, and suppose $\pop_{q,\tym}^{[1],(k)} \leq \pop_{q,\tym}^{[2],(k)}$ for some $\tym \in \{0,\dots,\lastTym - 1\}$. Then, 
%
\begin{align*}
\pop_{q,\tym + 1}^{[1],(k)} &= \pop_{q,\tym}^{[1],(k)} + h\,\left(\gompOne\Big(\pop_{q,\infty}^{(k)} - \pop_{q,\tym}^{[1],(k)}\Big) - \sum\limits_{d \in \dSet} \dCons_{d,q}\,\exp\big(-\dExp_{d,q}\,t(\tym)\big)\,\econ_{d,\tym}^{[1]}\right)\\
&\leq \pop_{q,\tym}^{[1],(k)} + h\,\left(\gompOne\Big(\pop_{q,\infty}^{(k)} - \pop_{q,\tym}^{[1],(k)}\Big) - \sum\limits_{d \in \dSet} \dCons_{d,q}\,\exp\big(-\dExp_{d,q}\,t(\tym)\big)\,\econ_{d,\tym}^{[2]}\right)\\
&= (1 - \gompOne h)\,\pop_{q,\tym}^{[1],(k)} + h\left(\gompOne\,\pop_{q,\infty}^{(k)} - \sum\limits_{d \in \dSet} \dCons_{d,q}\,\exp\big(-\dExp_{d,q}\,t(\tym)\big)\,\econ_{d,\tym}^{[2]}\right)\\
&\leq (1 - \gompOne h)\,\pop_{q,\tym}^{[2],(k)} + h\left(\gompOne\,\pop_{q,\infty}^{(k)} - \sum\limits_{d \in \dSet} \dCons_{d,q}\,\exp\big(-\dExp_{d,q}\,t(\tym)\big)\,\econ_{d,\tym}^{[2]}\right)\\
&= \pop_{q,\tym}^{[2],(k)} + h\,\left(\gompOne\Big(\pop_{q,\infty}^{(k)} - \pop_{q,\tym}^{[2],(k)}\Big) - \sum\limits_{d \in \dSet} \dCons_{d,q}\,\exp\big(-\dExp_{d,q}\,t(\tym)\big)\,\econ_{d,\tym}^{[2]}\right)\\
&= \pop_{q,\tym+1}^{[2],(k)}.
\end{align*}
Because $\pop_{q,0}^{[1],(k)} = \pop_{q,0}^{[2],(k)} = \ipop_{q}^{(k)}$, it follows by induction that $\pop_{q,\tym}^{[1],(k)} \leq \pop_{q,\tym}^{[2],(k)}$, $\forAll \tym \in \{0,\dots,\lastTym\}$, which also implies $\pop_{q,\lastTym}^{[1],(k)} \leq \pop_{q,\lastTym}^{[2],(k)}$.
$\myQED$

\begin{repeattheorem}[\thref{nestedEulerApprox}]
Consider the system of differential equations
\begin{subequations}
\begin{align*}
\dot{y}(t) &= f(t, y, z),~y(0) = y_{0},\\
\dot{z}(t) &= g(t, z),~z(0) = z_{0},
\end{align*}
\end{subequations}
and the Euler's approximation $\{(y_{s}, z_{s})\}_{s = 0}^{S}$ with step size $h$, given by
$y_{s+1} = y_{s} + h \, f\big(t(s), y_{s}, z_{s}\big)$ and $z_{s+1} = z_{s} + h \, g\big(t(s), z_{s}\big)$. 
Let $\lambda_{z} = \max\{|z_{s} - z_{0}|, s \in \{0,\dots,S\}\}$, and suppose $g$ is continuous in both variables and Lipschitz continuous in its second variable, i.e., there exists $L_{g} > 0$ such that for all $t \in [0, T]~\text{and}~u, v \in \mathbb{R}$ with $|u - z_{0}| \leq \lambda_{z}, |v - z_{0}| \leq \lambda_{z}$,
\begin{align*}
|g(t,u) - g(t,v)| \leq L_{g}|u - v|.
\end{align*} 
Similarly, suppose $f$ is continuous in all variables and Lipschitz continuous (with respect to the $\ell_{1}$ norm) in its second and third variables with constant $L_{f}$. Furthermore, suppose $y$ and $z$ are twice continuously differentiable. Then, for all $s \in \{0,\dots,S\}$,
\begin{align*}
|y_{s} - y(t(s))| \leq \frac{h}{2}\left(\frac{\alpha_{z}}{L_{g}}(e^{L_{g}T} - 1) + \frac{\alpha_{y}}{L_{f}}\right)(e^{L_{f}T} - 1),
\end{align*}
where $\alpha_{z} = \max\limits_{\tau \in [0,T]} |\ddot{z}(\tau)|$ and $\alpha_{y} = \max\limits_{\tau \in [0,T]} |\ddot{y}(\tau)|$. 
\end{repeattheorem}

\myProof
By the well-known single-stage Euler's method approximation analysis, 
\begin{align*}
|z_{s} - z(t(s))| \leq \frac{h\,\alpha_{z}}{2 L_{g}}(e^{L_{g}T} - 1),
\end{align*}
for all $s \in \{0,\dots,S\}$; see~\cite{Suli2003}. By Taylor's expansion, $y(t(s+1)) = y(t(s)) + h\, \dot{y}(t(s)) + \dfrac{h^{2}}{2} \ddot{y}(\xi)$, for some $\xi \in [t(s), t(s+1)]$. Thus, we have
\begin{align*}
\big|y_{s+1} - y(t(s+1))\big| &= \bigg\vert y_{s} - y(t(s)) + h\left(f(t(s), y_{s}, z_{s}) - \dot{y}(t(s)) - \frac{1}{2}h\,\ddot{y}(\xi)\right)\bigg\vert\\
&= \bigg\vert y_{s} - y(t(s)) + h\left(f(t(s), y_{s}, z_{s}) - f(t(s), y(t(s)), z(t(s))) - \frac{1}{2}h\,\ddot{y}(\xi)\right)\bigg\vert\\
&\leq \big|y_{s} - y(t(s))\big| + \bigg\vert h\Big(f(t(s), y_{s}, z_{s}) - f(t(s), y(t(s)), z(t(s)))\Big)\bigg\vert + \frac{h^{2}\alpha_{y}}{2}.
\end{align*}
Using the error bound for $z_{s}$ and the Lipschitzness of $f$, we have
\begin{align*}
\big|y_{s+1} - y(t(s+1))\big| &\leq \big|y_{s} - y(t(s))\big| + h\,L_{f}\big|\big|(y_{s},z_{s}) - (y(t(s)), z(t(s)))\big|\big|_{1} + \frac{h^{2}\alpha_{y}}{2}\\
&\leq \big|y_{s} - y(t(s))\big| + h\,L_{f}\big|y_{s} - y(t(s))\big| + h\,L_{f}\big|z_{s} - z(t(s))\big| + \frac{h^{2}\alpha_{y}}{2}\\
&\leq (1 + h\,L_{f})\big|y_{s} - y(t(s))\big| + \frac{h^{2}\,L_{f}\,\alpha_{z}}{2L_{g}}(e^{L_{g}T} - 1) + \frac{h^{2}\alpha_{y}}{2}.
\end{align*}
An induction argument on $s$ shows that
\begin{align*}
\big|y_{s+1} - y(t(s+1))\big| &\leq \frac{h^{2}}{2}\left(\frac{L_{f}\,\alpha_{z}}{L_{g}}(e^{L_{g}T} - 1) + \alpha_{y}\right)\sum\limits_{k = 0}^{s}(1 + hL_{f})^{k},
\end{align*}
which implies
\begin{align*}
\big|y_{s} - y(t(s))\big| &\leq \frac{h^{2}}{2}\left(\frac{L_{f}\,\alpha_{z}}{L_{g}}(e^{L_{g}T} - 1) + \alpha_{y}\right)\left(\frac{(1 + h\,L_{f})^{s} - 1}{h\,L_{f}}\right)\\
&= \frac{h}{2}\left(\frac{\alpha_{z}}{L_{g}}(e^{L_{g}T} - 1) + \frac{\alpha_{y}}{L_{f}}\right)\big((1 + hL_{f})^{s} - 1\big)\\
&\leq \frac{h}{2}\left(\frac{\alpha_{z}}{L_{g}}(e^{L_{g}T} - 1) + \frac{\alpha_{y}}{L_{f}}\right)\big(e^{hL_{f}s} - 1\big)\\
&\leq \frac{h}{2}\left(\frac{\alpha_{z}}{L_{g}}(e^{L_{g}T} - 1) + \frac{\alpha_{y}}{L_{f}}\right)\big(e^{L_{f}T} - 1\big)
\end{align*}
and completes the proof.
\myQED

\begin{repeattheorem}[\thref{nestedPopApprox}]
Let $\conc(t)$ and $\pop_{q}(t),\forAll q \in \typeSet$, be the state variable functions for drug concentration and cell population, respectively, in an optimal solution to the (single-drug) chemotherapy optimization problem~\eqref{exactForwardSingleStage} without the effective concentration and operational constraints. Furthermore, suppose that $\conc$ and $\pop_{q},\forAll q \in \typeSet$, are twice continuously differentiable, and let $\tilde{\conc}$ and $\tilde{\pop}_{q},\,\forAll q \in \typeSet$, be the corresponding Euler's approximations with time-step $h$. Then,
\begin{align} \label{eq:nestedP}
\bigg\vert \sum\limits_{q \in \typeSet} \tilde{\pop}_{q,S} - \sum\limits_{q \in \typeSet} \pop_{q}(T) \bigg\vert &\leq \sum\limits_{q \in \typeSet} \frac{h}{2}\left(\frac{\alpha_{\conc}}{|\concPar|}(e^{|\concPar|T} - 1) + \frac{\alpha_{q}}{\max\{|\dCons_{q}|, |\gompOne|\}}\right)(e^{\max\{|\dCons_{q}|, |\gompOne|\}T} - 1),
\end{align}
where $\alpha_{\conc} = \max\limits_{\tau \in [0,T]}|\ddot{\conc}(\tau)|$ and $\alpha_{q} = \max\limits_{\tau \in [0,T]} |\ddot{\pop}_{q}(\tau)|,\forAll q \in \typeSet$.
\end{repeattheorem}

\myProof
First note that, given a fixed (optimal) control function $\drug: [0,T] \mapsto \mathbb{R}$, the corresponding state variable functions $\conc$ and $\pop_q,\forAll q \in \typeSet$, are uniquely defined by Theorems~\ref{uniqueCSolution} and~\ref{uniquePSolution}. Observe that $g (t,C)=-\concPar \, \conc + \drug/\vol$ is continuous in both variables, and $L_{g} = |\concPar|$ is the Lipschitz constant with respect to the second variable. Moreover, $f_{q} (t,\pop_{q},\conc) = \gompOne\big(\pop_{q,\infty} - \pop_{q}\big) -  \dCons_{q}\exp(-\dExp_{q}\,t)\,\conc, \forAll q \in \typeSet$, is continuous in all variables, and $L_{f_{q}} = \max\{|\dCons_{q}|, |\gompOne|\}$ is the Lipschitz constant with respect to the second and third variables. In addition, $\conc$ and $\pop_{q}, \forAll q \in \typeSet$, are twice continuously differentiable. Then, by~\thref{nestedEulerApprox},
\begin{align*}
\big|\tilde{\pop}_{q,S} - \pop_{q}(T) \big| &\leq \frac{h}{2}\left(\frac{\alpha_{\conc}}{|\concPar|}(e^{|\concPar|T} - 1) + \frac{\alpha_{q}}{\max\{|\dCons_{q}|, |\gompOne|\}}\right)(e^{\max\{|\dCons_{q}|, |\gompOne|\}T} - 1),
\end{align*}
for each $q \in \typeSet$, which immediately implies that~\eqref{eq:nestedP} holds and completes the proof. 
\myQED
\newpage
\section{Estimating Initial Tumor Population via a Branching Process}
\label{app:branching}
\counterwithin{figure}{section}
\setcounter{figure}{0}
Cell mutation is commonly modeled by a branching process; see e.g.,~\citep{Kimmel2015}. Recall that we consider four cell types: (0) non-resistant, (1) capecitabine-resistant, (2) docetaxel-resistant, and (3) etoposide-resistant. Assume the tumor starts with a single, non-resistant cell at generation zero. In each subsequent generation, every tumor cell gives rise to two daughter cells. With probability $\para_{(0,q)}$, exactly one of the daughter cells mutates to a cell type $q \in \{1,2,3\}$. Note that this mutation only occurs from the non-resistant cell type to a single-drug resistant cell type. All other cell types give rise to exactly two identical daughter cells in each subsequent generation.

To generate scenarios describing tumor heterogeneity, we simulated this process and clustered similar replications into aggregate scenarios, following previous approaches in the stochastic programming literature~\citep{Gulpinar2004,Beraldi2014}. To this end, we generated 10,000 replications over 30 reproductive generations with $\para_{(0,q)} = 0.5\%,\forAll q \in \{1,2,3\}$. The end populations from the simulations were normalized according to the Studentized residual with respect to each cell type. Then, the vectors of normalized populations were clustered (via K-means) into ten aggregate scenarios, after constructing a sum of squared errors plot and observing only small reductions with more scenario clusters. Each scenario in Table~\ref{tab:simulatedScenarioGeneration} is a cluster centroid, and the corresponding probability is the size of the cluster divided by the number of trials.  Figure~\ref{fig:scenarioGeneration} illustrates the distribution of different cell types in the simulated outcomes and the sum of squared errors plot from K-means clustering. In this figure, $N \equiv$ non-resistant, $C \equiv$ capecitabine-resistant, $D \equiv$ docetaxel-resistant, and $E \equiv$ etoposide-resistant.
\begin{figure}[th]
\centering
\begin{subfloat}[Simulated initial cell populations (log)]{
\centering
\includegraphics[width = .50\textwidth]{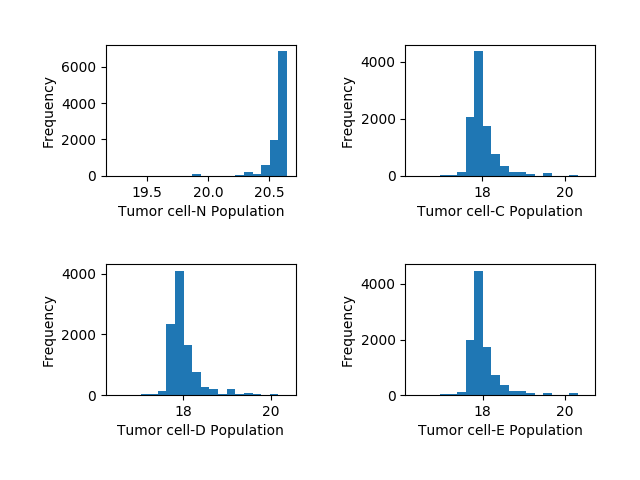}
\label{fig:histo}}
\end{subfloat}
\begin{subfloat}[Sum of squared errors from K-means clustering]{
\centering
\includegraphics[width = .45\textwidth]{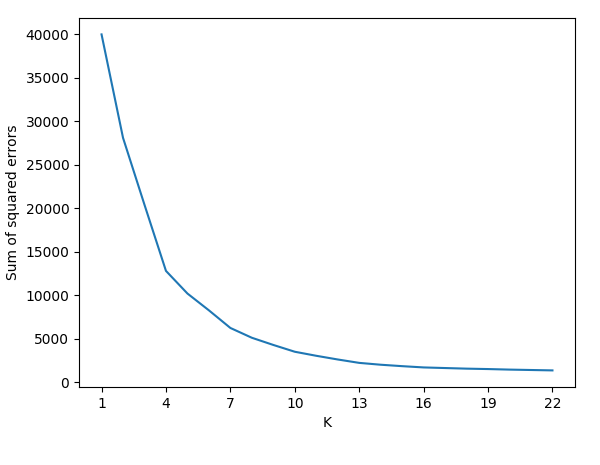}
\label{fig:Kmeans}}
\end{subfloat}
\caption{Scenario generation using a branching process and K-means clustering}
\label{fig:scenarioGeneration}
\end{figure}

Along this line, we present a general result concerning the expected cell type populations at an arbitrary generation. Let $\bipop(t) \in \mathbb{R}^{|\typeSet|}_{+}$ be a multivariate random variable representing the population of each tumor cell type at generation $t \in \mathbb{Z}_{+}$. Denote the vector of probabilities of outcomes for non-resistant cells by $\bpara_{(0)}$. Then, $\bipop(t+1) \sim \text{MD}(\ipop_{0}(t), \bpara_{(0)}) + 2\bipop(t) - \mathbf{e}_{0}\circ\bipop(t),$ where \text{MD}($n, \mathbf{p}$) denotes the multinomial distribution with $n$ trials and probabilities $\mathbf{p}$. We compute the expected cell populations of each type using the probability vector $\bpara_{(0)}$.
\begin{proposition}\thlabel{proposition:BranchingExpectation}
For $t \in \mathbb{Z}_{+}$, the expected cell type 0 (non-resistant) population at generation $t$ is $\mathbb{E}(\ipop_{0}(t)) = (\para_{(0,0)} + 1)^{t}$. The expected cell type $q$ (resistant) population at generation $t+1$ is $\mathbb{E}(\ipop_{q}(t+1)) = \sum\limits_{k = 0}^{t}2^{k}\para_{(0,q)}(\para_{(0,0)} + 1)^{t - k}$.
\end{proposition}

\myProof
Consider the non-resistant cell population. We prove the claim by induction, and one can observe that $t = 0$ and $t = 1$ readily serve as base cases. Assume for all $t \leq t_{0},~\mathbb{E}(\ipop_{0}(t)) = (\para_{(0,0)} + 1)^{t}$. Let $\mathbb{E}[\cdot]_{0}$ denote the expectation of the $0^{th}$ index. The expectation for generation $t_{0} + 1$ is
\begin{align*}
\mathbb{E}(\ipop_{0}(t_{0}+1)) &= \mathbb{E}[\text{MD}(\ipop_{0}(t_{0}),\bpara_{(0)}) + 2\bipop(t_{0}) - \mathbf{e}_{0}\circ\bipop(t_{0})]_{0}\\
&= \mathbb{E}[\text{MD}(\ipop_{0}(t_{0}),\bpara_{(0)})]_{0} + \mathbb{E}[2\bipop(t_{0}) - \mathbf{e}_{0}\circ\bipop(t_{0})]_{0}\\
&= \mathbb{E}[\text{MD}(\ipop_{0}(t_{0}),\bpara_{(0)})]_{0} + \mathbb{E}[\bipop(t_{0})]_{0}.
\end{align*}
Let $Z \sim [\text{MD}(\ipop_{0}(t_{0}), \bpara_{(0)})]_{0}$. Observe that $\bipop(t_{0})$ is a discrete random variable with finite support. Denote the outcomes of $\ipop_{0}(t_{0})$ by $\{\omega_{k}\}_{k = 1}^{K}$ with probabilities $p_{k} = \mathbf{Pr}\{\ipop_{0}(t_{0}) = \omega_{k}\},$ for all $k \in \{1,\dots,K\}$. Hence, by the law of total expectation 
\begin{align*}
\mathbb{E}[Z] &= \mathbb{E}[\mathbb{E}[Z \setbar \ipop(t_{0})]]\\
&= \sum\limits_{k = 1}^{K} p_{k}\mathbb{E}[\mathrm{MD}(\omega_{k}, \bpara_{(0)})]_{0}\\
&= \sum\limits_{k = 1}^{K}p_{k}\omega_{k}\para_{(0,0)}\\
&= \para_{(0,0)}\mathbb{E}[\ipop_{0}(t_{0})].
\end{align*}
Thus, by the induction hypothesis,
\begin{align*}
\mathbb{E}[\text{MD}(\ipop_{0}(t_{0}),\bpara_{(0)})]_{0} + \mathbb{E}[\bipop(t_{0})]_{0} &= \para_{(0,0)}\mathbb{E}[\ipop_{0}(t_{0})] + \mathbb{E}[\bipop(t_{0})]_{0}\\
&= \para_{(0,0)}(\para_{(0,0)} + 1)^{t_{0}} + (\para_{(0,0)} + 1)^{t_{0}}\\
&= (\para_{(0,0)} + 1)^{t_{0} + 1}.
\end{align*}
By induction, this proves the result for the non-resistant cell population, for all generation $t \in \mathbb{Z}_{+}$.

We now prove the resistant cell case ($q \neq 0$). For the base case of $t = 0$, we have $\mathbb{E}(\ipop_{q}(1)) = \para_{(0,q)}$, which is given by the definition of the multinomial random variable. In addition,
\begin{align*}
\sum\limits_{k = 0}^{0} 2^{k}\para_{(0,q)}(\para_{(0,0)} + 1)^{0 - k} = \para_{(0,q)},
\end{align*}
and the base case is satisfied. 
Next, assume for any $t < t_{0},$ we have $\mathbb{E}(\ipop_{q}(t+1)) = \sum\limits_{k = 0}^{t}2^{k}\para_{(0,q)}(\para_{(0,0)} + 1)^{t - k}$. We prove the case for $t_{0}$. Let $\mathbb{E}[\cdot]_{q}$ denote the expectation of the $q^{th}$ index, then the expectation at generation $t_{0} + 1$ is 
\begin{align*}
\mathbb{E}(\ipop_{q}(t_{0} + 1)) &= \mathbb{E}[\text{MD}(\ipop_{0}(t_{0}), \bpara_{(0)}) + 2\bipop(t_{0}) - \mathbf{e}_{0}\circ\bipop(t_{0})]_{q}\\
&= \mathbb{E}[\text{MD}(\ipop_{0}(t_{0}), \bpara_{(0)})]_{q} + \mathbb{E}[2\bipop(t_{0}) - \mathbf{e}_{0}\circ\bipop(t_{0})]_{q}.
\end{align*}
Similar to the non-resistant case, let $Z \sim [\text{MD}(\ipop_{0}(t_{0}), \bpara_{(0)})]_{q}$. By the law of total expectation, 
\begin{align*}
\mathbb{E}[Z] &= \mathbb{E}[\mathbb{E}[Z \setbar \ipop(t_{0})]]\\
&= \sum\limits_{k = 1}^{K} p_{k}\mathbb{E}[\mathrm{MD}(\omega_{k}, \bpara_{(0)})]_{q}\\
&= \sum\limits_{k = 1}^{K}p_{k}\omega_{k}\para_{(0,q)}\\
&= \para_{(0,q)}\mathbb{E}[\ipop_{0}(t_{0})].
\end{align*}
Thus, by the induction hypothesis, 
\begin{align*}
\mathbb{E}[\text{MD}(\ipop_{0}(t_{0}), \bpara_{(0)})]_{q} + \mathbb{E}[2\bipop(t_{0}) - \mathbf{e}_{0}\circ\bipop(t_{0})]_{q} 
= \ &\para_{(0,q)}\mathbb{E}[\ipop_{0}(t_{0})] + 2\mathbb{E}[\ipop_{q}(t_{0})]\\
= \ &\para_{(0,q)}(1 + \para_{(0,0)})^{t_{0}} + 2\sum\limits_{k = 0}^{t_{0} - 1}2^{k}\para_{(0,q)}(\para_{(0,0)} + 1)^{t_{0} - 1 - k}\\
= \ &\para_{(0,q)}(1 + \para_{(0,0)})^{t_{0}} + \sum\limits_{k = 0}^{t_{0} - 1}2^{k+1}\para_{(0,q)}(\para_{(0,0)} + 1)^{t_{0} - (k+1)} \\
= \ &\para_{(0,q)}(1 + \para_{(0,0)})^{t_{0}} + \sum\limits_{k = 1}^{t_{0}}2^{k}\para_{(0,q)}(\para_{(0,0)} + 1)^{t_{0} - k}\\
= \ &\sum\limits_{k = 0}^{t_{0}}2^{k}\para_{(0,q)}(\para_{(0,0)} + 1)^{t_{0} - k}.
\end{align*}
By induction, this proves the expectation for resistant cell type $q$, for all generations $t \in \mathbb{Z}_{+}$. 
$\myQED$
\newpage
\section{Model Parameters}
\label{app:estimate}
\counterwithin{table}{section}
\setcounter{table}{0}
Table~\ref{table:populationParameters} summarizes the cell population dynamics parameter values used in our numerical study, as described in Section~\ref{sec:calib}.
\begin{table}[ht]
\caption{Cell population dynamics parameters}
\small
\centering
\begin{tabular}{|c|c|c|c|}
\hline
Parameter Name & Symbol & Unit & Value \\
\hline
Initial cancer cell population & $\sum_{q \in \typeSet}\epop_{q, 0}$ & cell & $10^{9}$ \\
\hline
Cancer cell population limit & $\sum_{q \in \typeSet}\epop_{q,\infty}$ & cell & $10^{12}$ \\
\hline
Gompertz shape parameter & $\gompOne$ & $\text{day}^{-1}$ & $7 \bigcdot 10^{-4}$ \\
\hline
White blood cell initial population & $\epop_{w,0}$ & $\text{cell} \ \text{m}^{-3}$ & $8 \bigcdot 10^{12}$ \\
\hline
White blood cell turnover & $\nu_{w}$ & $\text{day}^{-1}$ & $0.15$ \\
\hline
White blood cell production rate & $\upsilon_{w}$ & $\text{cell} \ \text{m}^{-3} \ \text{day}^{-1}$ & $1.2 \bigcdot 10^{12}$ \\
\hline
\end{tabular}
\label{table:populationParameters}
\end{table}

\noindent For the pharmacokinetics parameters, i.e., elimination rate $\xi_d,\forAll d \in \dSet$, effect compartment $\vol$, and effectiveness threshold $\rhs_{d,\mathrm{eff}},\forAll d \in \dSet$, we used the values reported by~\cite{Iliadis2000} and~\cite{Frances2011}. Regarding the effect compartment $\vol$,~\cite{Iliadis2000} use a two-compartmental model for drug distribution, but we use a single-compartmental model assuming all dose goes through the first compartment into the second. To estimate the pharmacodynamics parameters, we used the following clinical administration regimens:
\begin{itemize}
\item Capecitabine~\citep{o2001randomized}: 
1255 $\text{mg}/\text{m}^2$ twice daily, 6 cycles of a two-week treatment period followed by a one-week rest period, response rate of 30\%,
\item Docetaxel~\citep{chan1999}: 
100 $\text{mg}/\text{m}^2$, 7 cycles of one-hour infusion every three weeks, response rate of 47\%,
\item Etoposide~\citep{yuan2015}: 60 $\text{mg}/\text{m}^2$ daily, 7 cycles of a 10-day treatment period followed by a 11-day rest period, response rate of 9\%.
\end{itemize}
Based on these regimens and the reported response rates, we estimated the fractional kill effect parameter of each drug for the non-resistant cell type, i.e., $\dCons_{d, 0}$, and set $\dCons_{d, q} = 0.25 \, \dCons_{d, 0},\forAll q \in \{1,2,3\}$, to account for drug-resistance in our numerical study. Although we use non-zero temporal resistance parameters $\dExp_{d,q},\forAll d \in \dSet,\, q \in \typeSet$, in our model, we set these parameters to zero for estimation of the fractional kill effect parameters, which is a conservative assumption. 


For a drug $d \in \dSet$, we simulated the corresponding administration regimen through $K$ trials; for each trial $k \in {\cK} = \{1,\dots,K\}$, we generated a kill parameter perturbation $\epsilon_{k}$ from a normal distribution with mean 0 and variance $\sigma^2$. The estimation of the fractional kill effect parameter $\dCons=\dCons_{d,0}$ is based on solving the following linear system with variables $\dCons$ and $\pop^{(k)}_s, \forAll k \in {\cK}, \, s \in \{0,\ldots,S'\}$:
\begin{equation}\label{estimateKillParams}
\begin{aligned}
&\pop^{(k)}_{\tym+1} = \pop^{(k)}_{\tym} + h\Big(\gompOne \big( \pop_{0,\infty} - \pop^{(k)}_{\tym}\big) - (\dCons + \epsilon_{k})\,\econ_{\tym}\Big), \forAll k \in {\cK}, \, s \in \{0,\ldots,S'\},\\
&\pop^{(k)}_{0} = \ln (N_{0}),\\
& \frac{1}{K}\sum\limits_{k = 1}^{K} \pop^{(k)}_{\lastTym'} = P' + \delta,
\end{aligned}
\end{equation}
where $\pop^{(k)}_s$ represents the tumor population at time-step $s \in \{0,\ldots,S'\}$ given the fractional kill parameter value $\dCons + \epsilon_{k}$, $\econ_{\tym}$ denotes drug effective concentration at time-step $s \in \{0,\ldots,S'\}$ dictated by the administration regimen, $P'$ is the target tumor size based on the definition of treatment response, i.e., cell population of a tumor (in logarithmic scale) with a diameter half the size of the original tumor, and $\delta$ is an adjustable parameter that we used to account for the reported response rate. Observe that, given perturbations $\{\epsilon_{k}\}_{k = 1}^{K}$, if the system~\eqref{estimateKillParams} has a solution, it is unique. By adjusting the value of $\delta$ and solving the system~\eqref{estimateKillParams} iteratively, for each drug $d \in \dSet$, we found the value of $\dCons=\dCons_{d, 0}$ that resulted in the response rate reported in the corresponding clinical trial. In these experiments, the number of trials was $K = 1,000$ and, for each drug, the standard deviation $\sigma$ was about 10\% of the final kill parameter value. Figure~\ref{fig:estimate} displays the simulation results. In this figure, $PRR$ denotes the treatment (partial) response rate. 
\begin{figure}[t]
\centering
\begin{subfloat}[Capecitabine]{
\centering
\includegraphics[width = .45\textwidth]{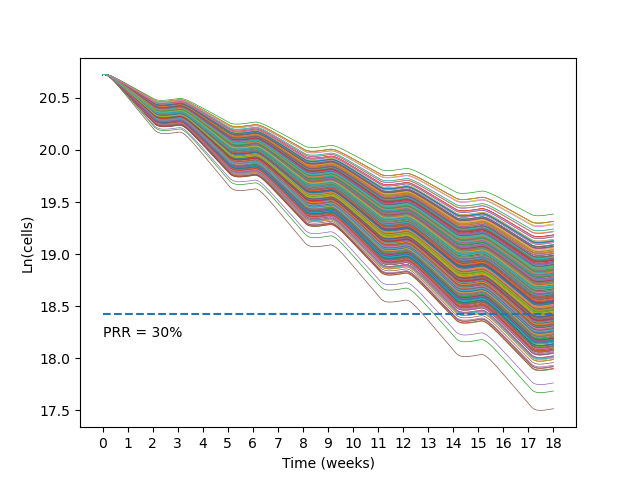}
\label{fig:est_C}}
\end{subfloat}
\begin{subfloat}[Docetaxel]{
\centering
\includegraphics[width = .45\textwidth]{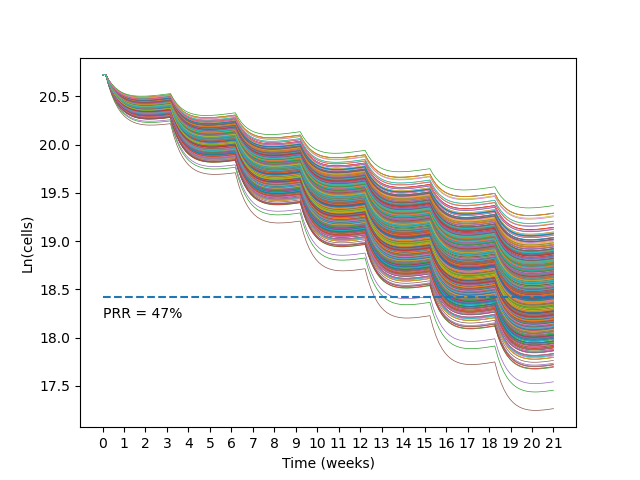}
\label{fig:est_D}}
\end{subfloat}
\begin{subfloat}[Etoposide]{
\centering
\includegraphics[width = .45\textwidth]{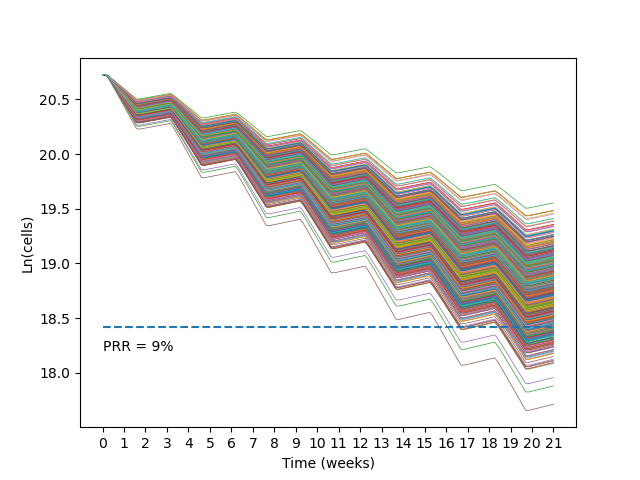}
\label{fig:est_E}}
\end{subfloat}
\caption{Estimation of the fractional kill effect parameters through simulation of clinical trials}
\label{fig:estimate}
\end{figure}

As mentioned earlier, we let $\dExp_{d,q} = 0,\forAll d \in \dSet,\, q \in \typeSet$, in estimating the fractional kill effect parameters. In the optimization model, we used the values reported by~\cite{Frances2011} for temporal resistance parameters of capecitabine and docetaxel. For etoposide, we could not find a reported value in the literature; thus, we made a conservative assumption that etoposide's parameter has a higher value than those of capecitabine and docetaxel. We also assumed a constant value for the temporal resistance parameters for each drug across the cancer cell types, i.e., for each $d\in \dSet$,~$\dExp_{d,0} = \dExp_{d,q},\forAll q \in \{1,2,3\}$. Note that, since the kill effect parameters are estimated assuming the absence of temporal resistance, any non-zero value for theses parameters in the optimization model generates a conservative solution. We acknowledge the need for future studies to best estimate the values of temporal resistance parameters, but we point out that the results of our sensitivity analysis show that the temporal resistance parameters are much less influential on the objective value of the optimization problems than the fractional kill effect parameters; see Section~\ref{sec:sensitivity}.

As stated in Section~\ref{sec:calib}, given the narrow therapeutic margin of cytotoxic drugs, we made a conservative assumption that the fractional kill effect of a drug on white blood cells is no less than its effect on non-resistant cancer cell types, i.e., $\dCons_{d,w} = \dCons_{d,0},\forAll d \in \dSet$. For the parameter representing the delay in the response of white blood cells to cytotoxic drugs, we used the value reported by~\cite{Iliadis2000}, i.e., $t_{w} = 5$ days. Table~\ref{table:PkPDParameters} displays the pharmacokinetics and pharmacodynamics parameter values used in our numerical study.
\begin{table}[ht]
\caption{Pharmacokinetics and pharmacodynamics parameters}
\small
\centering
\begin{tabular}{|c|c|c|c|}
\hline
Parameter Name & Symbol & Unit & Value (Capecitabine, Docetaxel, Etoposide) \\
\hline
Cancer kill effect & $\dCons_{d,0}$ & $\text{m}^{3} \ \text{gr}^{-1} \ \text{day}^{-1}$ & $(7.2\bigcdot 10^{-5},~8.0\bigcdot 10^{-3},~5.1\bigcdot 10^{-3} )$ \\
\hline 
White blood cell kill effect & $\dCons_{d,w}$ & $\text{m}^{3} \ \text{gr}^{-1} \ \text{day}^{-1}$ & $(7.2\bigcdot 10^{-5},~8.0\bigcdot 10^{-3},~5.1\bigcdot 10^{-3} )$ \\
\hline
Temporal resistance & $\dExp_{d,0}$ & $\text{day}^{-1}$ & $(5.7\bigcdot 10^{-3},~1.3\bigcdot 10^{-2},~1.4\bigcdot 10^{-2})$ \\
\hline
Elimination rate & $\xi_{d}$ & $\text{day}^{-1}$ & $(0.6,~0.2,~0.8)$ \\ 
\hline 
Effectiveness threshold & $\rhs_{d,\mathrm{eff}}$ & $\text{gr} \ \text{m}^{-3}$ & $(0.0,~0.0,~0.5)$ \\
\hline
Effect compartment volume & $\vol$ & $\text{m}^{3}$ & $15\bigcdot 10^{-3}$ \\
\hline
White blood cell delay & $t_{w}$ & day & 5\\
\hline
\end{tabular}
\label{table:PkPDParameters}
\end{table}

\phantom{empty line after table}

\noindent 
Finally, Table~\ref{table:oprParameters} summarizes the operational parameter values used in our numerical study. As stated in Section~\ref{sec:calib}, we used the simulation results of the abovementioned clinical trials to determine the operational parameters concerning maximum drug concentration, maximum infusion rate, and maximum cumulative daily dose. The values of neutropenia and lymphocytopenia thresholds are from the clinical literature~\citep{rosado2011hyper,mitrovic2012prognostic}.
\begin{table}[ht]
\caption{Operational parameters}
\small
\centering
\begin{tabular}{|c|c|c|c|}
\hline
Parameter Name & Symbol & Unit & Value (Capecitabine, Docetaxel, Etoposide) \\
\hline
Oral pill size & $\para_{d,\text{pill}}$ & \text{mg} & $(500,~\text{NA},~50)$ \\
\hline
Max drug concentration & $\rhs_{d, \text{conc}}$ & $\text{gr}/\vol$ & $(7.10,~0.17,~0.12)$ \\
\hline 
Max infusion rate (oral) & $\rhs_{d, \text{rate}}$ & $\text{gr m}^{-2}$ & $(1.25,~\text{NA},~0.03)$ \\
\hline 
Max infusion rate (intravenous)& $\rhs_{d, \text{rate}}$ & $\text{gr m}^{-2}\ \text{hr}^{-1}$ & $(\text{NA},~0.10,~\text{NA})$ \\
\hline
Max daily dose & $\rhs_{d, \text{cum}}$ & $\text{gr m}^{-2}$ & $(2.51,~0.10,~0.06)$ \\
\hline
Treatment rest & $\rhs_{d, \text{rest}}$ & day & $(\text{NA},~6,~\text{NA})$ \\
\hline
Treatment window & $\rhs_{d, \text{win}}$ & day & $(\text{NA},~1,~\text{NA})$ \\
\hline
Neutropenia threshold & $\rhs_{\text{neu}}$ & $\text{cell m}^{-3}$ & $2.5\bigcdot 10^{12}$ \\
\hline
Lymphocytopenia threshold & $\rhs_{\text{lym}}$ & $\text{cell m}^{-3}$ & $1\bigcdot 10^{12}$ \\
\hline
Neutrophil ratio & $\theta_{\text{neu}}$ & & $0.5$ \\
\hline
Lymphocyte ratio & $\theta_{\text{lym}}$ & & $0.3$ \\
\hline
\end{tabular}
\label{table:oprParameters}
\end{table}
\clearpage
\newpage
\section{Additional Figures}
\label{app:figs}
This section contains additional figures illustrating the results of our numerical study. 
Figures~\ref{fig:optAdmin}--\ref{fig:sto_remaining} provide further details on the optimal solutions of the proposed deterministic and stochastic models, i.e., formulations~\eqref{eq:detMILP} and \eqref{eq:stoMILP}. Figures~\ref{fig:prob_optAdmin}--\ref{fig:pro_remaining} concern the stochastic optimization model with a probability-based objective discussed in Section~\ref{sec:stochastic}. 

\newpage
\counterwithin{figure}{section}
\setcounter{figure}{0}
\begin{figure}[!ht]
\centering
\begin{subfloat}[Capecitabine]{
\centering
\includegraphics[width = .45\textwidth]{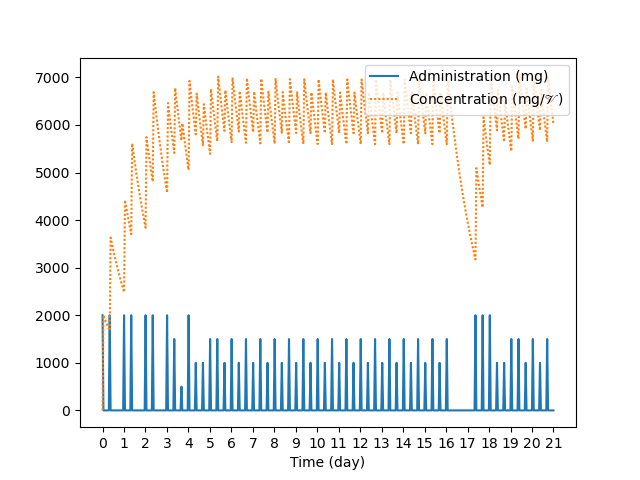}
\label{fig:opt_C}}
\end{subfloat}
\begin{subfloat}[Docetaxel]{
\centering
\includegraphics[width = .45\textwidth]{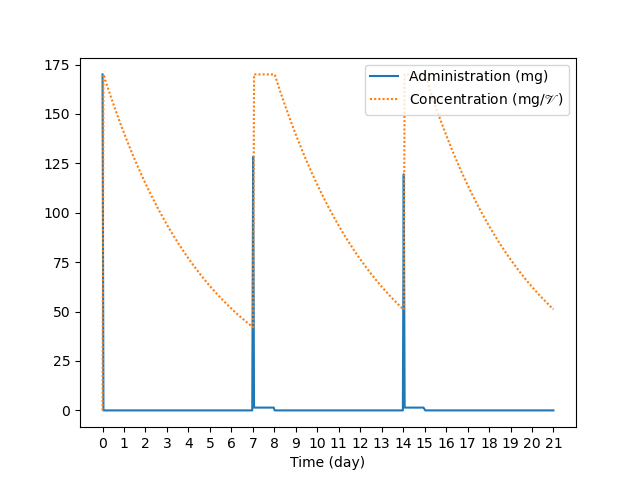}
\label{fig:opt_D}}
\end{subfloat}
\begin{subfloat}[Etoposide]{
\centering
\includegraphics[width = .45\textwidth]{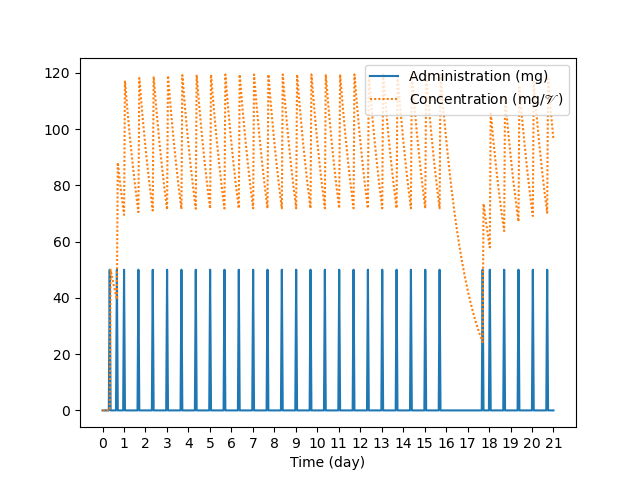}
\label{fig:opt_E}}
\end{subfloat}
\caption{Optimal administration and concentration for each drug, given by the deterministic model~\eqref{eq:detMILP}}
\label{fig:optAdmin}
\end{figure}
\begin{figure}[!ht]
\centering
\begin{subfloat}[Administration]{
\centering
\includegraphics[width = .45\textwidth]{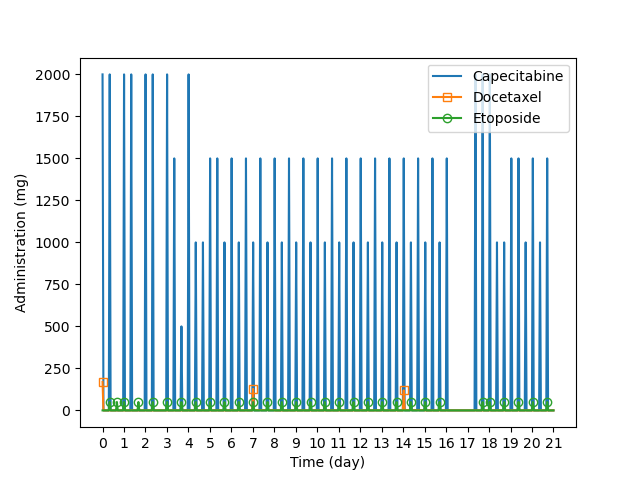}
\label{fig:sto_Administration}}
\end{subfloat}
\begin{subfloat}[White blood cell population]{
\centering
\includegraphics[width = .45\textwidth]{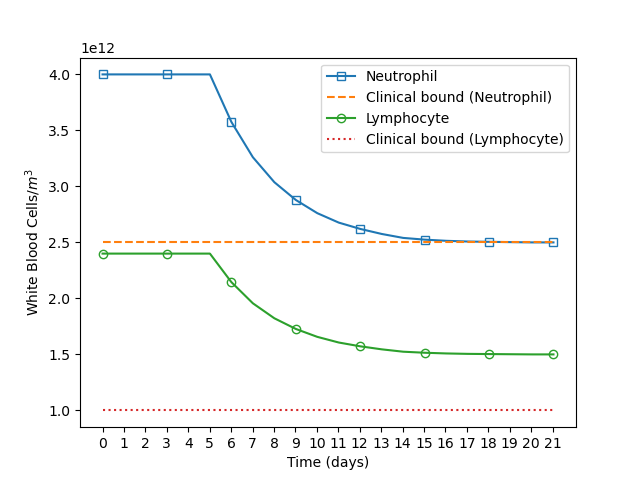}
\label{fig:sto_Concentration}}
\end{subfloat}
\caption{Optimal drug administration and white blood cell count, given by the chance-constrained model~\eqref{eq:stoMILP}}
\label{fig:sto_optAdmin}
\end{figure}
\begin{figure}
\centering
\begin{subfloat}{
\centering	
\includegraphics[width = .45\textwidth]{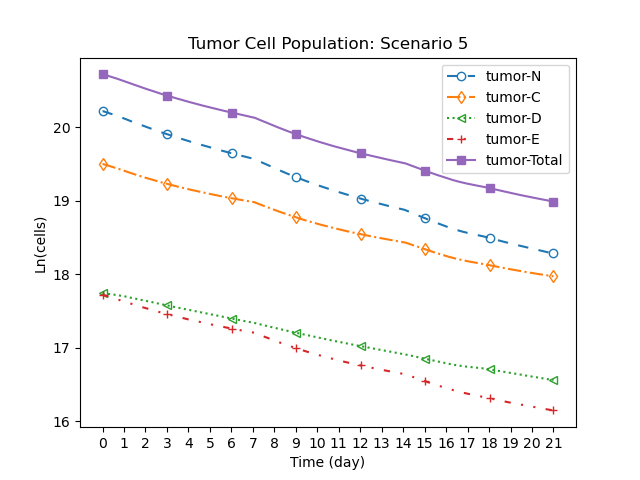}}
\end{subfloat}
\begin{subfloat}{
\centering
\includegraphics[width = .45\textwidth]{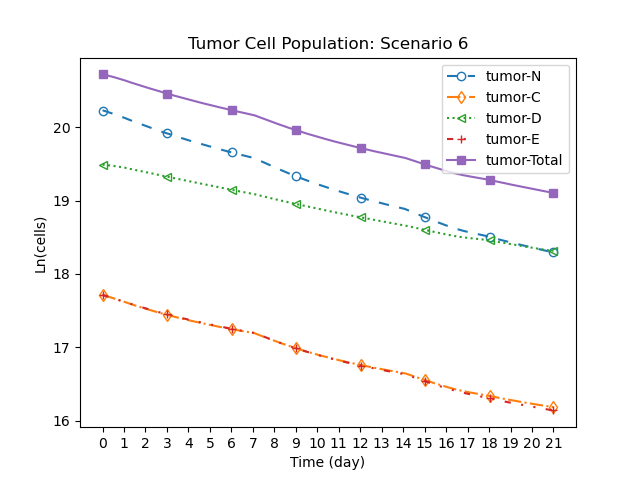}}
\end{subfloat}
\begin{subfloat}{
\centering	
\includegraphics[width = .45\textwidth]{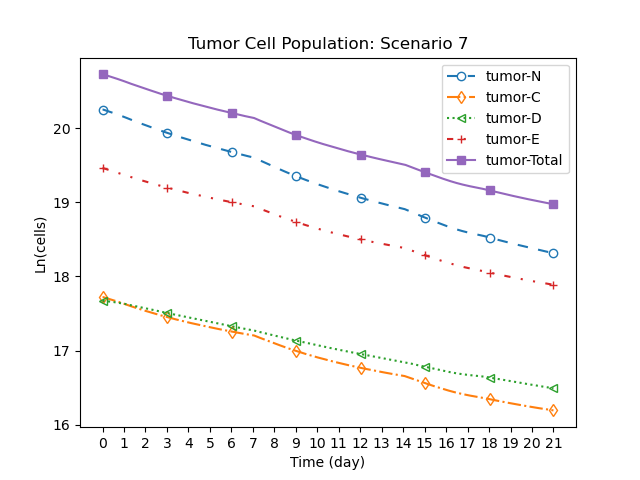}}
\end{subfloat}
\begin{subfloat}{
\centering
\includegraphics[width = .45\textwidth]{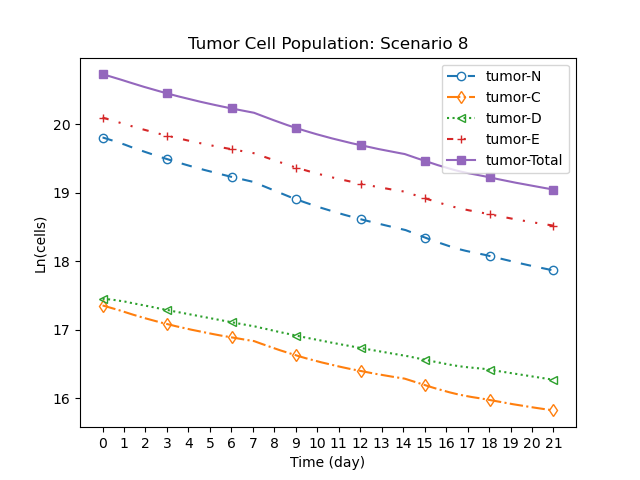}}
\end{subfloat}
\begin{subfloat}{
\centering	
\includegraphics[width = .45\textwidth]{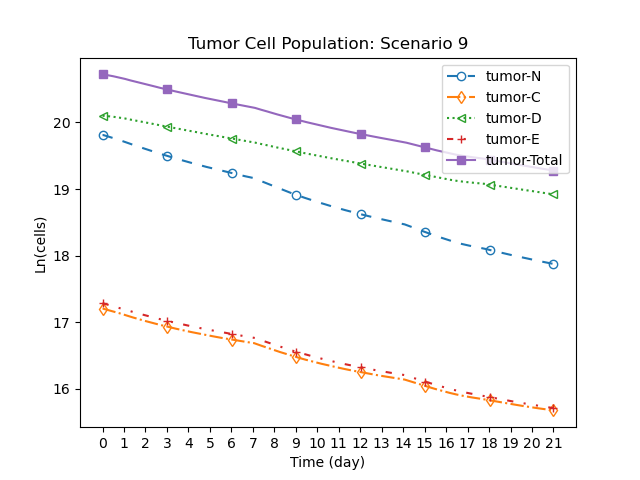}}
\end{subfloat}
\begin{subfloat}{
\centering
\includegraphics[width = .45\textwidth]{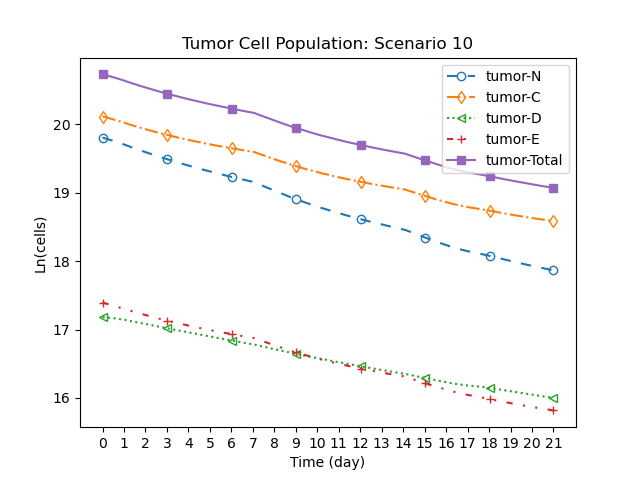}}
\end{subfloat}
\caption{Treatment effect on tumor cell populations under scenarios with the realization probability of 0.01 or less, given by the chance-constrained model~\eqref{eq:stoMILP}}
\label{fig:sto_remaining}
\end{figure}
\begin{figure}
\centering
\begin{subfloat}[Administration]{
\centering
\includegraphics[width = .45\textwidth]{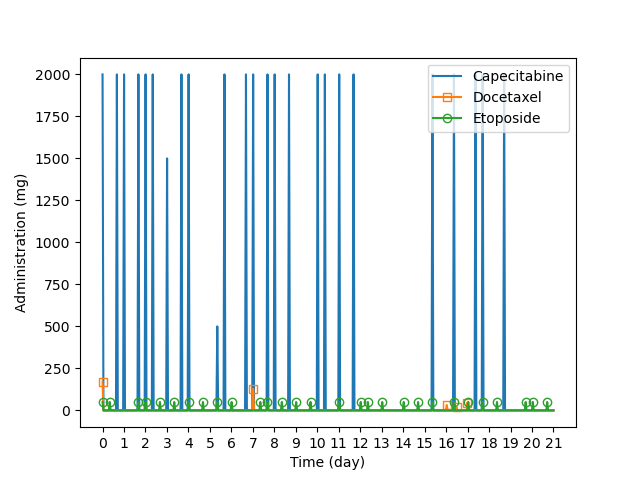}
\label{fig:prob_Administration}}
\end{subfloat}
\begin{subfloat}[White blood cell population]{
\centering
\includegraphics[width = .45\textwidth]{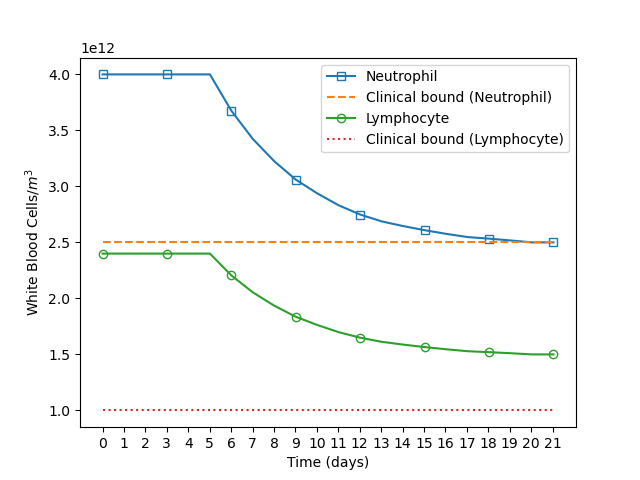}
\label{fig:prob_Concentration}}
\end{subfloat}
\caption{Optimal drug administration and white blood cell count, given by the chanced-constrained model with a probability-based objective (optimal objective value $\epsilon = 0$)}
\label{fig:prob_optAdmin}
\end{figure}
\begin{figure}
\centering
\begin{subfloat}{
\centering	
\includegraphics[width = .45\textwidth]{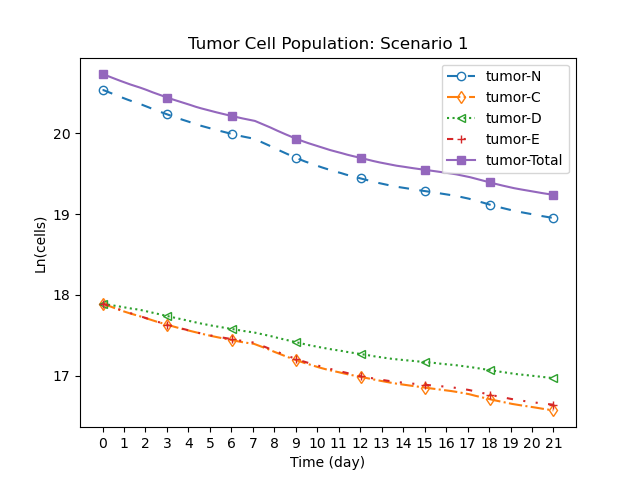}}
\end{subfloat}
\begin{subfloat}{
\centering
\includegraphics[width = .45\textwidth]{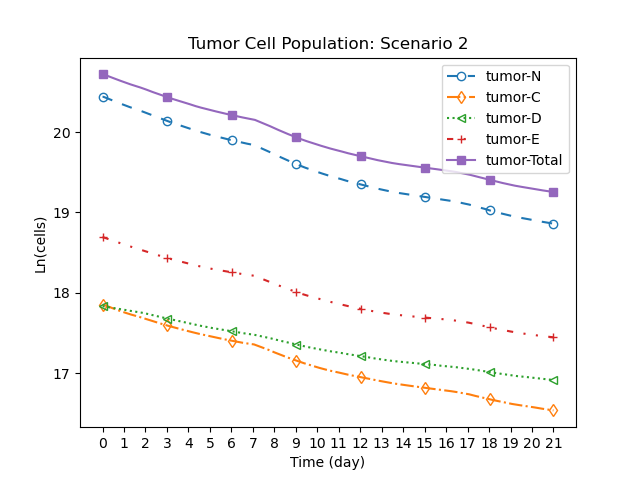}}
\end{subfloat}
\begin{subfloat}{
\centering	
\includegraphics[width = .45\textwidth]{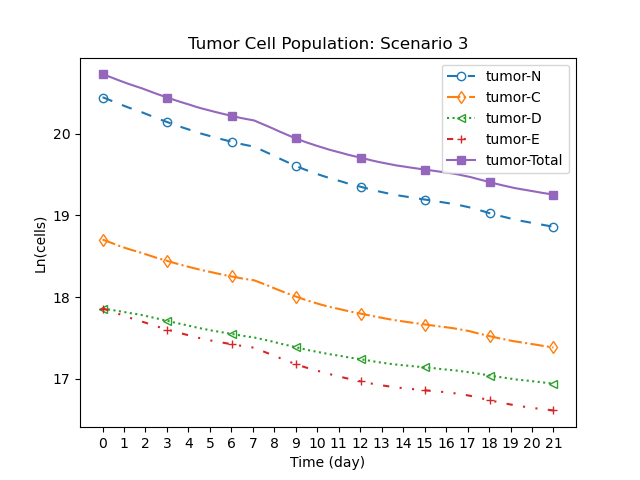}}
\end{subfloat}
\begin{subfloat}{
\centering
\includegraphics[width = .45\textwidth]{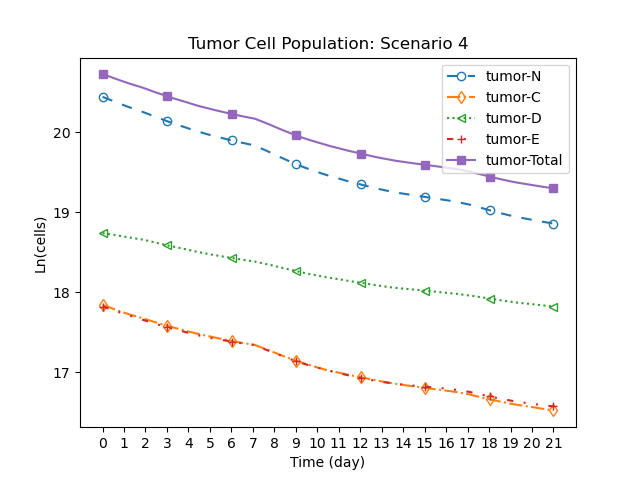}}
\end{subfloat}
\caption{Treatment effect on tumor cell populations under Scenarios 1--4, given by the chanced-constrained model with a probability-based objective (optimal objective value $\epsilon = 0$)}
\label{fig:pro_remaining}
\end{figure}
%
\end{APPENDICES}
\end{document}